\documentclass{imanum_arxiv}
\usepackage{amsfonts,amsmath,amsthm,amssymb}
\usepackage{stmaryrd}
\usepackage{mathtools}
\usepackage{url}
\usepackage{tcolorbox}
\usepackage{multicol}
\usepackage{tikz}
\usepackage{pgfplots}
\usepackage{dsfont}
\pgfplotsset{compat=1.7}
\definecolor{dkgreen}{rgb}{0,0.8,0.4}
\definecolor{mygreen}{HTML}{43a047}
\newcommand*{\llbrace}{\{\mskip-5mu\{}
\newcommand*{\rrbrace}{\}\mskip-5mu\}}
\DeclareMathOperator*{\esssup}{ess\,sup}

\newcommand{\vecc}[1]{\mathbf{#1}}

\newcommand{\dx}{\, \textup{d} x}

\renewcommand{\u}{\vecc{u}}
\newcommand{\uh}{\vecc{u}_h}
\renewcommand{\v}{\vecc{v}}

\newcommand{\e}{^{\textup{e}}}

\renewcommand{\t}{^{\textup{t}}}

\newcommand{\f}{^{\textup{f}}}

\newcommand{\R}{\mathbb{R}} % real
\newcommand{\C}{\mathbb{C}} % complex
\newcommand{\N}{\mathbb{N}} % natural
\newcommand{\Cinv}{C_{\textup{inv}}}
\newcounter{rownumber}

\newcommand{\GeD}{\Gamma^{\textup{e}}_{\textup{D}}}
\newcommand{\GaD}{\Gamma^{\textup{a}}_{\textup{D}}}
\newcommand{\GfD}{\Gamma^{\textup{f}}_{\textup{D}}}
\newcommand{\GtD}{\Gamma^{\textup{t}}_{\textup{D}}}

\newcommand{\GefI}{\Gamma^{\textup{e,f}}_{\textup{I}}}
\newcommand{\GeaI}{\Gamma^{\textup{e,a}}_{\textup{I}}}
\newcommand{\GaaI}{\Gamma^{\textup{a,a}}_{\textup{I}}}
\newcommand{\GftI}{\Gamma^{\textup{f,t}}_{\textup{I}}}
\newcommand{\Ha}{H^1_{\textup{D}}(\Omega_\textup{a})}

\newcommand{\Hf}{H^1_{\textup{D}}(\Omega_\textup{f})}
\newcommand{\He}{\boldsymbol{H}^1_{\textup{D}}(\Omega_\textup{e})}
\newcommand{\Ht}{{H}^1_{\textup{D}}(\Omega_\textup{t})}

\newcommand{\Fhaa}{\mathcal{F}_h^{\textup{a,a}}}
\newcommand{\producta}[2]{(#1,#2)_{\Omega_\textup{a}}} 
\newcommand{\producte}[2]{(#1,#2)_{\Omega_\textup{e}}} 
\newcommand{\productf}[2]{(#1,#2)_{\Omega_\textup{f}}} 
\newcommand{\productt}[2]{(#1,#2)_{\Omega_\textup{t}}} 

\newcommand{\productFa}[2]{\langle #1, \; #2 \rangle_{\mathcal{F}_h^{\textup{a,a}}}}
\newcommand{\productFea}[2]{\langle #1, \; #2 \rangle_{\mathcal{F}_h^{\textup{e,a}}}}
\newcommand{\productFf}[2]{\langle #1, \; #2 \rangle_{\mathcal{F}_h^{\textup{a,a}}}}

\newcommand{\productintea}[2]{(#1,#2)_{\Gamma_\textup{I}^{\textup{e}, \textup{a}}}} 
 
\newcommand{\norma}[1]{\|#1\|_{\Omega_\textup{a}}} 
\newcommand{\normE}[1]{\|#1\|_{\textup{E}}} 
\newcommand{\normf}[1]{\|#1\|_{\Omega_\textup{f}}} 
\newcommand{\normt}[1]{\|#1\|_{\Omega_\textup{t}}} 
 
\newcommand{\normEe}[1]{\|#1\|_{\textup{E}^\textup{e}}} 
\newcommand{\normLinfEe}[1]{\|#1\|_{L^{\infty}\textup{E}^\textup{e}}} 
\newcommand{\normEa}[1]{\|#1\|_{\textup{E}^\textup{a}}} 
\newcommand{\normLinfEa}[1]{\|#1\|_{L^{\infty}\textup{E}^\textup{a}}} 
\newcommand{\normEf}[1]{\|#1\|_{\textup{E}^\textup{f}}}

\newcommand{\normLinfE}[1]{\|#1\|_{L^{\infty}\textup{E}}} 
\newcommand{\norme}[1]{\|#1\|_{\Omega_\textup{e}}}
\newcommand{\normFa}[1]{\|\sqrt{\chi}\, \llbracket #1 \rrbracket \,\|_{\mathcal{F}_h^{\textup{a,a}}}}

\newcommand{\normsa}[1]{{\|#1\|}_{s, \textup{a}}}
\newcommand{\tripplenorme}[1]{{\vert\kern-0.25ex\vert\kern-0.25ex\vert #1 
		\vert\kern-0.25ex\vert\kern-0.25ex\vert}_{s, \textup{e}}}
\newcommand{\tripplenorma}[1]{{\vert\kern-0.25ex\vert\kern-0.25ex\vert #1 
		\vert\kern-0.25ex\vert\kern-0.25ex\vert}_{s, \textup{a}}}
\newcommand{\LtLt}[1]{\|#1\|_{L^2L^2}}	
\newcommand{\LinfLinf}[1]{\|#1\|_{L^{\infty}L^{\infty}}}
\newcommand{\LtLinf}[1]{\|#1\|_{L^2L^{\infty}}}	
\newcommand{\LinfLt}[1]{\|#1\|_{L^{\infty}L^2}}	

\newcommand{\forma}{a^{\textup{a}}}
\newcommand{\forme}{a^{\textup{e}}}

\newcommand{\formah}{a_h^{\textup{a}}}
\newcommand{\formeh}{a_h^{\textup{e}}}

\newcommand{\Ief}{\mathcal{I}}

\newcommand{\Iae}{\mathcal{I}}
\newcommand{\Iea}{\mathcal{I}}

\newcommand{\Veh}{\boldsymbol{V}_h^{\textup{e}} }
\newcommand{\Vth}{V_h^{\textup{t}} }
\newcommand{\Vah}{V_h^{\textup{a}} }
\newcommand{\Vfh}{V_h^{\textup{f}} }
\newcommand{\varrhoa}{\varrho^{\textup{a}}}
\newcommand{\varrhoe}{\varrho^{\textup{e}}}
\newcommand{\varrhof}{\varrho^{\textup{f}}}
\newcommand{\varrhot}{\varrho^{\textup{t}}}
\newcommand{\zetae}{\zeta^{\textup{e}}}
\newcommand{\zetat}{\zeta^{\textup{t}}}
\newcommand{\Omegaa}{\Omega_{\textup{a}}}
\newcommand{\Omegae}{\Omega_{\textup{e}}}
\newcommand{\Omegaf}{\Omega_{\textup{f}}}
\newcommand{\Omegat}{\Omega_{\textup{t}}}
\newcommand{\fe}{\boldsymbol{f}^{\textup{e}}}

\newcommand{\fbt}{\boldsymbol{f}^{\textup{t}}}
\newcommand{\ft}{{f}^{\textup{t}}}
\newcommand{\fth}{f^{\textup{t}}_h}
\newcommand{\fa}{f^{\textup{a}}} 
\newcommand{\ff}{f^{\textup{f}}} 
\newcommand{\ffh}{f^{\textup{f}}_h}

\newcommand{\fah}{f_h^{\textup{a}}} 
\newcommand{\normala}{\vecc{n}^{\textup{a}}} 
\newcommand{\normale}{\vecc{n}^{\textup{e}}} 
\newcommand{\normalf}{\vecc{n}^{\textup{f}}} 
\newcommand{\normalt}{\vecc{n}^{\textup{t}}} 
\usepackage{graphicx}

\jno{drnxxx}
\received{\today}
\revised{xxxxx}
\newcommand{\change}[2]{#2}
\newcommand{\version}[2]{#1}  % Arxiv version

\begin{document}

%\title{A hybrid discontinuous Galerkin method for nonlinear elasto-acoustic coupling}
\title{\change{blue}{A discontinuous Galerkin coupling for nonlinear elasto-acoustics}}
%\title{A discontinuous Galerkin type coupling for nonlinear elasto-acoustics}
%\title{A discontinuous Galerkin coupling method  for nonlinear elasto-acoustics}
%\title{A discontinuous Galerkin type coupling method  for nonlinear elasto-acoustics}
% Short title for running heads:
\shorttitle{A discontinuous Galerkin coupling for nonlinear elasto-acoustics}

\author{%
{\sc
Markus Muhr\thanks{Corresponding author. Email: muhr@ma.tum.de},
Barbara Wohlmuth\thanks{Email: wohlmuth@ma.tum.de}} \\[2pt]
Department of Mathematics, Technical University of Munich, Boltzmannstraße 3, 85748 Garching bei München, Germany\\[6pt]
{\sc and}\\[6pt]
{\sc Vanja Nikoli\a{'}c}\thanks{Email: vanja.nikolic@ru.nl}\\[2pt]
IMAPP – Department of Mathematics, Radboud University, Heyendaalseweg 135,
6525 AJ Nijmegen, The Netherlands
}
% Short list of authors for running heads:
\shortauthorlist{M. Muhr, V. Nikoli\a{'}c and B. Wohlmuth}

\maketitle

\begin{abstract}
% Body of abstract:
{Inspired by medical applications of high-intensity ultrasound, we study a coupled elasto-acoustic problem with general acoustic nonlinearities of quadratic type as they arise, for example, in the Westervelt and Kuznetsov equations of nonlinear acoustics. We derive convergence rates in the energy norm of a finite element approximation to the coupled problem in a setting that involves different acoustic materials and hence jumps within material parameters. A subdomain-based discontinuous Galerkin approach realizes the acoustic-acoustic coupling of different materials. At the same time, elasto-acoustic interface conditions are used for a mutual exchange of forces between the different models. Numerical simulations back up the theoretical findings in a three-dimensional setting with academic test cases as well as in an application-oriented simulation, where the modeling of human tissue as an elastic versus an acoustic medium is compared.}
% Keywords:
{nonlinear acoustics; elasto-acoustic coupling; discontinuous Galerkin methods; Westervelt's equation; Kuznetsov's equation; ultrasonic waves.}
\end{abstract}

\section{Introduction}\label{sec:Introduction}
This work is devoted to the numerical analysis and simulation of coupled linear elastic-nonlinear acoustic problems, which arise in a variety of medical and industrial applications of ultrasonic waves.\\
\indent Coupled problems in general play an essential role in different fields of application. Starting with fluid-structure interaction (FSI), where fluid-dynamical equations like Navier-Stokes equations are coupled to equations of solid mechanics \citep{FSIReview,wick2020optimization}, over electro/magneto-mechanical systems involving electromagnetic field equations \citep{MKaltenbacher}, to biomedical applications, such as the mathematical modeling of tumor growth and the simulation of thermo-acoustic effects \citep{shevchenko2012multi}, which couple wave and heat equations. All such problems have in common at least two different models - mostly in the form of partial differential equations - that describe, for example one a fluid part, the other a solid part of the overall problem. It is then necessary to couple the individual models to a global system. Therefore different techniques can be used, for example by individual (volumetric) source terms or factors mutually depending on each other such as temperature and speed of sound in \citep{shevchenko2012multi} or -- in case of spatially separated model-domains with some common interface -- via a coupling using Lagrange multipliers, as discussed in \citep{Li1998} for a simple model-problem, or in \citep{flemisch2006elasto} for an elasto-acoustic coupling. A direct exchange via boundary conditions can also be used, as discussed in \citep{felippa2011classification} in a FSI setting.\\
\indent The coupled problem considered in this work is the elasto-acoustic problem, where elastic subdomains are modeled by an elastodynamic wave equation, whereas on (mostly) fluid subdomains a nonlinear acoustic wave equation is employed to model the propagation of pressure waves. The choice of the subdomain models often depends on the concrete application of interest. Therefore, different coupling of elastic and acoustic equations have been investigated in the literature. %with the simplest case of linear elasticity and a linear, undamped wave equation being the most common. 
A Lagrange-multiplier approach for the coupling of a linear elastic equation and a linear undamped acoustic problem is considered in \citep{flemisch2006elasto}, whereas the approach of \citep{antonietti2018high} employs a direct exchange of forces via interface conditions.\\
\indent In this work, displacement-based linear elasticity is chosen to model the elastic parts of the domain, whereas a nonlinear wave equation with a general nonlinearity of quadratic type is used for the acoustic part:
\begin{align} \label{nlwave_eq}
\frac{1}{c^2}\ddot{\psi}- \Delta \psi -\frac{b}{c^2} \Delta \dot{\psi} = \frac{1}{c^2}(k_1 \dot{\psi} \ddot{\psi}+ k_2 \nabla \psi \cdot \nabla \dot{\psi}).
\end{align}
This equation is formulated in terms of the acoustic velocity potential $\psi= \psi(x,t)$, from which the acoustic pressure can be computed using the relation $p=\varrho\dot{\psi}$, where $\varrho$ is the mass density of the medium. In homogeneous media, depending on the choice of the parameters $k_1$ and $k_2$, the well-known Westervelt and Kuznetsov equations of nonlinear acoustics are obtained from \eqref{nlwave_eq} as special cases; see~\citep{kuznetsov1971equation, westervelt1963parametric, MKaltenbacher} and Section~\ref{Sec:ProblemSetting} below for more details on the modeling aspects. Such wave equations model the propagation of sound waves with sufficiently high amplitudes and frequencies through thermoviscous media. The thermoviscous dissipation is reflected in the presence of the strong damping $-\frac{b}{c^2}\Delta \dot{\psi}$ in the equation. The magnitude of the parameter $b$ plays a significant role in the analysis. If $b \rightarrow 0^+$, the hyperbolic character of $\frac{1}{c^2}\ddot{\psi}-\Delta\psi$ dominates the behavior of solutions. On the other hand, if the sound diffusivity $b$ is relatively large, the parabolic character of $\frac{1}{c^2}\ddot{\psi}-\frac{b}{c^2}\Delta\dot{\psi}$ is pronounced and, in homogeneous media, one can even expect exponential decay of the energy of solutions in time; see, for example,~\citep{mizohata1993global, kaltenbacher2009global}. Thus, the presence of this strong damping will be crucial in our energy arguments.\\ %which contributes to its parabolic-like character and will thus be crucial in our energy analysis.  \\ % In contrast to linear wave equations, the nonlinearities influence the steepening of the wave front and, when not much dissipation is present in the system, development of a shock; see, for example,~\citep{duck2002nonlinear,MuhrSO, MKaltenbacher} for a deeper insight into this process. \\
\indent A common medical use of high-intensity ultrasonic waves is in the non-invasive treatments of kidney stones \citep{lee2007ultrasound} and certain types of cancer \citep{kennedy2005high}. In such scenarios, high-intensity waves are generated, for example, by vibrating piezoelectric transducers \citep{MKaltenbacher} arranged on a part of the boundary of an acoustic medium, which might, for example, be a simple pipe filled with water. Due to the shape of the transducers, the ultrasonic waves induced into the acoustic medium are focused towards the central axis of the device, increasing the pressure amplitude there even more \citep{kaltenbacher2016shape,MuhrSO}. Finally, on the other end of the acoustic channel, the ultrasound waves propagate into human tissue, where they further travel towards, for example, a kidney stone or tumor. Due to the fact that the waves are focused and of high power, the waves reaching, e.g., the kidney stone have enough energy to break it apart into smaller debris, thus avoiding an open surgery for the patient; see, e.g.,~\citep{skolarikos2006extracorporeal}. \\ % Such treatments are referred to as \textit{extracorporeal shock wave lithotripsy} (ESWL)~\citep{skolarikos2006extracorporeal}.\\
\indent For solving the arising partial differential equations, finite elements are our method of choice. A full discontinuous Galerkin higher-order approach has been developed in \citep{antonietti2018high} for the linear elasto-acoustic problem. In \citep{antonietti2020high}, the acoustic propagation in homogeneous media without the quadratic gradient nonlinearity has been treated in a discontinuous Galerkin setting. \change{blue}{In \cite{antonietti2012non}, the purely elastic problem has been considered using the dG approach for the coupling of different subdomains; the developed approached was termed  the Discontinuous Galerkin Spectral Element Method (DGSEM)}. \\%, which we also adopt for this work.\\ %The elastic problem alone has been extensively studied as well in similar settings; see, for example, in \citep{AntoniettiAyusoMazzieriQuarteroni_2016,antonietti2017hp}.\\
\indent  To our best knowledge, rigorous numerical analysis of an interface coupling between linear elasticity and nonlinear acoustics has not been performed before. The available results in the literature on the numerical analysis seem to focus either on related linear coupled problems or on nonlinear (acoustic) wave propagation in homogeneous media. In particular, a priori analysis of a high-order discontinuous Galerkin method for a spatial discretization of the corresponding undamped linear problem ($b=k_1=k_2=0$) has been conducted in~\citep{antonietti2018high}. Error analysis of the (semi-)discrete Kuznetsov equation is seemingly still an open problem. A high-order discontinuous Galerkin method for the Westervelt equation (obtained when $k_2=0$ in  \eqref{nlwave_eq}) in homogeneous media has been analyzed in~\citep{antonietti2020high}. We also point out the results of~\citep{MaierDiss}, which as a particular case include rigorous analysis of the semi-discrete (based on a finite-element discretization) and fully discrete Westervelt equation in pressure form in homogeneous inviscid media, where $b=0$.\\
\indent The hybrid approach developed in this work combines advantages of both the conforming and discontinuous Galerkin framework. Within regions of constant material properties and simple geometry we use a conforming hexahedral mesh with its - conforming as well - nodal degrees of freedom located at the Gauss--Laguerre--Lobatto points known from spectral FEM. However, in order for more flexibility concerning the meshes being used, complex, possibly non-conforming, interfaces are resolved using a DG-paradigm, which also allows us to clearly separate blocks of different material properties, resolve jumps in the coefficients precisely and refine (material) subdomains of the overall model individually to their needs.\\
\indent The rest of the paper is structured as follows. In Section~\ref{Sec:ProblemSetting}, we discuss the modeling aspects of the problem, interface conditions, and the variational formulation as well as the discrete finite element setting. Section \ref{sec:SemiDiscProb} introduces the semi-discrete problem with the necessary notation for the hybrid DG-coupling approach. In Section \ref{sec:LinearErrorEstimate}, we prove stability of a linear version of the approximate problem and derive convergence rates of this numerical scheme in the energy norm, where the proof of the error-estimate is postponed to the appendix. Section~\ref{sec:NonlinearArgument} then presents a fixed-point argument, through which we transfer the results from the linear error estimator to the nonlinear case, under the assumption of sufficiently small data and the global mesh size. In particular, the main error estimate for the nonlinear problem is obtained in Theorem~\ref{Thm:Main} below.  Finally, Section~\ref{sec:NumericalSimulation} contains our 3D numerical examples that illustrate the convergence results numerically in certain academic test cases as well as in more application-oriented simulations.

\section{The nonlinear elasto-acoustic problem} \label{Sec:ProblemSetting}
\begin{figure}[h]
	\hspace*{-1cm}\hspace*{1cm}\includegraphics[scale=0.12]{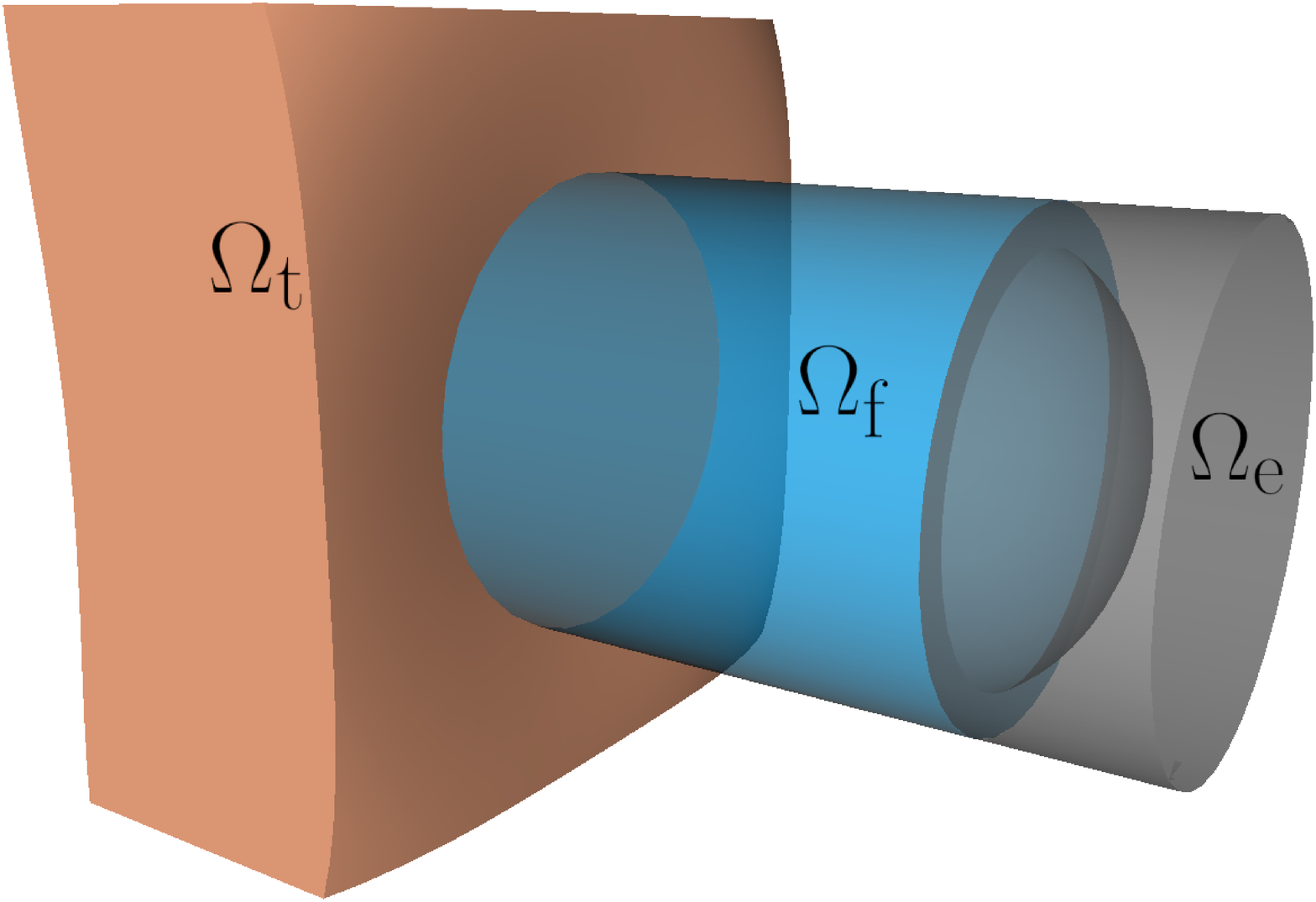}\hspace*{0.5cm}\raisebox{-3mm}{\includegraphics[scale=0.145]{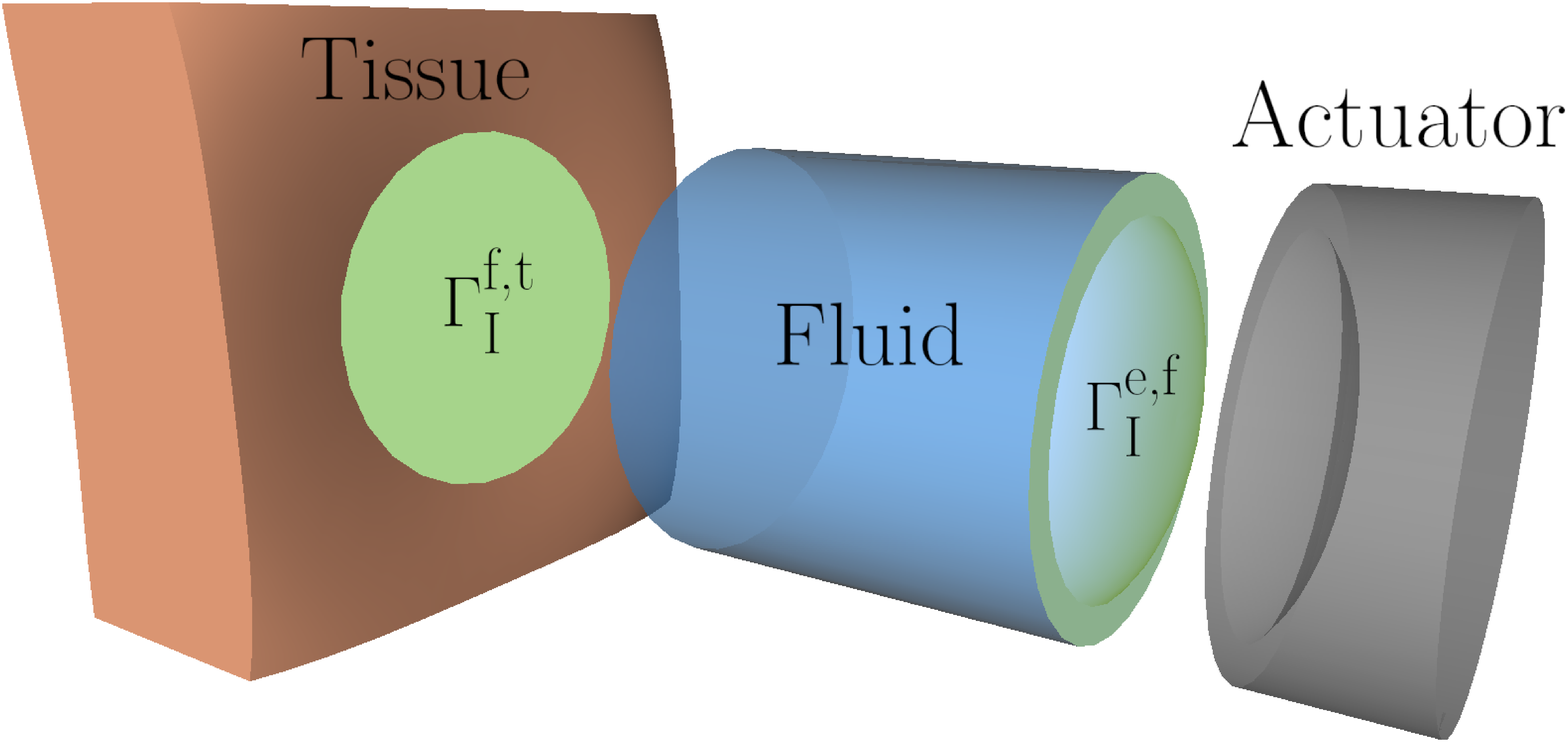}}\\
	\hspace*{1cm}\includegraphics[scale=0.12]{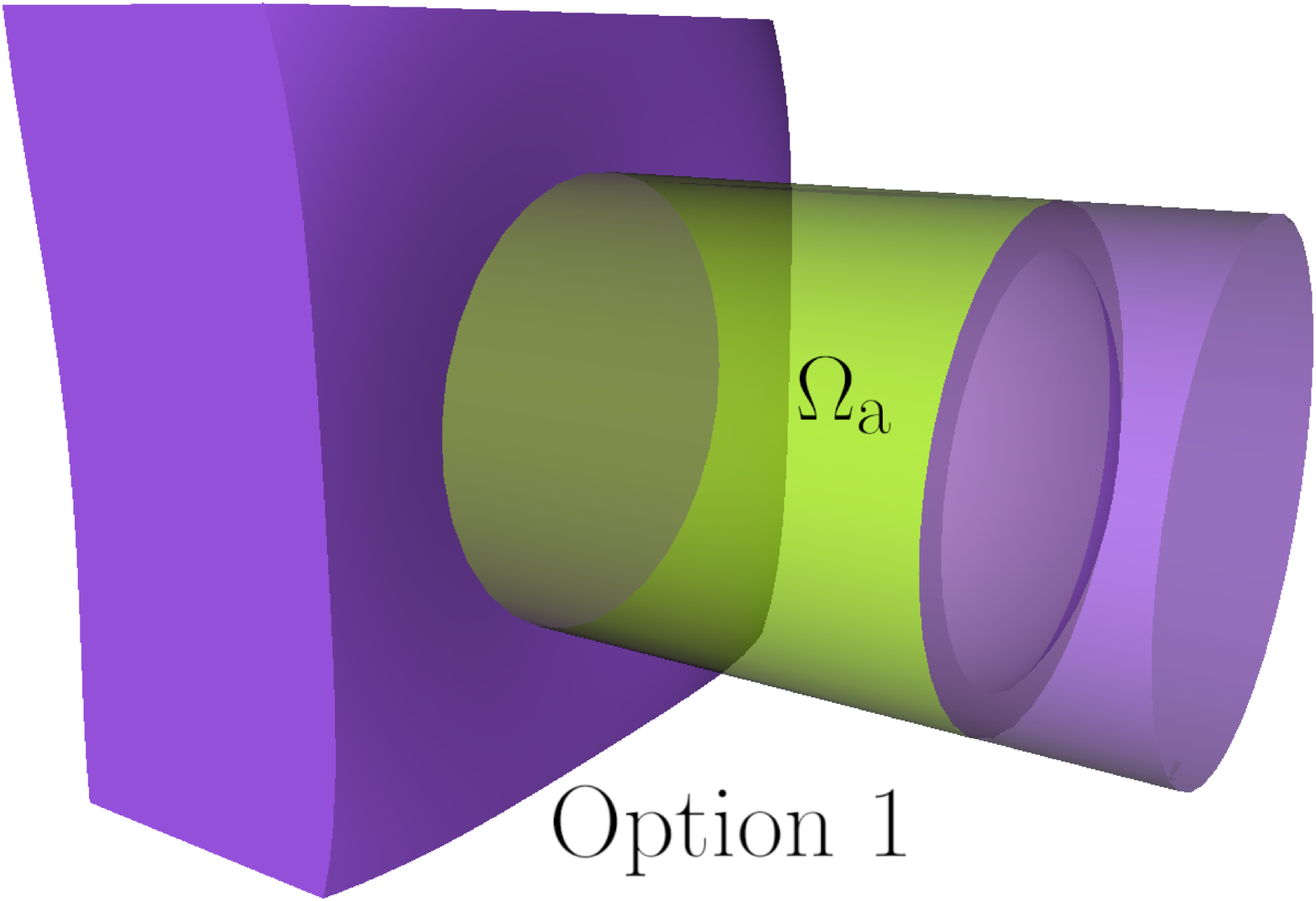}\hspace*{0.5cm}\includegraphics[scale=0.13]{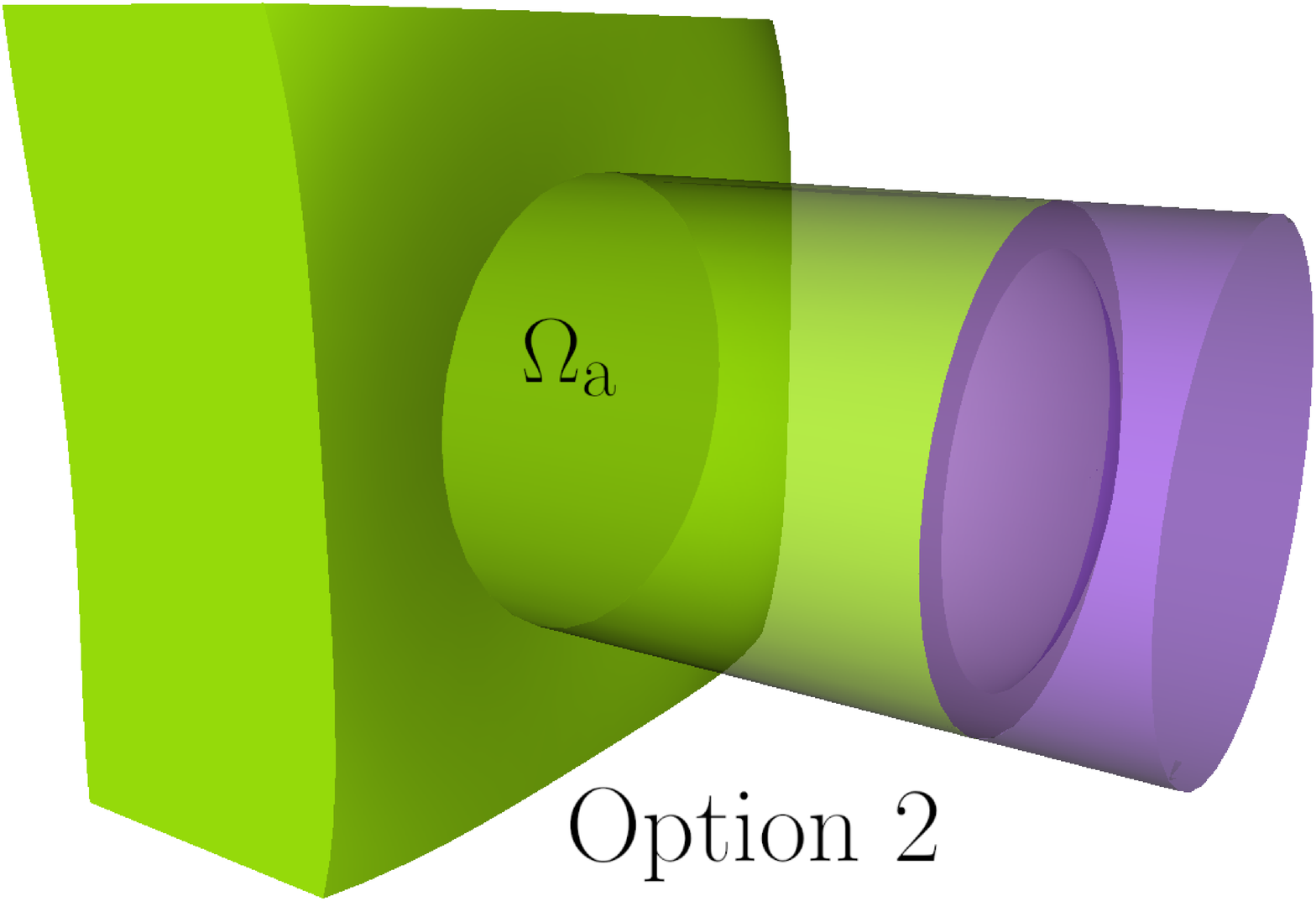}
	\caption{\label{fig:Domain}\textbf{(top left)} Actuator $\Omega_{\textup{e}}$, fluid $\Omega_{\textup{f}}$, and tissue $\Omega_{\textup{t}}$ domain in exemplary situation motivated by focused ultrasound applications. \textbf{(top right)} Material interfaces $\Gamma_{\textup{I}}^{\textup{e,f}}$, $\Gamma_{\textup{I}}^{\textup{f,t}}$ (highlighted in green). \textbf{(bottom left)} Option 1 with elastic tissue domain. The fluid forms the acoustic domain $\Omega_{\textup{a}}$ alone (highlighted in green). \textbf{(bottom right)} Option 2 with acoustic tissue domain. The acoustic domain is formed by $\Omega_{\textup{f}}$ and $\Omega_{\textup{t}}$ (highlighted in green) with an acoustic-acoustic interface $\Gamma_{\textup{I}}^{\textup{a,a}}=\Gamma_{\textup{I}}^{\textup{f,t}}$. % with material parameter jump in between. 
		Here the elasto-acoustic interface $\Gamma_{\textup{I}}^{\textup{e,a}}=\Gamma_{\textup{I}}^{\textup{e,f}}$.}
\end{figure}

We begin by stating the elasto-acoustic coupled problem in its strong form. The different domains of elastic as well as acoustic media are denoted using indices motivated by medical ultrasound applications as discussed in the introduction. Nevertheless, whenever the term ``excitator'' or ``actuator'' is used one might also think of some general elastic domain $\Omegae$\change{orange}{$\subset\R^d$} being coupled to some again  general ``fluid'' domain $\Omegaf$\change{orange}{$\subset\R^d$} representing an acoustic medium. What is here called ``tissue'' domain $\Omegat$\change{orange}{$\subset\R^d$} will be a domain where both an elastic or an acoustic model could be used, depending on the concrete application. We refer to Figure~\ref{fig:Domain} below for a graphical representation of the individual domains and interfaces.\\
\indent Throughout the paper, we assume the elastic and acoustic domains to be polygonal and convex, so that they can be discretized exactly. Furthermore, $T>0$ is a given fixed final time.\\
\indent The coupling terms as well as the analysis of the global problem depend on the choice of the material models. For the theoretical considerations, we restrict ourselves to the case of an acoustic model in the tissue domain. From the point of view of analysis, this is a more challenging problem since it involves a different type of interface coupling. However, our numerical experiments will deal with both situations and compare the two models in an applicational context. 
\begin{alignat*}{2}  
&\hspace*{-1cm} \textup{\bf Actuator / Mechanical excitator} \textup{ (elastic medium)} \hspace*{0cm}  &&\\[1mm]
&\varrhoe \, \ddot{\u}\e {+2\varrhoe \zetae\, \dot{\u}\e + \varrhoe(\zetae)^2\,\u\e} - \textup{div}(\,{\vecc{\boldsymbol{\sigma}}}(\u\e)\,) = \fe, && \hspace{0.5cm} \textup{in } \Omegae \times (0,T],\\
&\u\e =\, \vecc{0},   &&\hspace{0.5cm}\textup{on } \GeD\times (0,T],\\
&(\u\e, \dot{\u}\e)\vert_{t=0} =(\u^{\textup{e}}_0, \u^{\textup{e}}_1), && \hspace{0.5cm} \textup{in } \Omegae,\\[2mm]
&\hspace*{-1cm}\textup{\bf Fluid} \textup{ (acoustic medium)} \hspace*{5cm}  &&\\[1mm]
&(c\f)^{-2}\ddot{\psi}\f- \Delta \psi\f - \dfrac{b\f}{(c\f)^2}\, \Delta \dot{\psi}\f  =\ff(\dot{\psi}\f, \ddot{\psi}\f, \nabla \psi\f, \nabla \dot{\psi}\f), && \hspace{0.5cm} \textup{in } \Omegaf\times(0,T],\\
&\psi\f =0, && \hspace{0.5cm} \textup{on }\GfD\times (0,T],\\
&(\psi\f, \dot{\psi}\f)_{t=0} =(\psi\f_0, \psi\f_1), && \hspace{0.5cm} \textup{in } \Omegaf,\\[2mm]
&\hspace*{-1cm}\textup{\bf Tissue, Option 1}  \textup{ (elastic medium)} \hspace*{3cm}  &&\\[1mm]
&\varrhot \, \ddot{\u}\t +2\varrhot \zetat\, \dot{\u}\t + \varrhot(\zetat)^2\,\u\t - \textup{div}(\,{\vecc{\boldsymbol{\sigma}}}(\u\t)\,) = \fbt, && \hspace{0.5cm} \textup{in } \Omegat \times (0,T],\\
&\u\t = 0, && \hspace{0.5cm} \textup{on } \GtD\times (0,T],\\
&(\u\t, \dot{\u}\t)\vert_{t=0} =(\u^{\textup{t}}_0,\u^{\textup{t}}_1), && \hspace{0.5cm} \textup{in } \Omegat,\\[2mm]
&\hspace*{-1cm}\textup{\bf Tissue, Option 2} \textup{ (acoustic medium)} \hspace*{2.75cm}  &&\\[1mm]
&(c\t)^{-2}\ddot{\psi}\t- \Delta \psi\t - \dfrac{b\t}{(c\t)^2}\, \Delta \dot{\psi}\t  =\ft(\dot{\psi}\t, \ddot{\psi}\t, \nabla \psi\t, \nabla \dot{\psi}\t),&& \hspace{0.5cm} \textup{in } \Omegat\times(0,T],\\
&\psi\t  =0, && \hspace{0.5cm} \textup{on }\GtD\times (0,T],\\
&(\psi\t, \dot{\psi}\t)_{t=0} =(\psi\t_0, \psi\t_1), && \hspace{0.5cm} \textup{in } \Omegat.
\end{alignat*}
The acoustic right-hand side nonlinearity is defined as
\begin{equation}\label{eq:AcRHS}
\begin{aligned}
f^{\textup{i}}(\dot{\psi}^{\textup{i}}, \ddot{\psi}^{\textup{i}}, \nabla \psi^{\textup{i}}, \nabla \dot{\psi}^{\textup{i}})=&\, \dfrac{1}{(c^{\textup{i}})^2}\,\left( k^{\textup{i}}_1(\dot{\psi}^{\textup{i}})^2+k^{\textup{i}}_2|\nabla \psi^{\textup{i}}|^2\right)_t = \dfrac{2}{(c^{\textup{i}})^2}\,\left( k_1^{\textup{i}}\dot{\psi}^{\textup{i}} \ddot{\psi}^{\textup{i}}+k_2^{\textup{i}}\nabla \psi^{\textup{i}} \cdot \nabla \dot{\psi}^{\textup{i}}\right),
\end{aligned}
\end{equation}
where $\textup{i} \in \{\textup{f}, \textup{t}\}$. Here the notation $(\cdot)_t$ stands for the partial derivative with respect to time. The stress-tensor $\boldsymbol{\sigma}$ in the elastic domains is given by Hook's law in the framework of linear elasticity as \[\boldsymbol{\sigma}(\u) = \mathbb{C}\boldsymbol{\varepsilon}(\u)=\lambda\textup{div}(\u)\mathds{1}+2\mu\boldsymbol{\varepsilon}(\u),\] where $\vecc{u}$ stands for either $\vecc{u}\e$ or $\vecc{u}\t$. Furthermore, $\boldsymbol{\varepsilon}(\u)=\frac{1}{2}(\nabla\u+\nabla\u^{\top})$, $\mathbb{C}$ is the material stiffness tensor with $\lambda$ and $\mu$ being the Lam\a{'}e-parameters of the given material. \change{orange}{Further, $\zeta\e$ and $\zeta\t$ are damping parameters}. The source terms $\fe$ and $\fbt$ denote external body forces.\\

\paragraph{\bf Interface coupling.} We next discuss the coupling conditions over different interfaces. We introduce \[\GefI=\overline{\Omegae}\cap\overline{\Omegaf}\]
as the interface between the purely elastic and purely acoustic domains. Moreover, we set \[\GftI=\overline{\Omegaf}\cap\overline{\Omegat}\]
as the interface between the purely acoustic domain $\Omegaf$ and the domain $\Omegat$ for which the \emph{role} of the interface will differ, but not its geometry. \change{orange}{Finally we denote by $\vecc{n}^{\textup{i}},\textup{i}\in\lbrace\textup{e,t,f}\rbrace$ the outward unit normal vector field of the respective subdomain $\Omega_{\textup{i}}$.}\\

\noindent
\emph{Actuator-Fluid interface.} This is an elasto-acoustic interface that prescribes a normal stress/pressure on the elastic domain given by the acoustic pressure in the fluid domain, while the acoustic particle velocity is prescribed by the displacement velocity of the elastic medium
\begin{alignat*}{2}
\vecc{\boldsymbol{\sigma}}(\u\e)~\normale &= -\varrhof\,\left(\dot{\psi}\f+\frac{b\f}{(c\f)^2}\ddot{\psi}\f\right)~\normale, && \hspace{0.5cm} \textup{on } \GefI\times (0,T],\\
\frac{\partial\psi\f}{\partial\normalf} {+ \frac{b\f}{(c\f)^2}\,\frac{\partial\dot{\psi}\f}{\partial\normalf}}&= -\dot{\u}\e\cdot\normalf ,
&& \hspace{0.5cm} \textup{on } \GefI\times (0,T].
\end{alignat*}

\noindent\emph{Fluid-Tissue interface.} The type of the interface between the fluid and tissue depends on the choice of the material model for the tissue domain $\Omegat$.\\[1mm]
\noindent {\bf Option 1:  Elastic tissue}\\[1mm]
In case of an elastic material $\Omegat$, the interface has the same structure as the actuator-fluid interface.
\begin{alignat*}{2}
\vecc{\boldsymbol{\sigma}}(\u\t)~\normalt &= -\varrhof\,\left(\dot{\psi}\f+\frac{b\f}{(c\f)^2}\ddot{\psi}\f\right)~\normalt, && \hspace{0.5cm} \textup{on } \GftI\times (0,T],\\
\frac{\partial\psi\f}{\partial\normalf} {+ \frac{b\f}{(c\f)^2}\,\frac{\partial\dot{\psi}\f}{\partial\normalf}}&= -\dot{\u}\t\cdot\normalf,
&& \hspace{0.5cm} \textup{on } \Gamma_{\textup{I}}^{\textup{f}, \textup{t}}\times (0,T].
\end{alignat*}

\noindent {\bf Option 2: Acoustic tissue}\\[1mm]
In case of an acoustic material $\Omegat$, the type of the interface changes to the following transmission condition respecting jumps in the acoustic material parameters as well as enforcing continuity of the acoustic potential across the acoustic-acoustic interface in the trace sense.
\begin{alignat*}{2}
\left(\nabla\psi\t+\frac{b\t}{(c\t)^2}\nabla\dot{\psi}\t\right)\normalf&=\left(\nabla\psi\f+\frac{b\f}{(c\f)^2}\nabla\dot{\psi}\f\right)\normalf, && \hspace{0.5cm} \textup{on } \GftI\times (0,T],\\
\psi\t&=\psi\f, && \hspace{0.5cm} \textup{on } \GftI\times (0,T].
\end{alignat*}

\begin{comment}
\subsection{Notation}
We collect here the notation for different Hilbert spaces used through out the paper.
\begin{equation}
\begin{aligned}
& \Ha =\{v \in H^1(\Omega_\textup{a}): v=0 \ \text{on} \ \GeD\}\\
& \He =\{\vecc{v} \in \boldsymbol{H}^1(\Omega_\textup{e}): \boldsymbol{v}=0 \ \text{on} \ \GaD\}
\end{aligned}
\end{equation}
\end{comment}
As already discussed, we will focus on the stability and error analysis in the case of the tissue being an acoustic medium. The analysis when the tissue is elastic follows by analogous arguments, and so we omit it here. In the numerical experiments, we will extensively test both settings. We refer to Figure~\ref{fig:Domain} (bottom line) for a graphical description of the two options and corresponding domain setups. 
\begin{remark}[On the choice of the acoustic model] Our model for acoustic propagation in media with piecewise constant coefficients can be written in the strong form on $\Omegaa$ as
	\begin{equation} \label{general_wave}
	\begin{aligned}
	\tfrac{1}{c^2(x)}\ddot{\psi}-\Delta \psi- \textup{div}(\tfrac{b(x)}{c^2(x)} \nabla \dot{\psi}) = \tfrac{1}{c^2(x)}(k_1(x) \dot{\psi} \ddot{\psi}+ k_2(x) \nabla \psi \cdot \nabla \dot{\psi}).
	\end{aligned}
	\end{equation}
	This choice of the acoustic model is in part dictated by the physical interface conditions between the elastic and acoustic domain. Simplified versions of \eqref{general_wave} can be found in the nonlinear acoustics literature. In particular, equation
	\[
	\frac{1}{c^2(x)}\ddot{p}-\Delta p = \frac{\beta_{\textup{nl}}}{\varrho c^4} (p^2)_{tt}
	\]
	for the acoustic pressure is valid when the inhomogeneity of the medium varies in one direction; see~\citep[\S 5, Eq. (40)]{hamilton1998nonlinear}. In such cases, it can be assumed that this weak inhomogeneity of the medium results in changes only to the speed of sound, while the other parameters remain constant. Above, $\beta_{\textup{nl}}$ denotes the coefficient of nonlinearity in a given medium.  Integrating with respect to time and using the relation $p= \varrho \dot{\psi}$ gives
	\begin{equation} \label{particular_wave}
	\frac{1}{c^2(x)}\ddot{\psi}-\Delta \psi = \frac{\beta_{\textup{nl}}}{c^4} (\dot{\psi}^2)_t.
	\end{equation}
	Equation \eqref{general_wave} can be understood as a mathematical generalization of \eqref{particular_wave} that allows for a Kuznetsov-type nonlinearity and all coefficients to be piecewise constant functions. 
\end{remark}

\subsection{Notation}
Since in the analysis, we assume the tissue to be an acoustic medium, we can treat both the fluid and tissue as one acoustic domain with piecewise constant parameters. We thus denote this tissue-fluid acoustic region as
\[\overline{\Omegaa}= \overline{\Omegaf}\cup \overline{\Omegat}.\] 
The general elasto-acoustic interface $\GeaI$ and the acoustic-acoustic interface $\GaaI$ are then defined as \[\GeaI=\overline{\Omegaa}\cap\overline{\Omegae},\qquad \GaaI=\overline{\Omegaf}\cap\overline{\Omegat}.\]
We then define the parts of the boundary with homogeneous Dirichlet conditions prescribed as
\begin{equation*}
\begin{aligned}
\GaD= \, \partial \Omegaa \setminus  \Gamma_{\textup{I}}^{\textup{e,a}}, \qquad \GeD=\,  \partial \Omegae \setminus  \Gamma_{\textup{I}}^{\textup{e,a}}.
\end{aligned}
\end{equation*}
To facilitate the analysis, we introduce the following Hilbert spaces:
\begin{equation*}
\begin{aligned}
\Hf=& \,\{u \in H^1(\Omegaf): \ u=0 \ \text{on} \ \GfD\}, \qquad
\Ht= \,\{u \in H^1(\Omegat): \ u=0 \ \text{on} \ \GtD\}, \\
\Ha=& \,\{u \in H^1(\Omegaa): \ u=0 \ \text{on} \ \GaD\}, \qquad 
\He=\, \{\u \in \boldsymbol{H}^1(\Omegae): \ \u=\boldsymbol{0} \ \text{on} \ \GeD\},
\end{aligned}
\end{equation*}
\change{orange}{where $\boldsymbol{H}^1(\Omegae)$ stands for the vector-valued version of the $H^1(\Omegae)$ Sobolev space.}
%as well as the space
%\begin{equation}
%\begin{aligned}
%H^\Delta (\Omega_{\textup{i}})=& \,\{u \in L^2(\Omega_{\textup{i}}): \ \Delta u \in L^2(\Omega_{\textup{i}})\}, 
%\end{aligned}
%\end{equation}
%where $\textup{i} \in \{\textup{f}, \textup{t}, \textup{e}\}$. 
By $\norma{\cdot}$\change{orange}{and $\|\cdot\|_{\Omegae}$}, we denote the norm in $L^2(\Omegaa)$ \change{orange}{and $\boldsymbol{L}^2(\Omegae)$ again being the vector-valued version of $L^2(\Omegae)$:}
\[
\norma{\phi}=\left \{\int_{\Omegaa} |\phi|^2 \dx \right \}^{1/2},\qquad \change{orange}{\norme{\vecc{u}}=\left\{\int_{\Omegae} \|\vecc{u}\|^2\dx\right\}^{1/2}}.
\]
Similarly,  $\normf{\cdot}$ and  $\normt{\cdot}$ denote the norms in $L^2(\Omegaf)$ and $L^2(\Omegat)$, respectively.
\subsection{Assumptions on the medium parameters} We assume all the medium parameters to be piecewise constant functions: 
\begin{equation*}
\begin{aligned}
c(x) :=&\, \begin{cases}
c\f, ~~ x \in \Omegaf, \\
c\t, ~~ x \in \Omegat, \\
\end{cases} ~~ b(x) :=\begin{cases}
b\f, ~~ x \in \Omegaf, \\
b\t, ~~ x \in \Omegat, \\
\end{cases} ~~
k_1(x) :=&\,\begin{cases}
k_1\f, ~~ x \in \Omegaf, \\
k_1\t, ~~ x \in \Omegat, \\
\end{cases} ~~ k_2(x) :=\begin{cases}
k_2\f, ~~ x \in \Omegaf, \\
k_2\t, ~~ x \in \Omegat. \\
\end{cases}
\end{aligned}
\end{equation*}
Furthermore, we assume that $c\f$, $c\t >0$ and $b\f$, $b\t>0$. As mentioned in the introduction, the presence of strong $b$-damping in the nonlinear acoustic equation is crucial for the validity of our error estimates. The sign of $k$ does not play an important role in the analysis, and we can assume that $k_i \f$, $k_i \t \in \R$ for $i=1, 2$.  The piecewise acoustic nonlinearity is given by
\begin{equation}\label{eq:piecewiseAcousticForce}
f^{\textup{a}}(\dot{\psi},\psi,\nabla\psi,\nabla\dot{\psi}):=\begin{cases}
f^{\textup{f}}(\dot{\psi}^{\textup{f}},\psi^{\textup{f}},\nabla\psi^{\textup{f}},\nabla\dot{\psi}^{\textup{f}}), \quad x\in \Omegaf,\\
f^{\textup{t}}(\dot{\psi}^{\textup{t}},\psi^{\textup{t}},\nabla\psi^{\textup{t}},\nabla\dot{\psi}^{\textup{t}}), \quad x\in \Omegat.
\end{cases}
\end{equation} 
Since in the analysis we work with one elastic domain, from now on we drop the superscript $\textup{e}$ from the notation of the elastic solution and its derivatives, and use, for example, just $\vecc{u}$ in place of $\vecc{u}\e$. We also set $\zeta=\zeta\e >0$.\\
\indent The material density $\rho$ is the only parameter appearing in both elastic and acoustic media, hence we introduce the piecewise constant function
\begin{equation*}
\varrho(x)=\begin{cases}
\varrhoa(x)=\begin{cases}
\varrhof, \quad x\in\Omegaf,\\
\varrhot, \quad x\in\Omegat,
\end{cases}\\
\varrhoe,\qquad\qquad\hspace*{2mm}\quad x\in\Omegae,
\end{cases}
\end{equation*}
which is assumed to be positive as well.

\subsection{Weak formulation} We can now introduce the bilinear form \[\forma:\{\psi \in L^2(\Omegaa):  \psi_{\vert\Omegaf} \in \Hf, \ \psi_{\vert\Omegat} \in \Ht\}^2 \rightarrow \R\]  as
\begin{equation*}
\begin{aligned}
\forma(\psi, \phi)=\productf{\nabla \psi\f}{\nabla \phi\f}+\productt{\nabla \psi\t}{\nabla \phi\t},
\end{aligned}
\end{equation*}
where $\psi \vert_{\Omegaf}=\psi\f$ and  $\psi \vert_{\Omegat}=\psi\t$. The bilinear form $\forme:\He \times \He \rightarrow \R$ corresponding to the elastic medium is given by
\begin{equation*}
\begin{aligned}
\forme(\u, \v)=\producte{\mathbb{C}\boldsymbol{\varepsilon}(\u)}{ \boldsymbol{\varepsilon}(\v)}.
\end{aligned}
\end{equation*}
%where
%\begin{equation}
%\begin{aligned}
%\boldsymbol{\varepsilon}(\vecc{v}) =\frac12 \left(\nabla \vecc{v}+\nabla \vecc{v}^T\right), \quad \boldsymbol{\sigma}= \mathbb{C}\boldsymbol{\varepsilon}(\vecc{v}).
%\end{aligned}
%\end{equation}
Going forward, we assume that our problem has a solution in the following sense. Let the solution space for the displacement be given by
\begin{equation*}
\begin{aligned}
% X^{\textup{e}}= C([0,T]; \He \cap \boldsymbol{H}^2(\Omegae)) \cap C^1([0,T]; \He) \cap C^2([0,T]; L^2(\Omegae)), 
X^{\textup{e}}= C^1([0,T]; \He) \cap H^2(0,T;  \boldsymbol{H}^s(\Omegae)), 
\end{aligned}
\end{equation*}
and the solution space for the acoustic velocity potential by
\begin{equation*}
\begin{aligned}
X^{\textup{a}}=\{\psi \in C^1([0,T]; \Ha) \cap H^2(0,T; L^2(\Omegaa)):&\, \psi\vert_{\Omega_{\textup{i}}} \in H^2(0,T; H^s(\Omega_{\textup{i}}))\},
% X^{\textup{a}}=\{\psi \in C^1([0,T]; \Ha) \cap C^2([0,T]; L^2(\Omegaa)):&\, \psi\vert_{\Omega_{\textup{i}}} \in C([0,T]; H^s(\Omega_{\textup{i}})), \\
%&\dot{\psi}\vert_{\Omega_{\textup{i}}} \in L^2([0,T]; H^s(\Omega_{\textup{i}})), \\
%&\ddot{\psi}\vert_{\Omega_{\textup{i}}} \in L^2([0,T]; H^s(\Omega_{\textup{i}})), \ \textup{i} \in \{\textup{f}, \textup{t}\} \}.
\end{aligned}
\end{equation*}
where $s >1+d/2$, \change{orange}{$d$ being the spatial dimension of the problem}. We note that the choice of the regularity index $s$ is dictated by the nonlinear error analysis below; see Theorem~\ref{Thm:Main} for details.\\
\indent We assume that there exists $(\vecc{u}, \psi) \in X^{\textup{e}} \times X^{\textup{a}}$ such that
\begin{equation*} \label{weak_form}
\begin{aligned}
& \begin{multlined}[t]\producte{\varrhoe \ddot{\u}}{\v}+\producte{2 \varrhoe \zeta \dot{\u}}{\v}+\producte{\varrhoe \zeta^2 \u}{\v}+\forme(\u, \v)\\[0.2cm]
+\Iea(\varrhof(\dot{\psi}+\tfrac{b}{c^2} \ddot{\psi}), \v)-\Iae(\phi, \dot{\u})
+\producta{c^{-2} \ddot{\psi}}{\phi}+\forma(\psi+\tfrac{b}{c^2}\dot{\psi}, \phi)
\end{multlined}\\[0.2cm]
=&\, \producte{\fe}{\v}+ \producta{\fa(\dot{\psi}, \ddot{\psi}, \nabla \psi, \nabla \dot{\psi})}{\, \phi}
\end{aligned}
\end{equation*}
for a.e. $t \in (0,T)$ and all test functions $(\v, \phi) \in \He \times \{\phi \in L^2(\Omegaa):  \phi_{\vert\Omegaf} \in \Hf, \ \phi_{\vert\Omegat} \in \Ht\}$, supplemented by the initial conditions
\begin{equation*}
\begin{aligned}
(\u, \dot{\u})\vert_{t=0} =\, (\u_0, \u_1), \qquad
(\psi, \dot{\psi})_{t=0} =\, (\psi_0, \psi_1).
\end{aligned}
\end{equation*}

The interface form $\Iea:\Ha \times \He \rightarrow \R$ above is defined as
\begin{equation*}
\begin{aligned}
\Ief(\psi, \vecc{v})= (\psi \vecc{n}^\textup{e}, \v)_{\GeaI}.
\end{aligned}
\end{equation*}

Well-posedness of this coupled nonlinear problem appears to still be an open problem. The main difficulties in its analysis lie in the quasi-linear nature of the acoustic wave equation and the quadratic gradient nonlinearity, which, in general, requires the use of high-order energies in the analysis. However, some results on related linear coupled and nonlinear acoustic problems can be found in the literature. In particular, well-posedness of the linear undamped version ($b=k_1=k_2=0$) of the coupled problem with $\fa=\fa(x,t)$ is proven in~\citep[Theorem 1.1]{antonietti2018high}. Small-data global well-posedness of the Kuznetsov equation has been established in~\citep{mizohata1993global}, whereas the Westervelt equation expressed in terms of the acoustic pressure has been analyzed in, for example,~\citep{kaltenbacher2009global} and \citep{meyer2011optimal}. We also point out the analysis of a nonlinear elasto-acoustic problem in~\citep{brunnhuber2014relaxation}, where the acoustic field is modeled by the Westervelt equation with an additional strong nonlinear damping. 

\subsection{Discrete setting}
Before introducing the semi-discrete version of the above problem, we state the assumptions on the finite element mesh used in the discretization. Each subdomain $\Omegae$, $\Omegaf$, and $\Omegat$ is individually meshed in a conforming way by trilinear hexahedral Lagrange-elements, where the element-mapping for element $\kappa$ will be denoted by \[F_{\kappa}:(-1,1)^3\rightarrow\R^3, \quad F_{\kappa}((-1,1)^3)=\kappa.\] This gives rise to the subdomain-wise tessellations $\mathcal{T}_{h_{\textup{e}},\textup{e}}$, $\mathcal{T}_{h_{\textup{f}},\textup{f}}$, and $\mathcal{T}_{h_{\textup{t}},\textup{t}}$, where we have defined $h_{\textup{i}}:=\max_{\kappa\in\mathcal{T}_{h_{\textup{i}},\textup{i}}}h_{\kappa}$, $\textup{i}\in\lbrace \textup{e,f,t}\rbrace$, and the global mesh $\mathcal{T}_h$ with global mesh size \[h=\max\lbrace h_{\textup{e}},h_{\textup{f}},h_{\textup{t}}\rbrace.\] At the interfaces the subdomain meshes \emph{do not} have to match, which allows for a more flexible grid generation as well as different levels of refinement. For $\textup{i}\in\lbrace \textup{e,f,t}\rbrace$, we make the following assumptions on the families $\{\mathcal{T}_{h_{\textup{i}},\textup{i}}\}$:
\begin{itemize}
	\item Shape-regularity subdomain-wise: There exists $\sigma_{\textup{i}}>0$ such that each element $\kappa\in\mathcal{T}_{h_{\textup{i}},\textup{i}}$ contains a ball with the radius $\rho_{\kappa}\geq  \frac{h_{\kappa}}{\sigma_{\textup{i}}}$, where $h_{\kappa}$ denotes the diameter of the element $\kappa$. % and $\rho_{\kappa}$ the radius of the largest inscribed ball.
	\item Uniformity (subdomain-wise): There exist $\hat{\sigma}_{\textup{i}}>0$, such that: \(\frac{\displaystyle h_{\textup{i}}}{\displaystyle \min_{\kappa\in\mathcal{T}_{h_{\textup{i}},\textup{i}}}h_{\kappa}}\leq \hat{\sigma}_{\textup{i}}.\)
	\item Comparability of subdomains: There exists $\hat{\varsigma}_{\textup{i}}>0$ such that: \(h\leq \hat{\varsigma}_{\textup{i}} h_{\textup{i}}.\)
	\item Face non-degeneracy: Let $\kappa\in\mathcal{T}_{h_{\textup{i}},\textup{i}}$ and $F$ be any face of $\kappa$, then there exists some constant $\iota_{\textup{i}}$ such that $|F|\geq\frac{h_{\kappa}^2}{ \iota_{\textup{i}}}$.
\end{itemize}
\change{blue}{\begin{remark}
	We remark that the assumed face non-degeneracy is applied \textit{only} subdomain-wise and does \textit{not} exclude grids that are for example staggered by a length-parameter $\epsilon$ w.r.t. each other. The ``comparability of subdomain'' assumption can easily be satisfied even for individually generated sub-meshes as only the global mesh-size information needs to be exchanged.
\end{remark}}
Besides the computational domain's volumetric tessellation, we also define the sets of mesh faces that belong to the interfaces. As mentioned before, the analysis deals with Option 2, where the tissue is an acoustic medium. Hence, we classify the interfaces into elasto-acoustic and acoustic-acoustic. To deal with meshes that are non-matching at the interfaces, we follow ideas and notation of mortar-methods and classify master and slave sides for the subdomains meeting at the interfaces; see, for example,~\citep{BarbaraMortar}. In the case of the elasto-acoustic interface, the acoustic fluid region will be considered the master side, the elastic excitation domain the slave side. The fluid will be regarded as the master and the tissue domain as the slave side for the acoustic-acoustic interface again. \version{\change{blue}{For easier referencing} we define the interface-face-collections now \change{blue}{separately for both sides of the interface, where later summation/integration over the master side, of course, incorporates the other sides contributions as well, e.g. via jump-terms.}}{We define the interface-face-collection now only for the master side, where later summation/integration over them, of course, incorporates the respective slave sides contributions as well, e.g. via jump-terms.}
\version{\begin{align*}
\textup{Elasto-Acoustic interface faces (master side): }&\mathcal{F}_h^{\textup{e,a}}=\lbrace F : F\textup{~face of some~} \kappa\in\mathcal{T}_{h,\textup{f}} \wedge F\subset \Gamma_{\textup{I}}^{\textup{e,a}}\rbrace,\\
\change{blue}{\textup{Elasto-Acoustic interface faces (slave side): }}&\change{blue}{\mathcal{G}_h^{\textup{e,a}}=\lbrace G : G\textup{~face of some~} \kappa\in\mathcal{T}_{h,\textup{e}} \wedge G\subset \Gamma_{\textup{I}}^{\textup{e,a}}\rbrace,}\\
\textup{Acoustic-Acoustic interface faces: }&\mathcal{F}_h^{\textup{a,a}}=\lbrace F: F \textup{~face of some~} \kappa\in\mathcal{T}_{h,\textup{f}} \wedge F\subset \Gamma_{\textup{I}}^{\textup{a,a}}\rbrace.
\end{align*}}{
\begin{align*}
\textup{Elasto-Acoustic interface faces (master side): }&\mathcal{F}_h^{\textup{e,a}}=\lbrace F : F\textup{~face of some~} \kappa\in\mathcal{T}_{h,\textup{f}} \wedge F\subset \Gamma_{\textup{I}}^{\textup{e,a}}\rbrace,\\
\textup{Acoustic-Acoustic interface faces: }&\mathcal{F}_h^{\textup{a,a}}=\lbrace F: F \textup{~face of some~} \kappa\in\mathcal{T}_{h,\textup{f}} \wedge F\subset \Gamma_{\textup{I}}^{\textup{a,a}}\rbrace.
\end{align*}
}
If we talk about faces on a generic interface (e.g. any of elasto-acoustic or acoustic-acoustic), we use $\mathcal{F}_h^{\textup{I}}$.\\
\indent We can now use the following physical-domain elementwise ansatz-spaces:
\begin{align*}
\mathcal{Q}_p^F(\kappa)&=\lbrace \varphi:\kappa\rightarrow\R \,:\,\varphi=\hat{\varphi}\circ F_{\kappa}^{-1},~\textup{with some~}\hat{\varphi}\in \mathcal{Q}_p((-1,1)^3)\rbrace,\\
\underline{\mathcal{Q}}_p^F(\kappa)&=\lbrace v:\kappa\rightarrow\R^3 \,:\,v=\hat{v}\circ F_{\kappa}^{-1},~\textup{with some~}\hat{v}\in \underline{\mathcal{Q}}_p((-1,1)^3)\rbrace,
\end{align*}
where $\mathcal{Q}_p$ denotes the tensor-product polynomial space of degree $\leq p$ in each direction \change{orange}{and $\underline{\mathcal{Q}}_p$ its vector-valued version}. Furthermore, we introduce the finite-element spaces as
\begin{alignat*}{2}
V_h^{\textup{f}}&:=\lbrace \psi\in H^1_{\textup{D}}(\Omegaf) : \left.\psi\right|_{\kappa}\in \mathcal{Q}_{p}^F(\kappa)\,\forall \kappa\in \mathcal{T}_{h,\textup{f}}\rbrace, \qquad
&&V_h^{\textup{t}}:=\lbrace \psi\in H^1_{\textup{D}}(\Omegat) : \left.\psi\right|_{\kappa}\in \mathcal{Q}_{p}^F(\kappa)\,\forall \kappa\in \mathcal{T}_{h,\textup{t}}\rbrace,\\
V_h^{\textup{a}}&:=\lbrace \psi\in L^2(\Omegaa) : \left.\psi\right|_{\Omegaf}\in V_h^{\textup{f}}, \left.\psi\right|_{\Omegat}\in V_h^{\textup{t}}\rbrace, \qquad
&&\boldsymbol{V}_h^{\textup{e}}:=\lbrace v\in \boldsymbol{H}^1_{\textup{D}}(\Omegae) : \left.v\right|_{\kappa}\in \underline{\mathcal{Q}}_{p}^F(\kappa)\,\forall \kappa\in \mathcal{T}_{h,\textup{e}}\rbrace.
\end{alignat*}~\\[-1cm]
\begin{remark}
	Having different polynomial degrees $p_f$, $p_t$, and $p_e$ within the individual subdomains would be possible. However, for simplicity of notation, we restrict ourselves here to a single polynomial degree $p$.
\end{remark}

We use $x \lesssim y$ to denote $x \leq C y$, where the constant $C>0$ does not depend on the mesh size, however, might still depend on the material parameters and on the polynomial degree. For better readability we do not track this dependency on material parameters in this work. In an easier, purely acoustic setting such constant tracking was performed in \citep{antonietti2020high}.

\section{The semi-discrete problem}\label{sec:SemiDiscProb}
To state our hybrid semi-discrete problem, we introduce the following average and jump operators on the fluid-tissue interface. For sufficiently smooth $\psi$, we set the gradient average on $F \in \Fhaa$ and define the jump of the normal trace as follows:
\begin{equation*}
\begin{aligned}
\llbrace \nabla \psi \rrbrace =\frac{\nabla \psi \f +\nabla \psi \t}{2},\qquad \llbracket \psi \rrbracket = \psi \f \normalf +\psi \t \normalt.
\end{aligned}
\end{equation*}
%We also define the jump of the normal trace as
%\[
%\llbracket \psi \rrbracket = \psi \f \normalf +\psi \t \normalt.
%\]
Moreover, we introduce the short-hand notations
\begin{equation*}
\productFf{\psi}{\phi} = \sum_{F \in \Fhaa} (\psi, \phi)_{L^2(F)},\qquad \textup{and}\qquad \|\psi\|_{\Fhaa}= \productFf{\psi}{\psi}^{1/2}.
\end{equation*}
%and
%\begin{equation*}
%\|\psi\|_{\Fhaa}= \productFf{\psi}{\psi}^{1/2}.
%\end{equation*}
The semi-discrete problem is given by
\begin{equation} \label{weak_form_discrete}
\begin{aligned}
& \begin{multlined}[t]\producte{\varrhoe \ddot{\vecc{u}}_h(t)}{\vecc{v}_h}+\producte{2 \varrhoe \zeta \dot{\vecc{u}}_h(t)}{\vecc{v}_h}+\producte{\varrhoe \zeta^2 \vecc{u}_h(t)}{\vecc{v}_h}+\forme(\vecc{u}_h(t), \vecc{v}_h)\\[0.2cm]
+\producta{c^{-2} \ddot{\psi}_h(t)}{\phi_h}+\formah(\psi_h(t)+\tfrac{b}{c^2}\dot{\psi}_h(t), \phi_h)\\[0.2cm]+\Iea(\varrhof(\dot{\psi}_h(t)+\tfrac{b}{c^2} \ddot{\psi}_h(t)), \vecc{v}_h)-\Iae(  \phi_h, \dot{\vecc{u}}_h(t))
\end{multlined}\\[0.2cm]
=&\, \producte{\vecc{f}_\textup{e}(t)}{\vecc{v}_h}+ \producta{\fa_h(\dot{\psi}_h(t), \ddot{\psi}_h(t), \nabla \psi_h(t), \nabla \dot{\psi}_h(t))}{\, \phi_h}
\end{aligned}
\end{equation}
a.e. in time for all $(\vecc{v}_h, \phi_h) \in \Veh\times\Vah$ and supplemented with initial conditions
\begin{equation*}
\begin{aligned}
(\vecc{u}_h(0), \dot{\vecc{u}}_h(0), \psi_h(0), \dot{\psi}_h(0)) \in \Veh \times \Veh  \times \Vah \times  \Vah .
\end{aligned}
\end{equation*}
The acoustic gradient terms on the right-hand side of \eqref{weak_form_discrete} should be understood as
\begin{equation*}
\begin{aligned}
&\producta{\fa_h(\dot{\psi}_h(t), \ddot{\psi}_h(t), \nabla \psi_h(t), \nabla \dot{\psi}_h(t))}{\, \phi_h}\\
=&\,\productf{\ffh(\dot{\psi}_h\f(t), \ddot{\psi}_h\f(t), \nabla \psi_h\f(t), \nabla \dot{\psi}_h\f(t))}{\, \phi_h}+\productt{\fth(\dot{\psi}_h\t(t), \ddot{\psi}_h\t(t), \nabla \psi_h\t(t), \nabla \dot{\psi}_h\t(t))}{\, \phi_h},
\end{aligned}
\end{equation*}
\change{orange}{where $f_h^f$ and $f_h^t$ are defined analogously to $f^t$ and $f^t$ in \eqref{eq:AcRHS} just over the discrete spaces. This is used for a compact notation of the difference $\|f^a-f^a_h\|$ later in Sec. \ref{sec:LinearErrorEstimate} and \ref{sec:NonlinearArgument}.}\\
The discrete acoustic bilinear form $\formah: \Vah \times \Vah \rightarrow \R$ is given by
\begin{equation*}
\begin{aligned}
\formah(\psi_h, \varphi_h)
=&\, \begin{multlined}[t]\productf{\nabla \psi_h\f}{\nabla \varphi_h}+\productt{\nabla \psi_h\t}{\nabla \varphi_h}\\[1mm]
-\productFf{\llbrace \nabla \psi_h \rrbrace}{\llbracket \varphi_h \rrbracket} -\productFf{\llbrace \nabla \varphi_h \rrbrace}{\llbracket \psi_h \rrbracket} + \productFf{\chi \llbracket \psi_h \rrbracket}{\llbracket \varphi_h \rrbracket}. \end{multlined}
\end{aligned}
\end{equation*}
Finally, the stabilization parameter $\chi$ is defined face-wise on $F\in\mathcal{F}_{h}^{\textup{a,a}}$:
\begin{equation} \label{stabilization}
\begin{aligned}
\chi= \beta\frac{p^2}{h_F},\qquad h_F=\min\lbrace h_{\kappa}:\kappa\in\mathcal{T}_{h,\textup{t}},~\lambda^2(\overline{\kappa}\cap \overline{\kappa}_F)>0\rbrace,
\end{aligned}
\end{equation}
where $\kappa_F$ is the element that $F$ belongs to as a face and $\lambda^2(\cdot)$ is the two-dimensional Lebesgue measure. By that $h_F$ is the minimal $h$ of all elements from the slave-side (tissue domain) that have a non-trivial intersection with the face $F$, while $\beta>0$ is a suitable DG penalty parameter that will be chosen as sufficiently large to guarantee stability of the semi-discrete method. Finally we also incorporate the $p$-dependency in the stabilization parameter in the standard way via $p^2$, see~\citep{antonietti2018high,schotzau2002mixed,epshteyn2007estimation}, where $p$ is again the polynomial degree of the finite element ansatz functions.

\section{Stability analysis and a priori bounds of the linearized semi-discrete formulation}\label{sec:LinearErrorEstimate}
We first perform the stability analysis in the case
\begin{equation*}
\fah=\fah(x,t).
\end{equation*}
%We will only highlight the differences compared to the available analysis of the linearized problem in~\citep{antonietti2018high} and refer to~\citep{antonietti2018high} for further %details. In particular, 
To facilitate the later study of the nonlinear problem, the following points are important:
\begin{itemize}
	\item[(i)] Discretization errors $\fa-\fah$ of the acoustic source terms should be taken into account;
	\item[(ii)] An error estimate of $\ddot{\psi}_h$  and not only of the approximate potential field $\psi_h$ and its first time derivative $\dot{\psi}_h$ is needed. 
\end{itemize}
The second point implies that we have to involve the second time derivative of the approximate acoustic velocity potential as a test function in the energy analysis.\\
\indent Motivated by the analysis in~\citep{antonietti2020high}, we define the higher-order acoustic energy 
\begin{equation} \label{energy_acoustic}
\begin{aligned}
\normEa{\psi_h(t)}^2=&\,\begin{multlined}[t] \norma{\dot{\psi}_h(t)}^2+\int_0^t \norma{\ddot{\psi}_h(\tau)}^2 \, \textup{d}\tau
+ \norma{\nabla \tilde{\psi}_h(t)}^2 +\normFa{\tilde{\psi}_h(t)}^2 
, \end{multlined}
\end{aligned}
\end{equation}
where we have introduced the short-hand tilde notation \[\tilde{\psi}_h=\psi_h+\frac{b}{c^2}\dot{\psi}_h.\]
The gradient term should be understood in a broken sense as
\begin{equation*}
\begin{aligned}
\norma{\nabla \tilde{\psi}_h(t)}^2=&\, \forma(\tilde{\psi}_h(t), \tilde{\psi}_h(t))\\
=&\, \normf{\nabla (\psi_h\f+\frac{b\f}{(c\f)^{2}} \dot{\psi}_h\f)(t)}^2+\normt{\nabla (\psi_h\t+\frac{b\t}{(c\t)^{2}} \dot{\psi}_h\t)(t)}^2.
\end{aligned}
\end{equation*}
\change{orange}{For convenience of notation we also introduce the broken $H^s$-norm on the total acoustic domain as
	\begin{equation*}
	\normsa{\phi}^2 =\|\phi\|_{H^s(\Omegaf)}^2+\|\phi\|_{H^s(\Omegat)}^2,\quad s\geq 0.
	\end{equation*}
	%which will allow us to summarize the contributions of the two acoustic domains into one also for higher norms.} 
}We \change{orange}{further} recall that
\begin{equation*}
\begin{aligned}
 \normFa{\tilde{\psi}_h(t)}^2 =&\, \productFf{\chi \llbracket \tilde{\psi}_h \rrbracket}{\llbracket \tilde{\psi}_h \rrbracket}
\end{aligned}
\end{equation*}
with
\[
\llbracket \tilde{\psi}_h \rrbracket = (\psi_h \f+\tfrac{b \f}{(c \f)^2}\dot{\psi}_h\f) \normalf +(\psi_h \t+\tfrac{b \t}{(c \t)^2}\dot{\psi}_h\t) \normalt
\]
and the stabilization parameter $\chi$ defined in \eqref{stabilization}. The elastic energy is given by
\begin{equation} \label{energy_elastic}
\begin{aligned}
\normEe{\vecc{u}_h(t)}^2=&\,\begin{multlined}[t] \norme{\dot{\vecc{u}}_h(t)}^2+\norme{\uh(t)}^2
+ \norme{\boldsymbol{\varepsilon}(\uh(t))}^2. \end{multlined}
\end{aligned}
\end{equation}
We can then set the total energy as
\begin{equation*}
\normE{(\uh(t), \psi_h(t))}^2=\normEa{\psi_h(t)}^2+\normEe{\uh(t)}^2.
\end{equation*}
Note that except for the term $b/c^2 \dot{\psi}_h$ the norms within energies are chosen without scaling that would involve material parameters. For careful tracking of material parameters within the hidden constants in the numerical analysis, we refer to~\citep{antonietti2018high,antonietti2020high}.
\begin{remark} \label{Remark_b}
For the upcoming estimates, it is useful to note that a bound on the acoustic energy of a function $\psi$ at time $t$ will give us a bound on $\norma{\nabla \psi(t)}$ and $\norma{\nabla \dot{\psi}(t)}$ as well. Indeed, let
\[
\|\psi(t)\|^2_{\textup{E}^{\textup{a}}} \leq M
\]
for some $M>0$. Then by
		\begin{equation*}
	\begin{aligned}
	\norma{\nabla \tilde{\psi}(t)}^2=&\, \begin{multlined}[t] \norma{\nabla \psi(t)}^2 +\frac{(b\f)^2}{(c\f)^4}\normf{\nabla \dot{\psi}\f(t)}^2 +2\frac{b\f}{(c\f)^2}\productf{\nabla \psi_h\f(t)}{\nabla \dot{\psi}_h\f(t)} \\
+\frac{(b\t)^2}{(c\t)^4}\normf{\nabla \dot{\psi}\t(t)}^2 +2\frac{b\t}{(c\t)^2}\productt{\nabla \psi_h\t(t)}{\nabla \dot{\psi}_h\t(t)}	\leq M \end{multlined}
	\end{aligned}
	\end{equation*}
and H\"older's and Young's inequalities, we have
\begin{equation} \label{grad_est}
\begin{aligned}
\norma{\nabla \psi(t)}^2 +\frac{\underline{b}^2}{\overline{c}^4}\norma{\nabla \dot{\psi}(t)}^2 \leq&\, M+\frac{2\overline{b}}{\underline{c}^2}\norma{\nabla \psi(t)} \norma{\nabla \dot{\psi}(t)}\\
\leq&\,M+\frac{2}{\epsilon}\norma{\nabla \psi(t)}^2 + \epsilon \frac{\overline{b}^2}{2 \underline{c}^4}\norma{\nabla \dot{\psi}(t)}^2
\end{aligned}
\end{equation}
where $\underline{b}= \min \{b\f, b\t\}$,  $\overline{b}= \max \{b\f, b\t\}$ and similarly $\underline{c}= \min \{c\f, c\t\}$,  $\overline{c}= \max \{c\f, c\t\}$. We can choose $\epsilon>0$ small enough so that the last term on the right in \eqref{grad_est} is absorbed by the left side. By additionally relying on the bound
\[
\norma{\nabla \psi(t)}^2 \leq 2T \int_0^t \norma{\nabla \dot{\psi}(\tau)}^2 \textup{d}\tau+ 2 \norma{\nabla \psi(0)}^2 
\]
and Gronwall's inequality, we then have
	\begin{equation} \label{bound_grad}
\begin{aligned}
\norma{\nabla \psi(t)}^2 + \norma{\nabla \dot{\psi}(t)}^2 \leq \tilde{C}(T, b)  \left(M +\norma{\nabla \psi(0)}^2\right), \quad t\in [0,T],
\end{aligned}
\end{equation}
where the constant $\tilde{C}$ tends to infinity as $T \rightarrow \infty$ or $\underline{b} \rightarrow 0^+$, but does not depend on $h$.
\end{remark}	

\subsection{Preliminary theoretical results}\label{subsec:PrelimResults} In this subsection, we collect several well-known results on interpolation, trace inequalities, and standard estimates used in the (discontinuous) Galerkin framework that we will later need to employ in our proofs.
\begin{lemma}\label{lem:InverseEstimate}
	Let $p\in\N$ be a given polynomial degree and let $\kappa\in\mathcal{T}_h$. Then for $v\in\mathcal{Q}_p^F(\kappa)$, it holds
%	\begin{align} \label{eq:TraceInequality}
%	\|v\|_{L^2(\partial\kappa)}\lesssim h_{\kappa}^{-1/2}\|v\|_{L^2(\kappa)},\label{eq:TraceInequality} 
%	\end{align}	
	\begin{align} 		
	\|v\|_{L^{\infty}(\kappa)}&\leq \Cinv h_{\kappa}^{-d/2}\|v\|_{L^2(\kappa)},\label{inv_est}\\
		\|v\|_{L^2(\partial\kappa)}&\lesssim h_{\kappa}^{-1/2}\|v\|_{L^2(\kappa)}. \notag
	\end{align}
%	The hidden constant is then chosen as the largest of the ones coming from the three sub-domains.
\end{lemma}
%\begin{proof}
%	For the $L^{\infty}$-inverse estimate see Brenner Scott, Lemma 4.5.3, the trace $L^2$-estimate can be found in  \citep{DiPetroErnDG} Lemma 1.46. or \citep{evans2013explicit} Lemma 4.3. %Under imposed mesh assumptions, the statement follows as a special case of \citep[Lemma 1]{antonietti2020high}.
%\end{proof}
Note that the constants in Lemma~\ref{lem:InverseEstimate} depend on the polynomial degree of $v$. However, this dependency is suppressed here just as the material parameter dependencies.
%We first recall the following helpful inequalities, which we will rely on when proving the stability of our semi-discrete formulation.
\begin{lemma}\label{Lemma:ChiIneq} 	Let $\chi$ be the stabilization parameter defined in \ref{stabilization} with parameter $\beta$. Then it holds
	\begin{equation*}
	\|\chi^{-1/2}\llbrace\nabla \tilde{\psi}_h(t)\rrbrace\|_{\mathcal{F}_h^{\textup{a,a}}}\lesssim \frac{1}{\sqrt{\beta}}\|\nabla \tilde{\psi}_h(t)\|_{\Omegaa}, \quad t \in [0,T].
	\end{equation*}
\end{lemma}
\begin{proof}
	The proof follows analogously to the proof of~\citep[Lemma A.1]{antonietti2018high}.
\begin{comment}	 Indeed, owing to the trace inequality \eqref{eq:TraceInequality}, we have
	\begin{equation}
	\begin{aligned}
	\|\chi^{-1/2}\llbrace \nabla \tilde{\psi}_h(t) \rrbrace \|^2_{\Fhaa} =&\, \sum_{\kappa} \sum_{F \subset \partial \kappa} \|\chi^{-1/2}\llbrace \nabla \tilde{\psi}_h(t)\|^2_{L^2(F)}\\
	\lesssim&\,  \sum_{\kappa} \sum_{F \subset \partial \kappa} h_{F}^{-1} \|\chi^{-1/2}\llbrace \nabla \tilde{\psi}_h(t) \rrbrace\|^2_{L^2(\kappa)} \\
	\lesssim&\, \sum_{\kappa} \sum_{F \subset \partial \kappa} h_{F}^{-1} \frac{1}{\beta} \frac{p^2}{h_F} \|\llbrace \nabla (\tilde{\psi} \f + \tilde{\psi} \t)(t) \rrbrace\|^2_{L^2(\kappa)} \\
	\lesssim&\, \frac{1}{{\beta}}(\normf{\nabla \tilde{\psi}_h\f(t)}^2+\normt{\nabla \tilde{\psi}\t_h(t)}^2).
	\end{aligned}
	\end{equation}	
\end{comment}	
\end{proof}	
\begin{lemma}\label{Lemma:Ineq} \change{orange}{Let $\psi_h \in H^2(0,T; \Vah)$} and let $C_1$ and $C_2$ be two given positive constants. For a sufficiently large penalty parameter $\beta$ in \eqref{stabilization}, the following estimates hold:
	\begin{equation*}
	\begin{aligned}
	&\begin{multlined}[t]C_1 \normEa{\psi_h(t)}^2 -\productFa{\llbrace \nabla \tilde{\psi}_h(t) \rrbrace }{\llbracket \tilde{\psi}_h(t)  \rrbracket}
	\gtrsim \, \normEa{\psi_h(t)}^2, \end{multlined} \\
	&\begin{multlined}[t]C_2 \normEa{\psi_h(0)}^2 -\productFa{\llbrace \nabla \tilde{\psi}_h(0) \rrbrace }{\llbracket \tilde{\psi}_h(0)  \rrbracket}
	\lesssim \,  \normEa{\psi_h(0)}^2. \end{multlined}
	\end{aligned}
	\end{equation*}
\end{lemma}
\begin{proof}
\change{orange}{Note that $\psi_h \in H^2(0,T; \Vah) \hookrightarrow C^1([0,T]; \Vah)$.}	The proof can then be carried out analogously to the proof of Lemma A.2 in~\citep{antonietti2018high}.
	\begin{comment}
		By the Cauchy--Schwarz inequality and its discrete version, we have
	\begin{equation}
	\begin{aligned}	
	|\productFa{\llbrace \nabla \tilde{\psi}_h(t) \rrbrace }{\llbracket \tilde{\psi}_h(t)  \rrbracket} | \leq& \|\chi^{-1/2}\llbrace \nabla \tilde{\psi}_h(t) \rrbrace \|_{\Fhaa}\|\chi^{1/2}\llbracket \tilde{\psi}_h(t)  \rrbracket\|_{\Fhaa} \\
	\leq&\,  \frac{1}{4 \epsilon}\|\chi^{-1/2}\llbrace \nabla \tilde{\psi}_h(t) \rrbrace \|^2_{\Fhaa}+\epsilon \|\chi^{1/2}\llbracket \tilde{\psi}_h(t)  \rrbracket\|^2_{\Fhaa}.
	\end{aligned}
	\end{equation}
	Furthermore, thanks to the previous lemma, we have
	\begin{equation}
	\begin{aligned}
	\|\chi^{-1/2}\llbrace \nabla \tilde{\psi}_h(t) \rrbrace \|^2_{\Fhaa} \lesssim&\, \frac{1}{\beta}(\normf{\nabla \tilde{\psi}_h\f(t)}^2+\normt{\nabla \tilde{\psi}\t_h(t)}^2).
	\end{aligned}
	\end{equation}	
	Thus by setting $\epsilon=C_1/2$ and then choosing $\beta>0$ sufficiently large, the term $\productFa{\llbrace \nabla \tilde{\psi}_h(t) \rrbrace }{\llbracket \tilde{\psi}_h(t)  \rrbracket}$ can be absorbed by the acoustic energy, yielding the first bound. The second estimate is derived analogously.
\end{comment}	
\end{proof}~\\
The following lemmas summarize standard interpolation and stability estimates in the $L^2$, $H^1$, and $L^{\infty}$ norms.

\begin{lemma} \label{lem:InterpolationOperator} (see Theorem 4.6.14 in \citep{brenner2007mathematical})
	Let $\textup{i}\in\lbrace\textup{f,t}\rbrace$. There exists a subdomain-wise interpolation operator 
	$$\Pi^p:H^s(\Omega_{\textup{i}})\rightarrow V_h^{\textup{i}},$$
	which satisfies the following bounds:
	\begin{itemize}
		\item $\|\phi-\Pi^p\phi\|_{L^2(\Omega_{\textup{i}})}\leq C_{\textup{app2}} h^{s}\|\phi\|_{H^s(\Omega_{\textup{i}})}\qquad\quad\,\, \hspace*{2.5mm} \forall s=0,1,...,p+1$,\\
		\item $|\phi-\Pi^p\phi|_{H^1(\Omega_{\textup{i}})}\lesssim h^{s-1}\|\phi\|_{H^s(\Omega_{\textup{i}})}\qquad \qquad\,\,\, \quad\,\,\forall s=1,2,...,p+1$,\\ 
		\item $\|\phi-\Pi^p\phi\|_{L^{\infty}(\Omega_{\textup{i}})}\lesssim h^{s-d/2}\|\phi\|_{H^s(\Omega_{\textup{i}})}\qquad \qquad\,\,\, \forall~ d/2<s\leq p+1$.
	\end{itemize}
Component-wise application allows to extend the definition to vector-valued functions; i.e., there exists $\Pi^p:\boldsymbol{H}^s(\Omegae)\rightarrow\boldsymbol{V}_h^{e}$ with the same orders of approximation.
\end{lemma}
Due to the regularity assumptions on the exact solution to our problem, we can employ the Lagrange interpolation operator. Restriction to a single element $\kappa\in \mathcal{T}_h$  yields the following stability estimate.
%\begin{lemma} \label{lem:InterpolationOperatorOld} (see Theorem 4.6.11 in~\citep{brenner2007mathematical} and \citep{ciarlet1972interpolation}) Let $\kappa\in\mathcal{T}_h$ and $\phi\in H^s(\kappa)$, where $s\in\N_0$, and let $p\in\N$ be a given polynomial degree. Then there exists an interpolation operator \[\Pi_{\kappa}^p:H^s(\kappa)\rightarrow \mathcal{Q}^{F}_p(\kappa),\] which satisfies the following bounds:
%	\begin{itemize}
%		\item $\|\phi-\Pi_{\kappa}^p\phi\|_{L^2(\kappa)}\leq C_{\textup{app2}} h_{\kappa}^{s}\|\phi\|_{H^s(\kappa)}\qquad\quad\,\, \hspace*{2.5mm} \forall s=0,1,...,p+1$,\\
%		\item $|\phi-\Pi_{\kappa}^p\phi|_{H^1(\kappa)}\lesssim h_{\kappa}^{s-1}\|\phi\|_{H^s(\kappa)}\qquad \qquad\,\,\, \quad\,\,\forall s=1,2,...,p+1$,\\ 
%		\item $\|\phi-\Pi_{\kappa}^p\phi\|_{L^{\infty}(\kappa)}\lesssim h_{\kappa}^{s-d/2}\|\phi\|_{H^s(\kappa)}\qquad \qquad\,\,\, \forall~ d/2<s\leq p+1$.
%	\end{itemize}
%	Component-wise application allows to extend the definition to vector-valued functions; i.e., there exists $\Pi_{\kappa}^p:\boldsymbol{H}^s(\kappa)\rightarrow\boldsymbol{\mathcal{Q}}^F_p(\kappa)$ with the same orders of approximation.
%\end{lemma}
\begin{lemma}[See Lemma 4.4.1 in~\citep{brenner2007mathematical}]\label{lem:StabEst}
	The interpolation operator introduced in Lemma~\ref{lem:InterpolationOperator} fulfills the following stability estimate in the $W^{k,\infty}(\kappa)$ norm:
	$$
	\|\Pi^p\psi\|_{W^{k,\infty}(\kappa)}\leq C_{\textup{st}}\|\psi\|_{C^0(\overline{\kappa})}
	$$
\end{lemma}
\indent We also recall the following multiplicative trace inequality which will be used to derive error estimates.

\begin{lemma}\label{lem:MultiplicativeTraceInequality} Let $\kappa\in\mathcal{T}_h$ be a mesh element with diameter $h_{\kappa}$ satisfying our mesh assumptions. Furthermore, let $F$ be any face of $\kappa$. Then for $v\in H^1(\kappa)$ it holds
	\begin{equation*}
	\|v\|^2_{L^2(F)}\lesssim \|v\|_{L^2(\kappa)}\left(%2
	|v|_{H^1(\kappa)}+h_{\kappa}^{-1}\|v\|_{L^2(\kappa)}\right).
	\end{equation*}	
\end{lemma}
\begin{proof}
	The statement follows by Lemma 1.49 in \citep{DiPetroErnDG}.
\end{proof}
%For the later proof of the  error estimate for the linearization, we also recall the following two inequalities.
\begin{lemma}\label{lem:ChiGradientEstimate}
	Let $\phi_h\in\Vah$ and let $\chi$ be the stabilization parameter defined in \eqref{stabilization} with parameter $\beta$. Then for the global polynomial interpolant $\phi_I\in\Vah$ of degree $p$ of $\phi\in H^s(\Omegaa)$, $d/2<s\leq p+1$ it holds
	\begin{equation*}
	\|\chi^{-1/2}\llbrace\nabla(\phi-\phi_I)\rrbrace\|_{\mathcal{F}_h^{\textup{a,a}}}^2\lesssim h^{2(s-1)}\normsa{\phi}^2.
	\end{equation*}
\end{lemma}
\begin{proof}
The statement follows by Lemmas 3 and 4 in~\citep{antonietti2020high}.
\end{proof}

\begin{lemma}\label{lem:InterfaceEstimate}
	For any element $\kappa\in\mathcal{T}_h$ and polynomial degree $p\in\N$, let $\phi\in H^s(\kappa)$ for some $d/2<s\leq p+1$ and let $F$ be a face of $\kappa$. Then, the following interpolation estimate on the face $F$ holds true:
	\begin{equation*}
	\|\phi-\Pi^p\phi\|_{L^2(F)}\lesssim \sqrt{|F|}\,h_{\kappa}^{s-d/2}\,\|\phi\|_{H^s(\kappa)},\qquad\forall\, d/2<s\leq p+1
	\end{equation*}
\end{lemma}
\begin{proof}
The statement follows as a special case of the more general result of Lemma 4.2 in~\citep{CaGeHo2014};  see also Lemma 1.59 in~\citep{DiPetroErnDG}.
\end{proof}

With these technical results at hand, we now proceed with considering stability and error estimate for the linearization of the coupled problm.

\subsection{Stability of the semi-discrete formulation in the energy norm}\label{subsec:Stability}
%\begin{comment}
\noindent We first prove that our linearized semi-discrete approximation is \change{orange}{uniquely solvable and} stable in the energy norm.
\begin{proposition}\label{Prop:LinStability}
Let $\fe \in L^2(0,T; \boldsymbol{L}^2(\Omega))$, $\ffh \in L^2(0,T; \Vfh)$, and $\fth \in L^2(0,T; \Vth)$. The following estimate holds:
\begin{equation} \label{energy_est}
\begin{aligned}
\normE{(\uh(t), \psi_h(t))}^2 \lesssim&\, \begin{multlined}[t] \normE{(\uh(0), \psi_h(0)}^2+\int_0^t (\norme{\fe(\tau)}^2+\normf{\fah(\tau)}^2)\, \textup{d}\tau, \end{multlined}
\end{aligned}
\end{equation}
for a.e.\ $t \in [0,T]$, provided the penalty parameter $\beta$ in \eqref{stabilization} is sufficiently large. The hidden constant depends on the material parameters and the polynomial degree $p$, and tends to $\infty$ as as $T \rightarrow \infty$, but does not depend on the mesh size.
\end{proposition}
\begin{proof}
\change{orange}{We note that the existence of a unique $(\vecc{u}_h, \psi_h) \in H^1(0,T; \Veh) \times H^2(0,T; \Vah)$ follows by standard arguments for linear ODEs and the energy bounds derived below; see, for example,~\citep{nikolic2019priori}, \citep{antonietti2020high}. } \version{We test the problem by taking $\vecc{v}_h=\dot{\vecc{u}}_h$ and 
\[
\phi_h= \varrhof \dot{\tilde{\psi}}_h=\varrhof (\dot{\psi}_h+ \tfrac{b}{c^2}\ddot{\psi}_h).
\]
Note that we have the same factor $\varrho \f$ on both subdomains,
\[
\phi_h = \begin{cases}
\varrhof (\dot{\psi}_h\f+ \tfrac{b\f}{(c\f)^2}\ddot{\psi}_h\f) \quad \text{in } \Omegaf, \\[1mm]
\varrhof (\dot{\psi}_h\t+ \tfrac{b\t}{(c\t)^2}\ddot{\psi}_h\t) \quad \text{in } \Omegat.
\end{cases} 
\]
The reason for this choice of test functions is that they lead to the canceling out of the elasto-acoustic interface terms in \eqref{weak_form_discrete}:
\begin{equation*}
\begin{aligned}
& \productintea{\varrhof (\dot{\psi}\f_h+\tfrac{b\f}{(c\f)^2}\ddot{\psi}_h\f) \normale}{\dot{\vecc{u}}_h} = - \productintea{\dot{\vecc{u}}_h \cdot \normala}{\varrhof (\dot{\psi}_h\f+\tfrac{b\f}{(c\f)^2}\ddot{\psi}_h\f) }
\end{aligned}
\end{equation*}
because $\normala=-\normale$ on $\GeaI$. Moreover, scaling the acoustic test function by the constant $\varrhof$ (as opposed to $\varrho$) will not cause issues with the symmetry of dG terms across the acoustic interface.\\
\indent We are thus left with
\begin{equation*}
\begin{aligned}
&\begin{multlined}[t]\producte{\varrhoe \ddot{\vecc{u}}_h}{\dot{\vecc{u}}_h}+\producte{\varrhoe \zeta^2 \vecc{u}_h}{\dot{\vecc{u}}_h}+\producte{2 \varrhoe \zeta \dot{\vecc{u}}_h}{\dot{\vecc{u}}_h}+\producte{\boldsymbol{\sigma}(\vecc{u}_h)}{\boldsymbol{\varepsilon}(\vecc{\dot{u}}_h)} \\[1mm]
+\producta{c^{-2} \ddot{\psi}_h(t)}{\varrhof (\dot{\psi}_h+ \tfrac{b}{c^2}\ddot{\psi}_h)}+\formah(\psi_h(t)+\tfrac{b}{c^2}\dot{\psi}_h(t), \varrhof (\dot{\psi}_h+ \tfrac{b}{c^2}\ddot{\psi}_h))
 \end{multlined}\\[1mm]
=&\, \producte{\fe}{\dot{\vecc{u}}_h}+\producta{\fah}{\varrhof (\dot{\psi}_h+ \tfrac{b}{c^2}\ddot{\psi}_h)}.
\end{aligned}
\end{equation*}
Integrating the above identity over $(0,t)$ and employing integration by parts with respect to time leads to
\begin{equation*}
\begin{aligned}
& \begin{multlined}[t]
\frac12 \varrhoe \norme{\dot{\vecc{u}}_h(\tau)}^2 \Big \vert_0^t+\frac12 \varrhoe\zeta^2\norme{\vecc{u}_h(\tau)}^2 \Big \vert_0^t+2 \varrhoe \zeta\int_0^t \norme{\dot{\vecc{u}}_h}^2 \textup{d}\tau +\frac12 \norme{\mathbb{C}^{1/2}\boldsymbol{\varepsilon}(\vecc{u}_h(\tau))}^2 \Big \vert_0^t\\[1mm]
+\varrhof b\f (c\f)^{-4}\int_0^t \normf{\ddot{\psi}_h\f}^2 \textup{d}\tau+\frac12 \varrhof(c\f)^{-2} \normf{\dot{\psi}\f_h(\tau)}^2 \Big \vert_0^t+\frac12 \varrhof \normf{\nabla \tilde{\psi}_h\f(\tau)}^2 \Big \vert_0^t\\
+\varrhof  b\t (c\t)^{-4} \int_0^t \normt{\ddot{\psi}\t_h}^2 \textup{d}\tau+\frac12 \varrhof (c\t)^{-2} \normt{\dot{\psi}\t_h(\tau)}^2 \Big \vert_0^t+\frac12 \varrhof \norma{\nabla \tilde{\psi}\t_h(\tau)}^2 \Big \vert_0^t\\[1mm]
%-\varrhof\productFf{\llbrace \nabla  (\psi_h+\frac{b}{c^2}\dot{\psi}_h) \rrbrace}{\llbracket \dot{\psi}_h+\frac{b}{c^2}\ddot{\psi}_h \rrbracket} 
%\Big \vert_0^t+ \frac12\varrhof\productFf{\chi \llbracket \psi_h+\frac{b}{c^2}\dot{\psi}_h \rrbracket}{\llbracket \dot{\psi}_h+\frac{b}{c^2}\ddot{\psi}_h \rrbracket}\Big \vert_0^t\\[1mm]
{\change{orange}{-\varrhof\productFf{\llbrace \nabla  \tilde{\psi}_h(\tau) \rrbrace}{\llbracket \tilde{\psi}_h(\tau)\rrbracket} 
\Big \vert_0^t+ \frac12\varrhof\normFa{\tilde{\psi}_h(\tau)}^2\Big \vert_0^t}}
\end{multlined}\\
=&\, \int_0^t \producte{\fe}{\dot{\vecc{u}}_h} \textup{d}\tau +\int_0^t\producta{\fah}{\varrhof (\dot{\psi}_h+\frac{b}{c^2}\ddot{\psi}_h)} \textup{d}\tau.
\end{aligned}
\end{equation*}
By employing Young's inequality, we can estimate the last term on the right as follows:
\begin{equation*}
\begin{aligned}
&\int_0^t\producta{\fah}{\varrhof (\dot{\psi}_h+\frac{b}{c^2}\ddot{\psi}_h)} \textup{d}\tau \\
\leq&\, \begin{multlined}[t]\frac12  \int_0^t \norma{\fah}^2 \textup{d}\tau+ \frac12 (\varrhof)^2 \int_0^t \norma{\dot{\psi}_h}^2 \textup{d}\tau \\
+  \epsilon \left( \frac{(b \f)^2}{(c \f)^4} \int_0^t \normf{\ddot{\psi}_h}^2 \textup{d}\tau +\frac{(b \t)^2}{(c \t)^4} \int_0^t \normt{\ddot{\psi}_h}^2 \textup{d}\tau  \right)+  \frac{1}{4 \epsilon} \int_0^t \norma{\fah}^2 \textup{d}\tau.
\end{multlined}
%\lesssim&\, \begin{multlined}[t] C(\boldsymbol{\varepsilon}) \int_0^t \norma{\fah}^2 \textup{d}\tau+\boldsymbol{\varepsilon} \int_0^t \norma{\ddot{\psi}_h}^2 \textup{d}\tau
%+  \int_0^t \norma{\dot{\psi}_h}^2 \textup{d}\tau. \end{multlined}
\end{aligned}
\end{equation*}
Choosing $\epsilon >0$ small enough so that the $\ddot{\psi}_h$ terms on the right can be absorbed by the 5. and 8. term on the left, and recalling the definitions of the acoustic and elastic energies in \eqref{energy_acoustic} and \eqref{energy_elastic}, respectively, we arrive at
\begin{equation*}
\begin{aligned}
&\begin{multlined}[t]
\normEe{{\vecc{u}}_h(\tau)}^2 \Big \vert_0^t + \int_0^t \norme{\dot{\vecc{u}}_h}^2 \textup{d}\tau
+\normEa{\psi_h(\tau)}^2 \Big \vert_0^t -\productFa{\llbrace \nabla_h \tilde{\psi}_h(\tau) \rrbrace }{\llbracket \tilde{\psi}_h(\tau)  \rrbracket} \Big \vert_0^t
\end{multlined}\\
\lesssim&\, \begin{multlined}[t] \int_0^t \norme{ \fe}^2 \textup{d}\tau+\int_0^t \norme{\dot{\vecc{u}}_h}^2 \textup{d}\tau +\int_0^t \norma{ \fah}^2 \textup{d}\tau+ \int_0^t \norma{ \dot{\psi}_h}^2 \textup{d}\tau. \end{multlined}
\end{aligned}
\end{equation*}
Relying on the estimates of Lemma~\ref{Lemma:Ineq} yields
\begin{equation*}
\begin{aligned}
& \normEe{\vecc{u}_h(t)}^2 
+ \normEa{\psi_h(t)}^2 \\
\lesssim&\, \begin{multlined}[t] \int_0^t \norme{\fe}^2 \textup{d}\tau+\int_0^t \norme{\dot{\vecc{u}}_h}^2 \textup{d}\tau +\int_0^t \norma{\fah}^2 \textup{d}\tau+  \int_0^t \norma{ \dot{\psi}_h}^2 \textup{d}\tau. \end{multlined}
\end{aligned}
\end{equation*}
An application of Gronwall's inequality leads to the desired bound.}{\change{blue}{The proof then follows by testing the problem by $\vecc{v}_h=\dot{\vecc{u}}_h$ and 
	\[
	\phi_h= \varrhof \dot{\tilde{\psi}}_h=\varrhof (\dot{\psi}_h+ \tfrac{b}{c^2}\ddot{\psi}_h).
	\]
	Note that we have the same factor $\varrho \f$ on both subdomains,
	\[
	\phi_h = \begin{cases}
	\varrhof (\dot{\psi}_h\f+ \tfrac{b\f}{(c\f)^2}\ddot{\psi}_h\f) \quad \text{in } \Omegaf, \\[1mm]
	\varrhof (\dot{\psi}_h\t+ \tfrac{b\t}{(c\t)^2}\ddot{\psi}_h\t) \quad \text{in } \Omegat.
	\end{cases} 
	\]
	The reason for this choice of test functions is that they lead to the canceling out of the elasto-acoustic interface terms in \eqref{weak_form_discrete}:
	\begin{equation*}
	\begin{aligned}
	& \productintea{\varrhof (\dot{\psi}\f_h+\tfrac{b\f}{(c\f)^2}\ddot{\psi}_h\f) \normale}{\dot{\vecc{u}}_h} = - \productintea{\dot{\vecc{u}}_h \cdot \normala}{\varrhof (\dot{\psi}_h\f+\tfrac{b\f}{(c\f)^2}\ddot{\psi}_h\f) }
	\end{aligned}
	\end{equation*}
	because $\normala=-\normale$ on $\GeaI$. Moreover, scaling the acoustic test function by the constant $\varrhof$ (as opposed to $\varrho$) will not cause issues with the symmetry of dG terms across the acoustic interface. Standard computations then lead to \eqref{energy_est}. We omit the details here.}}
\end{proof}

\subsection{Error analysis of the linearized semi-discrete problem}\label{subsec:ErrorEstimate}
To facilitate the error analysis of the linearization, we define suitable norms in space and time as well as decompose the overall error into interpolation and discretization parts which are then estimated separately. We introduce the following norms:
\begin{equation} 
\normLinfEe{\v}= \displaystyle \change{dkgreen}{\esssup_{t\in (0,T)}}\normEe{\v(t)},\qquad \normLinfEa{\phi}=\sup_{t\in (0,T)}\normEa{\phi(t)}
\end{equation}
as well as the overall energy norm 
\begin{equation} \label{full_norms}
\normLinfE{(\v,\phi)}^2=\normLinfEe{\v}^2+\normLinfEa{\phi}^2.
\end{equation}
%\paragraph{\bf Error-decomposition:}
\noindent As standard, the total error between the solution and its approximation  \[e=(e_{\u},e_{\psi})=(\u-\u_h,\psi-\psi_h)\]
can be decomposed into two parts as follows:
\begin{align}\label{eq:ErrorDecomp}
{e}=(e_{\u},e_{\psi})&=(e_{\u,I}-e_{\u,h},e_{\psi,I}-e_{\psi,h})\\
&=((\u-\u_I)-(\u_I-\u_h),(\psi-\psi_I)-(\psi_I-\psi_h)), \notag
\end{align}
where $\vecc{u}_I$ and $\psi_I$ are the subdomain-wise defined, global interpolants of $(\vecc{u}, \psi) \in  X^{\textup{e}} \times  X^{\textup{a}}$ given by Lemma \ref{lem:InterpolationOperator}; i.e., $\left.\u_I\right|_{\Omegae} = \Pi^p\left.\u\right|_{\Omegae}$, $\left.\psi_I\right|_{\Omegaf}=\Pi^p\left.\psi\right|_{\Omegaf}$, and $\left.\psi_I\right|_{\Omegat}=\Pi^p\left.\psi\right|_{\Omegat}$.\\
\indent With these assumptions and the technical results from Subsection~\ref{subsec:PrelimResults} at hand, we can derive an approximation bound for the global interpolants. 
%To this end, we also introduce the norm
%\begin{equation*}
%\normsa{\phi}^2:=\|\phi\|_{H^s(\Omegaf)}^2+\|\phi\|_{H^s(\Omegat)}^2,
%\end{equation*}
%which allows us to summarize the contributions of the two acoustic domains into one.
\begin{lemma}\label{lem:InterpolationEstimate} Given $(\vecc{u}, \psi) \in  X^{\textup{e}} \times  X^{\textup{a}}$, the global interpolants satisfy the following error estimates: 
	\begin{align*}
	\normLinfEa{\psi-\psi_I}^2 &\lesssim h^{2(s-1)} \left\{\sup_{t\in (0,T)}\left(\normsa{\psi(t)}^2+\normsa{\dot{\psi}(t)}^2\right)+\int_0^t\normsa{\ddot{\psi}(\tau)}^2~\textup{d}\tau \right\},\\
	\normLinfEe{\u-\u_I}^2&\lesssim h^{2(s-1)}\sup_{t\in (0,T)}\left( \|\u(t)\|_{H^s(\Omegae)}^2+ \|\dot{\u}(t)\|_{H^s(\Omegae)}^2\right).
	\end{align*}
\end{lemma}
\begin{proof}
	\version{Employing Lemma~\ref{lem:InterpolationOperator} yields
	\begin{alignat*}{2}
	\normLinfEa{\psi-\psi_I}^2 \lesssim&\, \begin{multlined}[t] \sup_{t\in (0,T)} \left(h^{2(s-1)}\left(\normsa{\dot{\psi}(t)}^2+\int_0^t\normsa{\ddot{\psi}(\tau)}^2~\textup{d}\tau+\normsa{\tilde{\psi}(t)}^2\right)\right.\\
	\left. +\normFa{\tilde{\psi}(t)-\tilde{\psi}_I(t)}^2 \vphantom{h^{2(s-1)}}\right).\end{multlined}
	\end{alignat*}
	The integral term over the interface can be estimated using Lemma~\ref{lem:MultiplicativeTraceInequality} as well as again the interpolation estimates from Lemma~\ref{lem:InterpolationOperator}. In this manner, we obtain
	\begin{align*}
	\normFa{\tilde{\psi}-\tilde{\psi}_I}^2\lesssim&\, \norma{\sqrt{\chi}(\tilde{\psi}-\tilde{\psi}_I)}\left(2|\sqrt{\chi}(\tilde{\psi}-\tilde{\psi}_I)|_{H^1(\Omegaa)}+h^{-1}\norma{\sqrt{\chi}(\tilde{\psi}-\tilde{\psi}_I)}\right)\\
	\lesssim&\, h^{-1/2}_{F} h^{s}\normsa{\tilde{\psi}}\left(2h^{-1/2}_{F}  h^{s-1}\normsa{\tilde{\psi}}+h^{-1/2}_{F}  h^{-1}h^{s}\normsa{\tilde{\psi}}\right)\\
	\lesssim&\, \frac{h}{h_F} h^{2(s-1)}\normsa{\tilde{\psi}}^2,
	\end{align*}
	where we recall that $\chi=\beta\frac{p^2}{h_F}$; cf. \eqref{stabilization}. On account of our mesh assumptions, we then have
	\begin{align*}
	\normFa{\tilde{\psi}-\tilde{\psi}_I}^2 
	\lesssim\, h^{2(s-1)}\normsa{\tilde{\psi}}^2
	\end{align*}
	everywhere in time. Recalling also that $\tilde{\psi}=\psi+\frac{b}{c^2}\dot{\psi}$ yields the claimed estimate for the acoustic part. With analogous arguments using $\norme{\boldsymbol{\varepsilon}(\u-\u_I)}\lesssim \norme{\nabla(\u-\u_I)}$, we arrive at the respective result for the elastic part.}{Directly follows from Lemma ~\ref{lem:MultiplicativeTraceInequality} and the interpolation estimates from Lemma~\ref{lem:InterpolationOperator}.}
\end{proof}

We approximate initial conditions by applying the subdomain-wise global interpolation operators to the given data. In other words, we take $(\vecc{u}_h(0), \dot{\vecc{u}}_h(0), \psi_h(0), \dot{\psi}_h(0)) \in \Veh \times \Veh  \times \Vah \times  \Vah $ such that
\begin{equation} \label{approx_IC}
\begin{aligned}
& \u_h(0)=\u_{0,I}, \quad \dot{\vecc{u}}_h(0)=\u_{1,I}, \\
&\psi_h(0)=\psi_{0,I} , \quad \dot{\psi}_h(0)=\psi_{1,I}.
\end{aligned}
\end{equation}
We can now state the error bound in the energy norm for the linearized problem.  
\begin{theorem} \label{Thm:LinError} Let  $(\vecc{u}, \psi) \in  X^{\textup{e}} \times  X^{\textup{a}}$ and let the discrete initial conditions be obtained by interpolation of the exact ones; that is, let \eqref{approx_IC} hold. Then the approximation error $e=(e_{\u},e_{\psi})=(\u-\u_h,\psi-\psi_h)$ can be bounded as follows:
	\begin{equation} \label{lin_error_est}
	\begin{aligned}
	\|e\|_{L^{\infty}E}^2&\lesssim h^{2(s-1)}\bigg[\sup_{t\in(0,T)}\left(\|{\u}(t)\|_{H^{s}(\Omegae)}^2+\normsa{{{\psi}(t)}}^2+\normsa{\dot{{\psi}}(t)}^2\right)\\
	&\qquad\qquad\qquad+\int_0^T\bigg( \|\dot{\u}\|_{H^s(\Omegae)}^2+\|\ddot{\u}\|_{H^s(\Omegae)}^2+\normsa{\ddot{\psi}}^2\bigg)~\textup{d}\tau\bigg]
+\int_0^T\|\fa-\fah\|_{\Omegaa}^2~\textup{d}\tau,
	\end{aligned}
	\end{equation}
	provided the penalty parameter $\beta$ in \eqref{stabilization} is sufficiently large. The hidden constant in the estimate tends to infinity as $T \rightarrow \infty$, but does not depend on $h$.	
\end{theorem}
The proof follows by combining the arguments from the proof of Theorem 5.2 in \citep{antonietti2018high} (with respect to treating the elastic and coupling terms) with our particular choice of test functions, as already seen in Proposition \ref{Prop:LinStability}, and allowing for the error in the source term. \version{For self-consistency, we have worked out the technical details in Appendix A.}{For a worked out version of the proof, see Appendix A. in \cite{arxivVersion}.}

\section{A priori analysis of the nonlinear coupled problem} \label{sec:NonlinearArgument}
We next analyze the nonlinear problem by relying on the Banach fixed-point theorem; see, for example,~\citep{antonietti2020high, ortner2007discontinuous} for similar techniques used in the numerical analysis of nonlinear wave equations. To this end, we introduce the mapping
	\[\mathcal{S}: \mathcal{B}_h \ni  (\tilde{\vecc{u}}_h, w_h) \mapsto (\vecc{u}_h, \psi_h),\]
where $(\vecc{u}_h, \psi_h)$ solves the linear problem
\begin{equation} \label{linearization}
\begin{aligned}
& \begin{multlined}[t]\producte{\varrhoe \ddot{\u}_{h}(t)}{\v_h}+\producte{2 \varrhoe \zeta \dot{\u}_h(t)}{\v_h}+\producte{\varrhoe \zeta^2 \u_h(t)}{\v_h}+\forme(\u_h(t), \v_h)\\[0.1cm]
+\producta{c^{-2} \ddot{\psi}_h(t)}{\phi_h}+\formah(\psi_h(t)+\tfrac{b}{c^2}\dot{\psi}_h(t), \phi_h)\\
+\mathcal{I}(\varrhof(\dot{\psi}_h(t)+\tfrac{b}{c^2}\ddot{\psi}_h), \v_h)-\mathcal{I}(\phi_h,\dot{\u}_h(t))
\end{multlined}\\[0.1cm]
=&\, \producte{\fe(t)}{\v_h}+ \producta{\fah(\dot{w}_h(t),\ddot{w}_h(t), \nabla w_h(t), \nabla \dot{w}_h(t))}{\, \phi_h},
\end{aligned}
\end{equation}
a.e.\ in time for all test functions $(\vecc{v}_h, \phi_h) \in \Veh \times \Vah $ and supplemented with initial conditions \eqref{approx_IC}. Recall that
\begin{equation*}
\fah(\dot{w}_h,\ddot{w}_h, \nabla w_h, \nabla \dot{w}_h)= \frac{2}{c^2}\left(k_1 \dot{w}_h \ddot{w}_h+k_2 \nabla w_h \cdot \nabla \dot{w}_h\right)
\end{equation*}
with material parameter jumps allowed at the fluid-tissue interface; cf. \eqref{eq:piecewiseAcousticForce}.\\
\indent Furthermore, given $C_{\star}>0$, $\mathcal{B}_h$ is the ball containing all $(\tilde{\vecc{u}}_h, w_h) \in H^1(0,T; \Veh) \times H^2(0,T; \Vah)$ such that
\begin{equation*}
\begin{aligned}
\change{dkgreen}{\|(\u-\tilde{\vecc{u}}_h,\psi-w_h)\|_{L^\infty E}}
 \leq\, \begin{multlined}[t] C_\star h^{s-1} \left\{ \tripplenorme{\u}
+   \tripplenorma{\psi}\right\},  \end{multlined}
\end{aligned}
\end{equation*}
with initial conditions
\begin{equation*} \label{IC_nlproblem}
\begin{aligned}
(w_h(0), \dot{w}_h(0))=(\psi_{{0, I}}, \psi_{{1,I}}), \qquad (\tilde{\vecc{u}}_h(0), \dot{\tilde{\vecc{u}}}_h(0))=(\u_{0, {I}}, \u_{1,{I}}).
\end{aligned}
\end{equation*}
Above, we have introduced the following short-hand notation:
\begin{equation*}
\tripplenorme{\u}^2=\sup_{t\in (0,T)}\|\u(t)\|^2_{H^s(\Omegae)}+\int_0^T (\|\dot{\u}\|^2_{H^s(\Omegae)}+\|\ddot{\u}\|^2_{H^
	s(\Omegae)})\textup{d}\tau ,
\end{equation*}
and
\begin{equation*}
\tripplenorma{\psi}^2=\sup_{t\in (0,T)}\left(\normsa{\psi(t)}^2+\normsa{\dot{\psi}(t)}^2\right)+\int_0^T \normsa{\ddot{\psi}}^2\textup{d}\tau,
\end{equation*}
where we recall that the $s$-regularity above should be understood subdomain-wise:
\begin{equation*}
\normsa{\phi}^2=\|\phi\|_{H^s(\Omegaf)}^2+\|\phi\|_{H^s(\Omegat)}^2.
\end{equation*}
The constant $ C_\star=C_\star(\vecc{u}, \psi)$ in the error estimate will be specified below. The set $\mathcal{B}_h$ is non-empty because the solution of the linear problem belongs to it when $\fah=\fa$ for a suitably chosen $C_\star$.\\
\indent On account of the existence and uniqueness result for the linear problem, this mapping is well-defined. Furthermore, the space $(\mathcal{B}_h, d)$ is complete with respect to the metric $d((\vecc{u}, \psi), (\vecc{v}, \phi))=\normE{(\vecc{u}-\vecc{v}, \psi-\phi)}$.\\
\indent We first determine sufficient conditions for $\mathcal{S}$ to be a self-mapping. To simplify notation, we use  $\|{\cdot}\|_{L^pL^q}$ below in place of $\|\cdot\|_{L^p(0,T; L^q(\Omega_{\textup{i}}))}$ for $p$, $q\in\lbrace 2,\infty\rbrace$. We will also rely on the continuous embedding \[H^s(\Omega_{\textup{i}}) \hookrightarrow W^{1, \infty}(\Omega_{\textup{i}}), \quad s >1+d/2,\] with the embedding constant denoted by $C_{\textup{emb}}.$

%%%%%%%%%%%%%%%%%%%%%%%%%%%%%%%%%%%%%%%%%%%%%%%%%%%%%%%%%%%%%%%%%%%%%%%%%%%%%%%%%%%%%%%%%%%%%%%%%%%%%%%%
\begin{proposition} \label{Prop:SelfMapping}
Let $1 + d/2 <s \leq  p+1$ and $\fe \in L^2(0,T; \boldsymbol{L}^2(\Omega))$. Assume that the penalty parameter $\beta$ in \eqref{stabilization} is chosen as sufficiently large according to Proposition~\ref{Prop:LinStability} and Theorem~\ref{Thm:LinError}. Then there exist $\overline{h}>0$ and $\delta>0$, such that for
%\begin{equation} \label{fixed_point_norm}
%\begin{aligned}
%\|(\vecc{u}, \psi)\|_X =\begin{multlined}[t] \tripplenorme{\u} +   \tripplenorma{\psi}+\|\dot{\psi}\|_{L^\infty L^{\infty}}+\| \ddot{\psi}\|_{L^2 L^\infty}\\+ %\|\dot{\psi}\|_{L^2 L^{\infty}}+\|\nabla\dot{\psi}\|_{L^2 L^{\infty}}+\|\psi\|_{L^2 L^\infty} \leq \delta
%+\sup_{t\in (0,T)}(\normsa{\psi}+\normsa{\dot{\psi}})+(\int_0^T\normsa{\ddot{\psi}}^2~\textup{d}s))^{1/2}
%\end{multlined}
%\end{aligned}
%\end{equation}
\begin{equation*} \label{fixed_point_norm}
\begin{aligned}
\tripplenorme{\u} +   \tripplenorma{\psi} \leq \delta
\end{aligned}
\end{equation*}
and $0<h \leq \overline{h}$, the mapping $\mathcal{S}$ satisfies \[\mathcal{S}(\mathcal{B}_h) \subset \mathcal{B}_h.\]
\end{proposition}
\begin{proof}
	For a given $ (\tilde{\vecc{u}}_h, w_h) \in \mathcal{B}_h$ and $(\vecc{u}_h, \psi_h)$ solving \eqref{linearization}, we know from the linear result that there exists $C_{\textup{lin}}>0$ such that
	\begin{equation} \label{lin_est}
	\begin{aligned}
	\normLinfE{(\u-\u_h, \psi-\psi_h)} \leq \,  C_{\textup{lin}} \left(h^{s-1}\left( \tripplenorme{\u}+\tripplenorma{\psi}\right)+\LtLt{\fa-\fah}\right)
	\end{aligned}
	\end{equation}
	provided $\beta>0$ in \eqref{stabilization} is large enough. Since 
	\begin{equation*}
	\begin{aligned}
	\fa-\fah=&\, \frac{2}{c^2}\left(k_1(\dot{\psi}\ddot{\psi}-\dot{w}_h\ddot{w}_h)+k_2(\nabla \psi \cdot \nabla \dot{\psi}-\nabla w_h \cdot \nabla \dot{w}_h)\right) \\
	=&\, \begin{multlined}[t]\frac{2}{c^2}\left(k_1(\dot{\psi}(\ddot{\psi}-\ddot{w}_h)+\ddot{w}_h(\dot{\psi}-\dot{w}_h)) \right.  \left. +k_2(\nabla w_h \cdot \nabla (\dot{\psi}-\dot{w}_h) + \nabla \dot{\psi}\cdot \nabla (\psi - w_h)\right),\end{multlined}
	\end{aligned}
	\end{equation*}
	we can estimate the error in the acoustic source term \change{orange}{using $|\overline{k_i}|:=\max\lbrace|k_i\f|,|k_i\t|\rbrace, i=1,2$} as follows: % using $\dot{a}\ddot{a}-\dot{b}\ddot{b}=\dot{a}(\ddot{a}-\ddot{b})+\ddot{b}(\dot{a}-\dot{b})$ twice:
	\begin{equation*} 
	\begin{aligned}
	\LtLt{\fa-\fah}
	\leq &\,  \begin{multlined}[t]  \frac{2|\change{orange}{\overline{k_1}}|}{\underline{c}^2}\left(\LtLt{\dot{\psi}(\ddot{\psi}-\ddot{w}_h)}+\LtLt{\ddot{w}_h(\dot{\psi}-\dot{w}_h)}\right)\\ 
	+\frac{2|\change{orange}{\overline{k_2}}|}{\underline{c}^2}\left(\LtLt{\nabla w_h \cdot \nabla (\dot{\psi}-\dot{w}_h)}+\LtLt{\nabla \dot{\psi} \cdot \nabla (\psi - w_h)}\right),  \end{multlined}
	\end{aligned}
	\end{equation*}	
	where the gradient should be understood in a broken $\Omegaf \cup \Omegat$ sense. Recall that $\underline{c}=\min\{c\f, c\t\}$. Therefore, we infer
	\begin{equation}  \label{est_f-fh}
	\begin{aligned}
	& \LtLt{\fa-\fah} \leq &\, \begin{multlined}[t] \frac{2|\change{orange}{\overline{k_1}}|}{\underline{c}^2}\left(\LinfLinf{\dot{\psi}}\LtLt{\ddot{\psi}-\ddot{w}_h}+\LtLinf{\ddot{w}_h}\LinfLt{\dot{\psi}-\dot{w}_h}\right)\\ 
	+\frac{2|\change{orange}{\overline{k_2}}|}{\underline{c}^2}\left(\LtLinf{\nabla w_h}\LinfLt{\nabla (\dot{\psi}-\dot{w}_h)}+\LtLinf{\nabla \dot{\psi}}\LinfLt{\nabla (\psi - w_h)}\right). \end{multlined}
	\end{aligned}
	\end{equation}	
	We know the approximation error of $(\tilde{\vecc{u}}_h, w_h) \in \mathcal{B}_h$, so we can further deduce  
	\begin{equation*} %\label{est_f-fh}
	\begin{aligned}
	\LtLt{\fa-\fah} \leq\, \frac{2\overline{k}}{\underline{c}^2} \left(\LinfLinf{\dot{\psi}}+\LtLinf{\ddot{w}_h}+\LtLinf{\nabla w_h} +\LtLinf{\nabla \dot{\psi}} \right) 
	  C_\star h^{s-1} \left( \tripplenorme{\u}
	+  \tripplenorma{\psi}\right), 
	\end{aligned}
	\end{equation*}
	where $\overline{k}= \max \{\change{orange}{|\overline{k_1}|, |\overline{k_2}|}\}$. From \eqref{lin_est} we see that for the self-mapping property to hold, we have to guarantee that 
	\begin{equation*}
	\check{C}:= C_{\textup{lin}}\left[1+\frac{2\overline{k}}{\underline{c}^2} \bigg(\LinfLinf{\dot{\psi}}+\LtLinf{\ddot{w}_h}+\LtLinf{\nabla w_h} + \LtLinf{\nabla \dot{\psi}}\bigg)C_\star \right] \leq C_\star.
	\end{equation*}
We will next further bound the $w_h$ terms by relying on the inverse estimates together with the stability and approximation properties of the interpolant. Note first that
\[
\|\dot{\psi}\|_{L^\infty L^\infty} \leq C_{\textup{emb}}\tripplenorma{\psi}, \qquad \|\nabla \dot{\psi}\|_{L^2 L^\infty} \leq \sqrt{T}C_{\textup{emb}}\tripplenorma{\psi}.
\] 
Let $\kappa \in \mathcal{T}_{h_{\textup{f}},\textup{f}} \cup \mathcal{T}_{h_{\textup{t}},\textup{t}}$ be the element such that
	\[
\begin{aligned}
		\|\nabla w_h\|_{L^2L^\infty} =&\, \, \|\nabla w_h\|_{L^2L^\infty(\kappa)} \leq&\, \|\nabla w_h- \nabla \Pi^{p} \psi\|_{L^2L^\infty(\kappa)} +\|\nabla \Pi^{p} \psi\|_{L^2 L^\infty(\kappa)}.
	\end{aligned}
	\]
	 Then thanks to the inverse estimate \eqref{inv_est}, we find that 
	\begin{equation} \label{est_nabla_phi_h}
	\begin{aligned}
	\LtLinf{\nabla w_h} \leq&\, h^{-d/2}\Cinv \|\nabla w_h -\nabla \Pi^{p} \psi\|_{L^2L^2(\kappa)} +\|\nabla \Pi^{p} \psi\|_{L^2 L^\infty(\kappa)}\\
%	\leq&\, h^{-d/2}{\Cinv} (\|\nabla w_h -\nabla \psi_I \|_{L^2L^2(\Omegaf)}+\|\nabla w_h -\nabla \psi_I \|_{L^2L^2(\Omegat)}) +\|\nabla \Pi^{p} \psi\|_{L^2 L^\infty(\kappa)}\\
	\leq&\, \tilde{C}\Cinv h^{-d/2} \| w_h - \psi_I \|_{L^2 \textup{E}^{\textup{a}}}+\|\nabla \Pi^{p} \psi\|_{L^2 L^\infty(\kappa)},
	\end{aligned}
	\end{equation}
	where in the second line we have used the bound \eqref{bound_grad} on the gradient via the acoustic energy and relied on our choice of approximate initial data. Recall that the constant $\tilde{C}= \tilde{C}(T, b)$ tends to infinity as $T \rightarrow +\infty$ or $\underline{b} \rightarrow 0^+$; cf. Remark~\ref{Remark_b}. From here we have
	\begin{equation} \label{est_nabla_phi_h}
	\begin{aligned}
	\LtLinf{\nabla w_h} \leq&\, \begin{multlined}[t]  \tilde{C}\Cinv h^{-d/2} \|w_h-\psi + \psi -\psi_I\|_{L^2 \textup{E}^{\textup{a}}} +\|\nabla \Pi^{p} \psi\|_{L^2 L^\infty(\kappa)}\end{multlined}\\
	\leq&\, \begin{multlined}[t] \tilde{C}\Cinv h^{-d/2} \sqrt{T} C_\star h^{s-1} ( \tripplenorme{\u}
	+   \tripplenorma{\psi}) \\+ \tilde{C}\Cinv h^{-d/2} \sqrt{T} C_{\textup{app}} h^{s-1}\tripplenorma{\psi} +C_{\textup{st}}\|\psi\|_{L^2 L^\infty},\end{multlined}
	\end{aligned}
	\end{equation}
	where we have relied on the approximation properties of the global interpolant in the energy norm with $C_{\textup{app}}$ being the hidden constant therein; see Lemma~\ref{lem:InterpolationEstimate}. Also we have used the stability of the local interpolant in the  $W^{1, \infty}$ norm; see Lemma \ref{lem:StabEst}. Finally,
	\[
        C_{\textup{st}}\|\psi\|_{L^2 L^\infty} \leq \sqrt{T}C_{\textup{st}} C_{\textup{emb}}\tripplenorma{\psi}.
	\]
	Similarly, it holds that
	\begin{equation}\label{est_ddotpsih}
	\begin{aligned}
		\|\ddot{w}_h\|_{L^2L^{\infty}}\leq&\, \Cinv h^{-d/2}\left(\|\ddot{w}_h-\ddot{\psi}\|_{L^{2}L^2}+\|\ddot{\psi}-(\ddot{\psi})_I\|_{L^{2}L^2}\right)+\|\Pi^p\ddot{\psi}\|_{L^2 L^\infty(\kappa)}\\
		\leq&\, \begin{multlined}[t] \Cinv h^{-d/2}\big(C_\star  h^{s-1} ( \tripplenorme{\u}
		+   \tripplenorma{\psi})
		 +C_{\textup{app}} h^{s-1}\tripplenorma{\psi}\big)+C_{\textup{st}}\|\ddot{\psi}\|_{L^2 L^\infty},\end{multlined}
	\end{aligned}	
	\end{equation}
%	\begin{equation} 
%	\begin{aligned}
%	\|\ddot{w}_h\|_{L^2 L^\infty} \leq& \, \begin{multlined}[t] h^{-d/2} C_\star(\vecc{u}, \psi)  h^{s-1} ( \tripplenorme{\u}
%	+   \tripplenorma{\psi})\\  + h^{-d/2}C_{\textup{app2}} h^{s}(\|\ddot{\psi}\|_{L^2 H^{s}(\Omegaf)}+\|\ddot{\psi}\|_{L^2 H^{s}(\Omegat)})\\+C_{\textup{st}}\|\ddot{\psi}\|_{L^2 L^\infty}.\end{multlined}
%	\end{aligned}
%	\end{equation}
\begin{comment}	
Altogether, we have
	\begin{equation}
	\begin{aligned}
	\tilde{C} \leq &\, \begin{multlined}[t] C_{\textup{lin}}\bigg[1+\frac{2k}{c^2} \bigg(\|\dot{\psi}\|_{L^\infty L^{\infty}}+ \|\nabla \dot{\psi}\|_{L^2 L^{\infty}}+ \tilde{C} h^{-d/2} \sqrt{T} C_\star h^{s-1} ( \tripplenorme{\u} +   \tripplenorma{\psi}) \\+ \tilde{C}h^{-d/2} C_{\textup{app}} h^{s-1}(\|\psi\|_{L^2 H^{s}(\Omegaf)}+\|\psi\|_{L^2 H^{s}(\Omegat)}) +C_{\textup{st}}\|\nabla \psi\|_{L^2 L^\infty}\\
	+ h^{-d/2}C_\star h^{s-1} ( \tripplenorme{\u}
	+   \tripplenorma{\psi}) \\+ h^{-d/2} C_{\textup{app}} h^{s}(\|\ddot{\psi}\|_{L^2 H^{s}(\Omegaf)}+\|\ddot{\psi}\|_{L^2 H^{s}(\Omegat)}) \\+C_{\textup{st}}\| \ddot{\psi}\|_{L^2 L^\infty} \bigg)C_\star\bigg]. \end{multlined}
	\end{aligned}
	\end{equation}
\end{comment}	
where the last term can be further bounded as follows:
\[C_{\textup{st}}\|\ddot{\psi}\|_{L^2 L^\infty} \leq C_{\textup{st}}C_{\textup{emb}}\tripplenorma{\psi}. \]
Let us collect the embedding, stability, and approximation constants appearing above into one constant given by
\[
\mathcal{C}(T) = C_{\textup{emb}}(C_{\textup{st}}+1)(1+\sqrt{T})+\Cinv(C_{\textup{app}}+1)(1+\tilde{C}\sqrt{T}).
\]
Altogether, we then have 
\begin{comment}
\begin{equation}
\begin{aligned}
\check{C} \leq &\, \begin{multlined}[t] C_{\textup{lin}}\left\{1+\frac{2\overline{k}}{\underline{c}^2} ( \tripplenorme{\u}
+   \tripplenorma{\psi})\left[C_{\textup{emb}}(C_{\textup{st}}+1)(1+\sqrt{T}) + \tilde{C}\Cinv h^{s-1-d/2}\sqrt{T}  C_\star(\vecc{u}, \psi) \right. \right. \\ \left. \left.+ \tilde{C}\Cinv h^{s-1-d/2} \sqrt{T} C_{\textup{app}} 
+\Cinv h^{s-1-d/2} C_\star(\vecc{u}, \psi) + \Cinv h^{s-d/2} C_{\textup{app{2}}}\right]C_\star\right\}. \end{multlined}
\end{aligned}
\end{equation}
\end{comment}
\begin{equation*}
\begin{aligned}
\check{C} \leq &\, \begin{multlined}[t] C_{\textup{lin}}\left\{1+\frac{2\overline{k}}{\underline{c}^2} ( \tripplenorme{\u}
+   \tripplenorma{\psi})\mathcal{C}(T)\left[1 + h^{s-1-d/2} C_\star +  h^{s-1-d/2}
 \right]C_\star\right\}. \end{multlined}
\end{aligned}
\end{equation*}
Observe that if $h^{s-1-d/2} \leq 1$ and $ h^{s-1-d/2}C_\star \leq 1$, then 
\begin{comment}
	\begin{equation}
	\begin{aligned}
	\check{C} \leq &\, \begin{multlined}[t] C_{\textup{lin}}\left \{1 +\frac{2\overline{k}}{\underline{c}^2}  \delta(C_{\textup{emb}}(C_{\textup{st}}+1)(1+\sqrt{T})+ \Cinv +(1+C_{\textup{app}})\sqrt{T}\tilde{C}\Cinv +C_{\textup{app}{2}}\Cinv)   C_\star\right \}  . \end{multlined}
	\end{aligned}
	\end{equation}
\end{comment}
	\begin{equation*}
\begin{aligned}
\check{C} \leq &\, \begin{multlined}[t] C_{\textup{lin}}\left \{1 +\frac{2\overline{k}}{\underline{c}^2}  ( \tripplenorme{\u}
+   \tripplenorma{\psi}) \cdot 3\mathcal{C}(T)   C_\star\right \}\leq C_{\textup{lin}}\left(1+\frac{6\overline{k}}{\underline{c}^2}\delta \mathcal{C}(T)C_{\star}\right)  . \end{multlined}
\end{aligned}
\end{equation*}	
	In order to fulfill the desired property that $\check{C}\leq C_{\star}$, we demand the bound $\delta$ for $\tripplenorme{\u}
	+   \tripplenorma{\psi}$ to be small enough so that
	\begin{comment}	
	\begin{equation} \label{smallness_1}
		1-C_{\textup{lin}} \frac{2\overline{k}}{\underline{c}^2}\delta(C_{\textup{emb}}(C_{\textup{st}}+1)(1+\sqrt{T})+ \Cinv+(1+C_{\textup{app}})\sqrt{T}\tilde{C}\Cinv +C_{\textup{app}{2}}\Cinv +2C_{\textup{st}}) >0
	\end{equation}

{\color{lightgray}
	\[1-C_{\textup{lin}}\left(1 +\frac{2k}{c^2} \|(\vecc{u}, \psi)\|_{{X}}(3+ \sqrt{T}\tilde{C} +(1+\tilde{C}{\sqrt{T}})C_{\textup{app}}+2C_{\textup{st}})  \right)>0\]} 
\end{comment}
	\begin{equation} \label{smallness_1}
1-6 \frac{\overline{k}}{\underline{c}^2} C_{\textup{lin}} \mathcal{C}(T)\delta >0
\end{equation}
and then set
	\begin{equation} \label{const_Cstar}
	\begin{aligned}
		C_{\star}(\delta, T)=&\,\frac{C_{\textup{lin}}}{1-6 \dfrac{\overline{k}}{\underline{c}^2} C_{\textup{lin}} \mathcal{C}(T)\delta}.
		\end{aligned}
	\end{equation}
	\begin{comment}
{\color{lightgray}
	\begin{equation} \label{const_Cstar}
	\begin{aligned}
	&C_\star(\vecc{u}, \psi)\\
	=&\, \frac{C_{\textup{lin}}}{1-C_{\textup{lin}}\left(1 +\frac{2k}{\underline{c}^2} \|(\vecc{u}, \psi)\|_{{X}}(3+ \sqrt{T}\tilde{C} +(1+\tilde{C}\sqrt{T})C_{\textup{app}}+2C_{\textup{st}})  \right)}.
	\end{aligned}
	\end{equation}}
\end{comment}	
	Finally, we take $\overline{h}$ small enough so that
	\[\overline{h}^{s-1-d/2} \leq \min \left\{1, \frac{1}{C_\star(\delta, T)}\right\}.\]
	Altogether, for $1+d/2 < s \leq p+1$, we then have 
	\begin{equation*}
	\begin{aligned}
	\sup_{t \in (0,T)}\|(\vecc{u}-\vecc{u}_h, \psi-\psi_h)(t)\|_{E}  \leq \, C_\star(\delta, T)  h^{s-1} ( \tripplenorme{\u}
	+   \tripplenorma{\psi}).
	\end{aligned}
	\end{equation*}
%Note that above we have chosen $1$ as the bound of $\overline{h}^{s-1-d/2}$ for simplicity of exposition; the arguments follow analogously with $\overline{h}^{s-1-d/2} \leq \bar{C}$ for a given fixed positive constant $\bar{C}$.
\end{proof}

\noindent The smallness condition \eqref{smallness_1} is mitigated in practice by the fact that the factor $\frac{\overline{k}}{\underline{c}^2}$ is quite small; see Section~\ref{sec:NumericalSimulation} for typical values of the acoustic medium parameters. We next provide sufficient conditions under which $\mathcal{S}$ is strictly contractive.
\begin{proposition} \label{Prop:Contractivity}
	Let $ 1+ d/2<s \leq p+1$. Assume that the penalty parameter $\beta$ in \eqref{stabilization} is chosen sufficiently large according to Proposition~\ref{Prop:LinStability} and Theorem~\ref{Thm:LinError}. There exist $\overline{h}>0$ and $\delta>0$, such that for $0< h\leq \overline{h}$ and $\tripplenorme{\u} +   \tripplenorma{\psi}  \leq \delta$, the mapping $\mathcal{S}$ is strictly contractive on $\mathcal{B}_h$ in the topology induced by the $L^\infty \textup{E}$ norm; \change{blue}{cf. \eqref{full_norms}}.
\end{proposition}
\begin{proof}
	To prove contractivity, take $(\tilde{\vecc{u}}^{(1)}_h, w^{(1)}_h)$, $( \tilde{\vecc{u}}^{(2)}_h, w^{(2)}_h) \in \mathcal{B}_h$ and set  \[(\vecc{u}^{(1)}_h, \psi_h^{(1)})=\mathcal{S}(\tilde{\vecc{u}}^{(1)}_h, w^{(1)}_h), \qquad (\vecc{u}^{(2)}_h, \psi^{(2)}_h)=\mathcal{S}(\tilde{\vecc{u}}^{(2)}_h, w^{(2)}_h).\] The difference $(\bar{\vecc{u}}_h, \bar{\psi}_h)$, where $\bar{\psi}_h=\psi^{(1)}_h-\psi^{(2)}_h$ and  $\bar{\vecc{u}}_h=\vecc{u}^{(1)}_h-\vecc{u}^{(2)}_h$, then satisfies the weak form
	\begin{equation*} \label{weak_form_discrete_contractivity}
	\begin{aligned}
& \begin{multlined}[t]\producte{\varrhoe \ddot{\bar{\vecc{u}}}_h(t)}{\vecc{v}_h}+\producte{2 \varrhoe \zeta \dot{\bar{\vecc{u}}}_h(t)}{\vecc{v}_h}+\producte{\varrhoe \zeta^2 \bar{\vecc{u}}_h(t)}{\vecc{v}_h}+\forme(\bar{\vecc{u}}_h(t), \vecc{v}_h)\\[0.2cm]
+	\producta{c^{-2} \ddot{\bar{\psi}}_h(t)}{\phi_h}+\formah(\bar{\psi}_h(t)+\tfrac{b}{c^2}\dot{\bar{\psi}}_h(t), \phi_h)\\[0.2cm]+\Iea(\varrhof(\dot{\bar{\psi}}_h(t)+\tfrac{b}{c^2} \ddot{\bar{\psi}}_h(t)), \vecc{v}_h)-\Iae(  \phi_h, \dot{\bar{\vecc{u}}}_h(t))
\end{multlined}\\[0.2cm]
	=&\, \begin{multlined}[t] \producta{\fa_h(\dot{w}^{(1)}_h(t), \ddot{w}^{(1)}_h(t), \nabla w^{(1)}_h(t), \nabla \dot{w}^{(1)}_h(t))}{\, \phi_h}\\ - \producta{\fa_h(\dot{w}^{(2)}_h(t), \ddot{w}^{(2)}_h(t), \nabla w_h^{(2)}(t), \nabla \dot{w}^{(2)}_h(t))}{\, \phi_h} \end{multlined}
	\end{aligned}
	\end{equation*}
	for all $(\vecc{v}_h, \phi_h) \in \Veh \times \Vah$ a.e. in time, and supplemented with zero initial conditions. We can then rely on the linear stability result of Proposition~\ref{Prop:LinStability} with $\fe$ and the initial conditions set to zero, and the right-hand side taken as the above difference of the $\fah$ terms. This immediately yields
	\begin{equation*}
	\begin{aligned}
	& \|\mathcal{S}(\tilde{\vecc{u}}^{(1)}_h, w^{(1)}_h)(t)-\mathcal{S}( \tilde{\vecc{u}}^{(2)}_h, w^{(2)}_h)(t)\|_{\textup{E}}\\
	  \lesssim&\,\begin{multlined}[t] \frac{1}{\underline{c}^2}\|2k_1 (\dot{w}^{(1)}_h\ddot{w}^{(1)}_h-\dot{w}^{(2)}_h\ddot{w}^{(2)}_h)+2k_2(\nabla w^{(1)}_h \cdot \nabla \dot{w}^{(1)}_h-\nabla w^{(2)}_h \cdot \nabla \dot{w}^{(2)}_h)\|_{L^2(0,t; L^2)}.\end{multlined}
	\end{aligned}
	\end{equation*}
	for all $t \in (0,T)$. Then analogously to deriving \eqref{est_f-fh}, we have the estimate 
	\begin{equation*}
	\begin{aligned}
			& \|\mathcal{S}(\tilde{\vecc{u}}^{(1)}_h, w^{(1)}_h)(t)-\mathcal{S}( \tilde{\vecc{u}}^{(2)}_h, w^{(2)}_h)(t)\|_{\textup{E}}\\
		 \lesssim&\, \begin{multlined}[t]
		 \|\dot{w}_h^{(1)}\|_{L^{\infty}L^{\infty}}\|\ddot{w}_h^{(1)}-\ddot{w}_h^{(2)}\|_{L^2L^2}+\|\ddot{w}_h^{(2)}\|_{L^2L^{\infty}}\|\dot{w}_h^{(1)}-\dot{w}_h^{(2)}\|_{L^{\infty}L^2}\\
		+\|\nabla\dot{w}_h^{(1)}\|_{L^2L^{\infty}}\|\nabla w_h^{(1)}-\nabla w_h^{(2)}\|_{L^{\infty}L^2}+\|\nabla{w}_h^{(2)}\|_{L^2L^{\infty}}\|\nabla \dot{w}_h^{(1)}-\nabla \dot{w}_h^{(2)}\|_{L^{\infty}L^2}. \end{multlined}
	\end{aligned}
	\end{equation*}	
\begin{comment}
	{\color{lightgray}\begin{equation} \label{final_contractivity}
	\begin{aligned}
	& \|\mathcal{S}(\vecc{u}_h, w^{(1)}_h)(t)-\mathcal{S}(\vecc{u}_h, w^{(2)}_h)(t)\|_{\textup{E}}\\
	\lesssim&\,
	\begin{multlined}[t] \left(\|\dot{w}^{(1)}_h\|_{L^\infty L^\infty}+\|\ddot{w}^{(2)}_h\|_{L^2 L^\infty}\right. \left.
	+\|\nabla w^{(2)}_h\|_{L^2L^\infty} +\|\nabla \dot{w}^{(1)}_h\|_{L^2 L^\infty} \right)\\
	\times (\|\ddot{w}^{(1)}_h-\ddot{w}^{(2)}_h\|_{L^2L^2}+\|\dot{w}^{(1)}_h-\dot{w}^{(2)}_h\|_{L^\infty L^2}+\|\nabla(\dot{{w}}^{(1)}_h-\dot{w}^{(2)}_h)\|_{L^\infty L^2}\\+\|\nabla(w^{(1)}_h-w^{(2)}_h)\|_{L^\infty L^2}). \end{multlined}
	\end{aligned}
	\end{equation}}
\end{comment}	
	By taking the supremum over $t \in (0,T)$, we obtain
	\begin{equation} \label{final_contractivity}
	\begin{aligned}
	& \sup_{t \in (0,T)}\|\mathcal{S}(\tilde{\vecc{u}}^{(1)}_h, w^{(1)}_h)(t)-\mathcal{S}( \tilde{\vecc{u}}^{(2)}_h, w^{(2)}_h)(t)\|_{\textup{E}}\\
	\lesssim&\,
	\begin{multlined}[t] \left(\|\dot{w}^{(1)}_h\|_{L^\infty L^\infty}+\|\ddot{w}^{(2)}_h\|_{L^2 L^\infty}\right. \left.
	+\|\nabla w^{(2)}_h\|_{L^2L^\infty} +\|\nabla \dot{w}^{(1)}_h\|_{L^2 L^\infty} \right)\\
	\times \sup_{t \in (0,T)} \|(\tilde{\vecc{u}}^{(1)}_h-\tilde{\vecc{u}}^{(2)}_h, w^{(1)}_h-w^{(2)}_h)\|_{\textup{E}}, \end{multlined}
	\end{aligned}
	\end{equation}
	where we have again relied on estimate \eqref{bound_grad}. In view of estimates  \eqref{est_nabla_phi_h}--\eqref{est_ddotpsih} and analogous ones that can be derived for $\|\dot{w}^{(1)}_h\|_{L^\infty L^\infty}$ and $\|\nabla \dot{w}^{(1)}_h\|_{L^2L^\infty}$, we can reduce the terms in the bracket on the right-hand side of \eqref{final_contractivity} by assuming smallness of $\tripplenorme{\u} +   \tripplenorma{\psi}$ and $h$. In this way, we obtain strict contractivity of the mapping $\mathcal{S}$ in the $L^\infty \textup{E}$ norm, as claimed.
\end{proof}

Similarly to before, the hidden constant in \eqref{final_contractivity} has the form $\frac{\overline{k}}{\underline{c}^2} \cdot C$, which helps to fulfill the strict contractivity condition in more realistic ultrasonic settings.  By virtue of the previous two results and \change{orange}{the Banach fixed-point theorem}, we obtain a unique approximate solution $(\vecc{u}_h, \psi_h)$ in $\mathcal{B}_h$.
\begin{theorem} \label{Thm:Main}
Under the assumptions of Propositions~\ref{Prop:SelfMapping} and~\ref{Prop:Contractivity}, there exist $\overline{h}>0$ and $\delta>0$, such that for $0<h\leq \overline{h}$ and $\tripplenorme{\u} +   \tripplenorma{\psi} \leq \delta$, approximate solution $(\vecc{u}_h, \psi_h)$ of the nonlinear elasto-acoustic problem \eqref{approx_IC} satisfies the following error bound:
\begin{equation*}
\begin{aligned}
\change{dkgreen}{\|(\u-\tilde{\vecc{u}}_h,\psi-w_h)\|_{L^\infty E}}  \leq \, C_\star(\delta, T)  h^{s-1} ( \tripplenorme{\u}+   \tripplenorma{\psi}),
\end{aligned}
\end{equation*}
where the constant $C_\star(\delta, T)$ is given in \eqref{const_Cstar}. 	
\end{theorem}

\section{Numerical simulation}\label{sec:NumericalSimulation}
In this chapter, we perform numerical simulations with the method proposed and analyzed before. We begin with synthetic experiments to back up the proven convergence results with numerical data. Therefore, we conduct several mesh analysis scenarios, where we test the numerical solution against a known analytical one. Later we will come back to the initial motivation for this work and use the proposed method to simulate ultrasound excitation, propagation, and transition into human tissue in a more natural setting with physical parameters and more realistic domains. \change{dkgreen}{The simulations are conducted with the software \textit{SPEED}, \cite{mazzieri2013speed, antonietti2012non, stupazzini2009near}, in detail its elasto-acoustic development branch \cite{antonietti2020high}.} \version{}{\change{orange}{In order to focus on the numerical results, for precise details e.g. about the used material parameter values, artificial solutions, boundary data and geometry measures we refer to the extended arxiv version \cite{arxivVersion}.}}

\subsection{Numerical experiment 1: Test against analytical solution} 
\begin{figure}[h!]
	\begin{center}
		\includegraphics[scale=0.675]{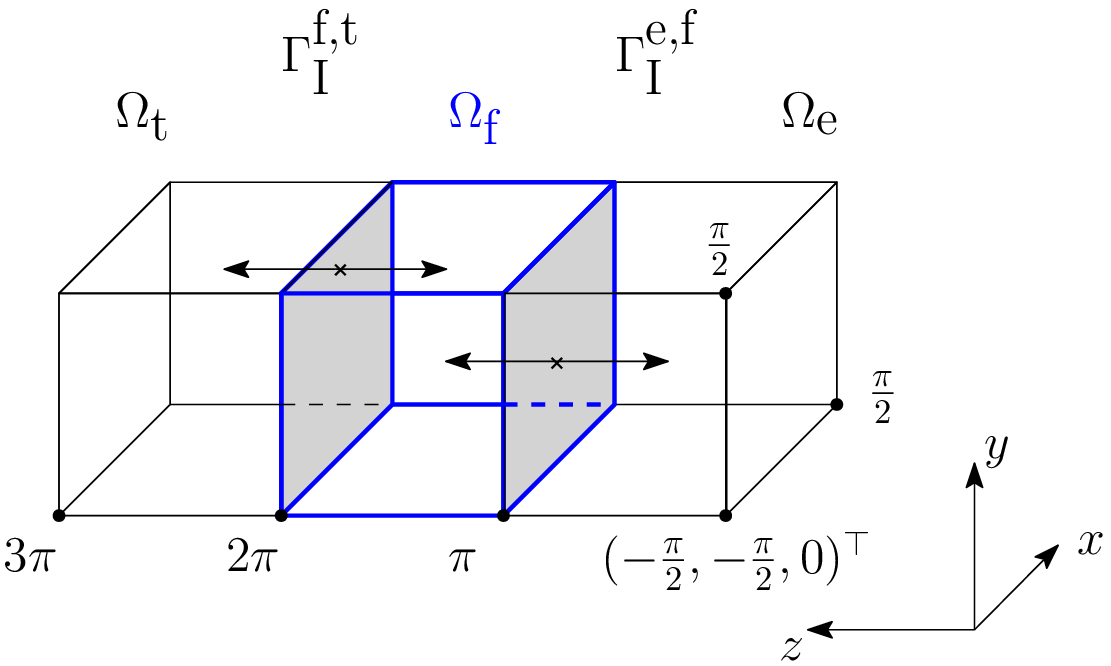}\hspace*{1.5cm}\raisebox{8mm}{\includegraphics[scale=0.125]{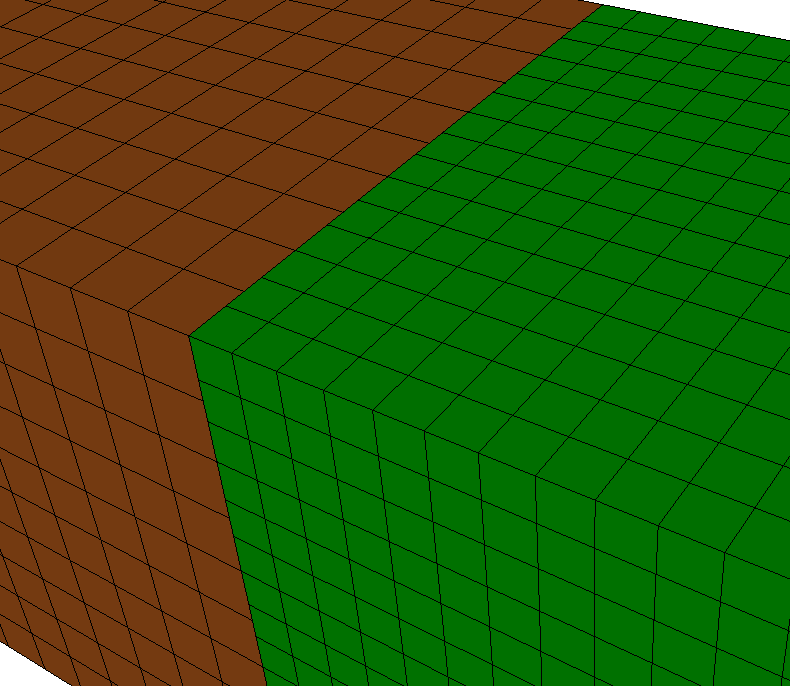}}
	\end{center}~\\[-15mm]
	\caption{\label{fig:CubesDomain} \textbf{(left)} ``\textit{Three stacked cubes}'' artificial domains for elasto-acoustic coupled problem with analytical solution available. \change{blue}{\textbf{(right)} Exemplary non-conforming mesh at the interfaces between the subdomains.}}
\end{figure}

In this first numerical experiment, we employ artificial domain sizes, boundary conditions, and external forces in order to enforce a known analytical solution to the full coupled problem. We then conduct a convergence study comparing our numerical solutions against that analytical one.\\
\indent We consider the following domains: $\Omegae:=(-\frac{\pi}{2},\frac{\pi}{2})^2\times(0,\pi)$, $\Omegaf:=(-\frac{\pi}{2},\frac{\pi}{2})^2\times(\pi,2\pi)$, and $\Omegat:=(-\frac{\pi}{2},\frac{\pi}{2})^2\times(2\pi,3\pi)$, which are depicted in Figure \ref{fig:CubesDomain}.
~\\[-10mm]

\paragraph{\bf Option 1: Elastic tissue} As the first option, we choose the tissue domain $\Omega_{\textup{t}}$ to be an elastic material. On the three domains, we then define the following analytical solutions:
\change{dkgreen}{\begin{alignat*}{2}
	\u^{\textup{i}}(x,y,z,t)&=\begin{pmatrix}
	\sin(x)\cos(y)\sin(z)\\
	\cos(x)\sin(y)\sin(z)\\
	\cos(x)\cos(y)\cos(z)
	\end{pmatrix}a\e(t)&&\qquad \textup{in } \Omega_{\textup{i}}\times (0,T],\quad \textup{i}\in\lbrace\textup{e,t}\rbrace\\
	\psi(x,y,z,t)&=\hspace*{4mm}\cos(x)\cos(y)\sin(z)\hspace*{3mm}a^{\textup{a}}(t)&&\qquad \textup{in } \Omegaf\times (0,T]
\end{alignat*}}
where the time dependent amplitudes of the elastic and acoustic fields $a\e(t)$ and $a^{\textup{a}}(t)$ are given by \change{dkgreen}{\[a^{\textup{i}}(t)=E^{\textup{i}}\sin(t)+D^{\textup{i}}\cos(t), \quad E^{\textup{i}},D^{\textup{i}}\in \R,\quad  \textup{i}\in\lbrace\textup{e,a}\rbrace.\]} \version{Indeed, those functions and vector fields pose an analytical solution to the coupled problem on the given domains, provided that
\begin{align*}
	-\varrhof &= \lambda\e-2\mu\e=\lambda\t-2\mu\t,\\
	D^{\textup{e}}&=E^{\textup{a}}-D^{\textup{a}}\frac{b\f}{(c\f)^2},\qquad E^{\textup{e}}=-D^{\textup{a}}-E^{\textup{a}}\frac{b\f}{(c\f)^2}
\end{align*}
which can be enforced by an adequate choice of synthetic parameters, and given the following right-hand side forcing terms for $\textup{i}\in\lbrace \textup{e,t}\rbrace$:
{\small
\change{orange}{
\begin{align*}
	\boldsymbol{f}^{\textup{i}}=&\left\{\left[2\varrho^{\textup{i}}\zeta^{\textup{i}}\dot{a}^{\textup{e}}(t)+\varrho^{\textup{i}}((\zeta^{\textup{i}})^2-1)a^{\textup{e}}(t)\right]\mathds{1}+\begin{pmatrix}
	\lambda^{\textup{i}}+4\mu^{\textup{i}}&  &  \\
	 & \lambda^{\textup{i}}+4\mu^{\textup{i}} &  \\
	 &  & 2\mu^{\textup{i}}-\lambda^{\textup{i}}
	\end{pmatrix}a^{\textup{e}}(t)\right\}\begin{pmatrix}
	\sin(x)\cos(y)\sin(z)\\
	\cos(x)\sin(y)\sin(z)\\
	\cos(x)\cos(y)\cos(z)
	\end{pmatrix}
\end{align*}	
}
\begin{align*}
%	\boldsymbol{f}^{\textup{i}}=&\begin{pmatrix}
%	\sin(x)\cos(y)\sin(z)\\
%	\cos(x)\sin(y)\sin(z)\\
%	\cos(x)\cos(y)\cos(z)
%	\end{pmatrix}\hspace*{-1mm}\left\{\left[\varrho^{\textup{i}}((\zeta^{\textup{i}})^2-1)+\hspace*{-1mm}\left.\begin{cases}
%	\lambda^{\textup{i}}+4\mu^{\textup{i}},&\hspace*{-2mm}x-\textup{coord.}\\
%	\lambda^{\textup{i}}+4\mu^{\textup{i}},&\hspace*{-2mm}y-\textup{coord.}\\
%	2\mu^{\textup{i}}-\lambda^{\textup{i}},&\hspace*{-2mm}z-\textup{coord.}\\
%	\end{cases}\right\}\right]a\e(t)+\left[2\varrho^{\textup{i}}\zeta^{\textup{i}}\right]\dot{a}\e(t)\right\},\\
	\fa&=\ff=\frac{1}{c^2}\bigg\{\left[(3c^2-1)a^{\textup{a}}(t)+3b\dot{a}^{\textup{a}}(t)\right]\cos(x)\cos(y)\sin(z)\\
	&\qquad\qquad\qquad +2\left[k_1\cos^2(x)\cos^2(y)\sin^2(z)-k_2\left(\sin^2(x)\cos^2(y)\sin^2(z)\right.\right.\\
	&\qquad\qquad\qquad\quad \left.\left.+\cos^2(x)\sin^2(y)\sin^2(z)+\cos^2(x)\cos^2(y)\cos^2(z)\right)\right]{a}^{\textup{a}}(t)\dot{a}^{\textup{a}}(t)\bigg\}
%	\fa&=\ff=\frac{1}{c^2}\bigg\{s^{\textup{a}}\cos(x)\cos(y)\sin(z)\left[(3c^2-1)\sin(t)+3b\cos(t)\right]\\
%	&\qquad \qquad -2(s^{\textup{a}})^2\left[k_2\sin^2(z)\left(\sin^2(x)\cos^2(y)+\cos^2(x)\sin^2(y)\right)\right.\\
%	&\qquad \qquad +\left.\cos^2(x)\cos^2(y)\left(k_2\cos^2(z)-k_1\sin^2(z)\right)\right]\cos(t)\sin(t)\bigg\}
\end{align*}}
\noindent The boundary and initial conditions are taken from the analytical solutions. The material parameters that are used can be found in Table \ref{tab:ArtMatParam}, where for Option 1 in the column of $\Omega_{\textup{t}}$ only the elastic values are relevant.}{In order for these solutions to fulfill the coupling conditions on the elasto-acoustic interfaces, certain conditions on the choice of the material parameters and the prescribed solution amplitudes have to be considered. The right hand sides $f^{\textup{i}},\textup{i}\in\lbrace\textup{e,t,f}\rbrace$ enforcing the given solutions can be obtained by insertion into the PDEs.}\\
\indent To analyze convergence of the proposed numerical scheme, we conduct the simulation on a sequence of meshes with $h$ tending to zero. Time integration is always conducted with a final time $T=2\pi$ and small enough time step size, such that on all meshes in use, the overall errors are all dominated by their spatial components. For the time integration, we employ the Newmark scheme in its predictor-corrector form using $\beta=0.25$ and $\gamma=0.5$ as in \citep[\S 5]{MKaltenbacher} for the nonlinear acoustic components, while the leapfrog scheme is used for the elastic ones. \version{Figure \ref{fig:ElasticTissueSimulation} \change{dkgreen}{(left)} shows the numerical solution on a selected mesh at time $t=0.625\cdot 2\pi$.}{} \change{blue}{The meshes are chosen to be \textit{non-conforming} at the subdomain interfaces \textit{on purpose} for this artificial data experiment to show convergence also in this general situation. In order to guarantee the non-conformance, the whole mesh-sequence is generated with a 3:2 ratio in mesh size between the central and the two outer blocks c.f. Fig.\,\ref{fig:CubesDomain}.}
\version{
\begin{figure}[h!]
	\begin{center}
		\includegraphics[trim=1mm 1mm 1mm 1mm, clip, scale=0.12]{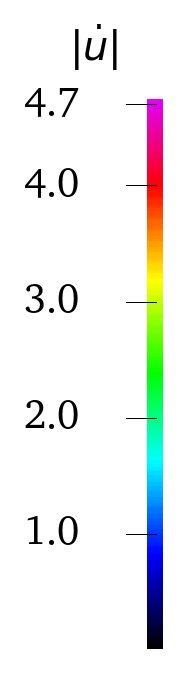}~~\includegraphics[scale=0.12]{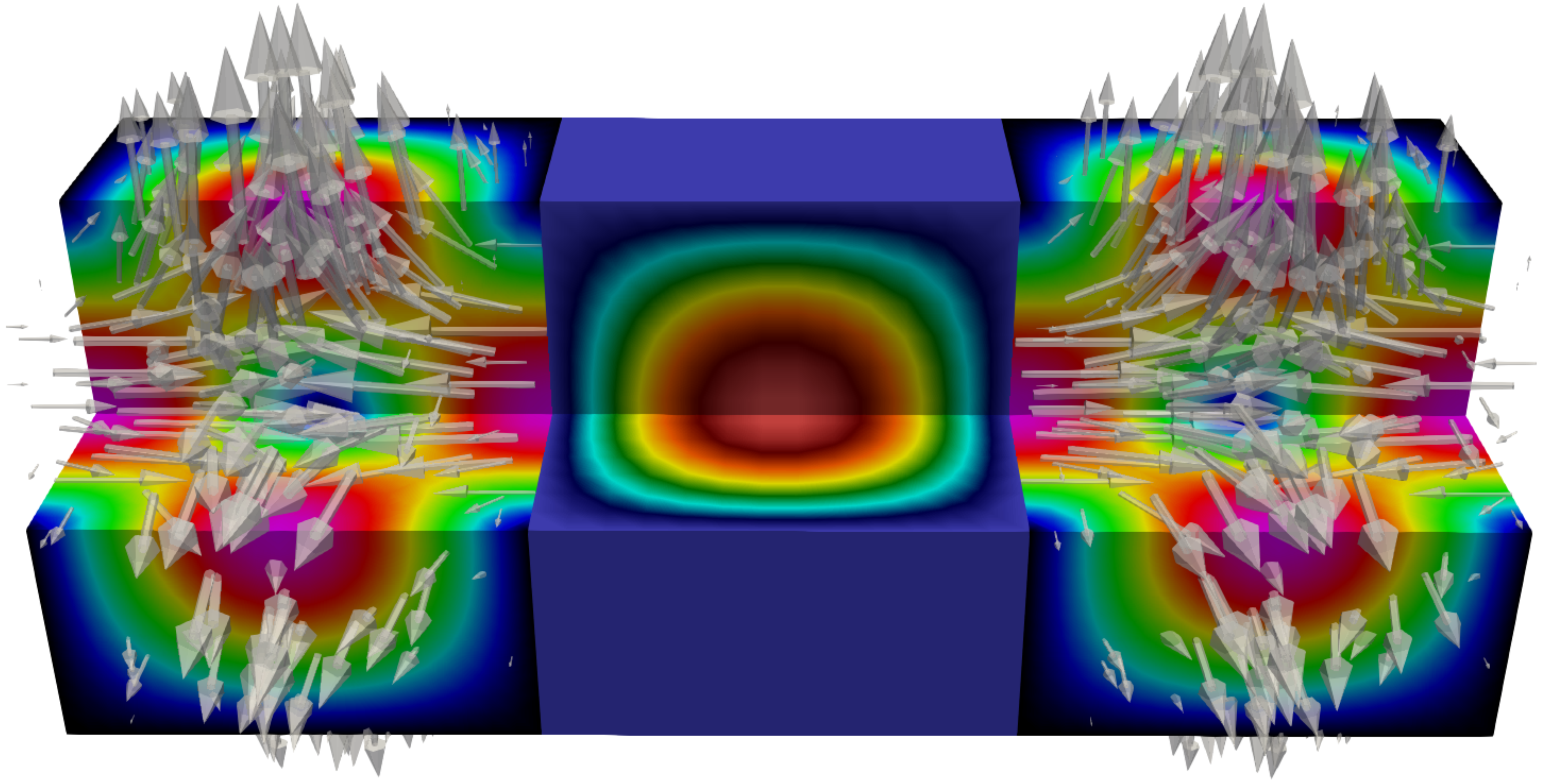}\hspace*{0cm}\includegraphics[trim=1mm 1mm 1mm 1mm, clip, scale=0.12]{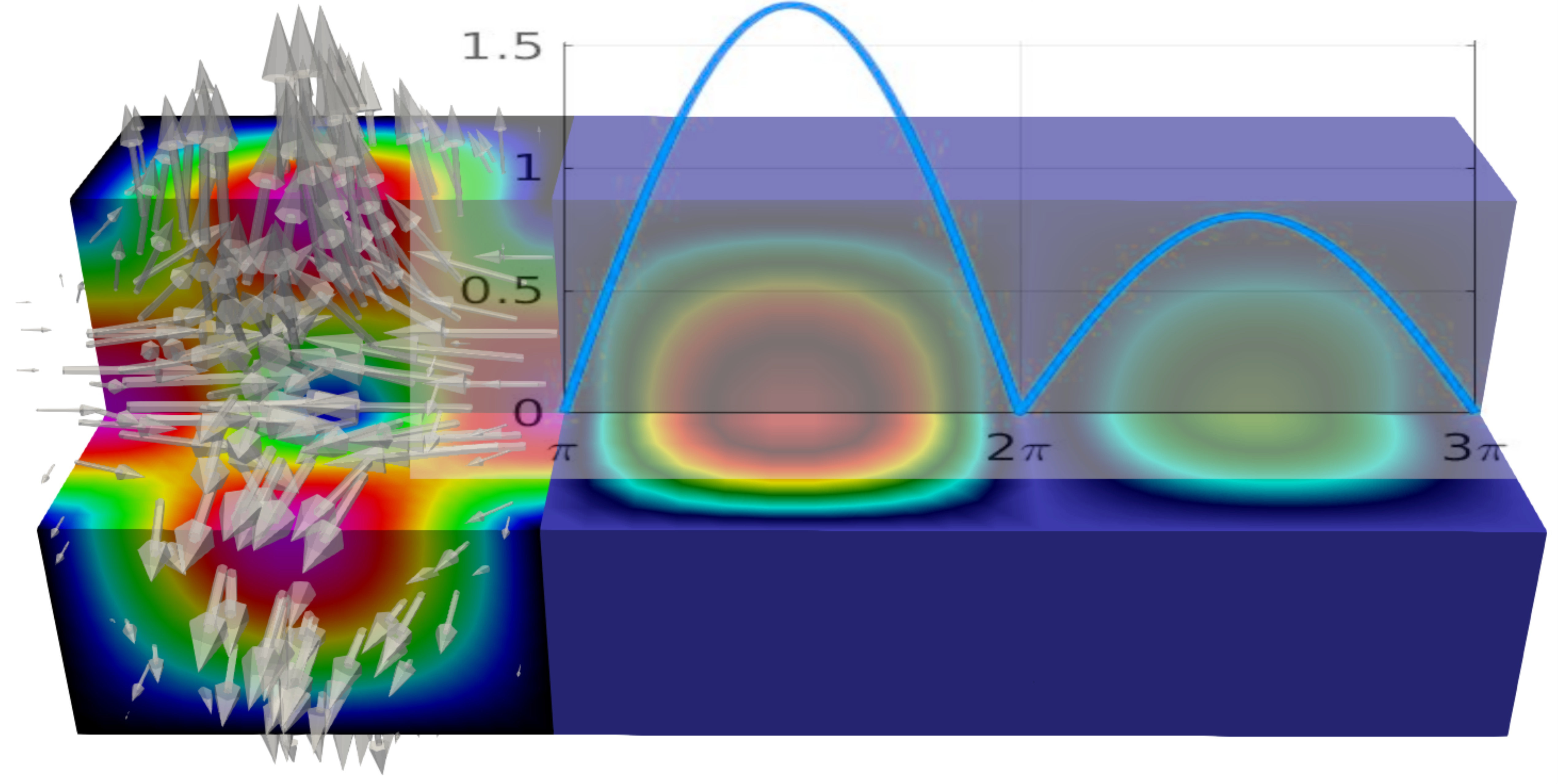}~~\includegraphics[scale=0.12]{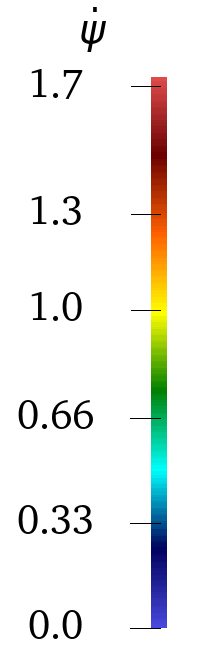}\\
	\end{center}
	\caption{\label{fig:ElasticTissueSimulation} \change{dkgreen}{Numerical solution using Material set 1 with \textbf{(left)} Option 1; i.e., with elastic tissue domain on the right. The middle domain is the only acoustic one here displaying $\dot{\psi}_h$. In the left and right elastic domains the vector field $\dot{\u}_h$ is shown in addition to the magnitude $\|\dot{\u}_h\|$. \textbf{(right)} Option 2; i.e., with acoustic tissue domain on the right.}}
\end{figure}~\\}{~\\}
\indent Convergence results are collectively plotted at the end of this subsection in order to be directly compared with modeling Option 2.

\paragraph{\bf Option 2: Acoustic tissue} We use the same geometry and subdomains as before (cf. Figure~\ref{fig:CubesDomain}), only with the tissue domain being modeled as an acoustic medium this time. This corresponds to the case analyzed in the theoretical part of this work. The analytical solutions are defined as before \change{dkgreen}{just with $\u$ (dropping the index \textup{e} for the elastic displacement $\u$) only being defined on $\Omegae$, while $\psi$ extends to $\Omegaa$ now.}
%\begin{alignat*}{2}
%\u(x,y,z,t)&=\begin{pmatrix}
%\sin(x)\cos(y)\sin(z)\\
%\cos(x)\sin(y)\sin(z)\\
%\cos(x)\cos(y)\cos(z)
%\end{pmatrix}a\e(t)&&\qquad \textup{in } \Omegae\times (0,T],\\
%\psi(x,y,z,t)&=\hspace*{4mm}\cos(x)\cos(y)\sin(z)\hspace*{3mm}a^{\textup{a}}(t)&&\qquad \textup{in } \Omegaa\times (0,T],
%\end{alignat*}
The amplitudes are again scaled with $a^{\textup{a}}(t)=E^{\textup{a}}\sin(t)+D^{\textup{a}}\cos(t)$ and $a^{\textup{e}}(t)=E^{\textup{e}}\sin(t)+D^{\textup{e}}\cos(t)$. However, this time the acoustic amplitudes are allowed to attain different values in the two individual acoustic domains; i.e.,
\begin{equation*}
	E^{\textup{a}}=\begin{cases}
	E^{\textup{f}}, & \textup{in } \Omegaf\\
	E^{\textup{t}}, & \textup{in } \Omegat
	\end{cases},\qquad D^{\textup{a}}=\begin{cases}
	D^{\textup{f}}, & \textup{in } \Omegaf\\
	D^{\textup{t}}, & \textup{in } \Omegat.
	\end{cases}
\end{equation*}
In order for these functions to be indeed solutions to the coupled problem, we employ the same conditions and right-hand side forcing terms as before, but in addition, to satisfy the acoustic-acoustic interface flux condition, \version{we demand that
\begin{align*}
	E\f-\frac{b\f}{(c\f)^2}D\f &= E\t-\frac{b\t}{(c\t)^2}D\t,\\
	D\f+\frac{b\f}{(c\f)^2}E\f &= D\t+\frac{b\t}{(c\t)^2}E\t.
\end{align*}
These conditions can be fulfilled, for example, by choosing $D\t=(c\t)^2$, $D\f=(c\f)^2$,  \[E\t=((c\f)^4(c\t)^2-(c\f)^2(c\t)^4+(b\f)^2(c\t)^2-b\f b\t(c\t)^2)/((c\f)^2b\t-(c\t)^2b\f),\]
and  $E\f=E\t-b\t+b\f$.}{additional constraints on the material parameters are imposed.} It should be noted that due to $\Omegaf$ and $\Omegat$ having \textit{different} material parameters, the forcing term $\fa$ also differs/jumps in between the acoustic subdomains. \version{The concrete material parameters that are used can be found again in Table \ref{tab:ArtMatParam}, where in the column $\Omegat$ only the acoustic entries are relevant this time.}{}\version{Once more, Figure \change{dkgreen}{\ref{fig:ElasticTissueSimulation} (right)} shows the numerical solution for the given situation at time $t=0.625\cdot 2\pi$ on a selected mesh, while convergence results can be found in Figure \ref{fig:ConvergenceStudy}}{Convergence results can be found in Figure \ref{fig:ConvergenceStudy}}. We note that the simulation with polynomial degree $p=1$ is not included in our theory. However the numerical experiment suggests that, at least in this synthetic setting, the convergence rate $\mathcal{O}(h^{s-1})$ can also be obtained for $p=1$.
%\begin{figure}[h!]
%	\begin{center}
%		\hspace*{1.25cm}\includegraphics[scale=0.22]{Acoust_Overlay.pdf}
%	\end{center}
%	\caption{\label{fig:AcousticTissueSimulation} Numerical solution using Material set 1 with Option 2 i.e. with acoustic tissue domain (back/right). The blue curve shows $\dot{\psi}$ across the central axis of the geometry in order to visualize the kink at the acoustic-acoustic interface.}
%\end{figure}
\version{
\begin{table}[h!]\renewcommand{\arraystretch}{1.4}
	\normalsize\setlength\tabcolsep{10pt}
	\begin{tabular}{|c|c||c|c|c||c|}
	\cline{3-5}
	\multicolumn{2}{c||}{\textbf{Material Set} 1/{\color{blue}\textit{(2)}}}&\multicolumn{3}{|c||}{\textbf{material block}} & \multicolumn{1}{c}{}\\ \hline
	\textbf{type} & \textbf{parameter} & \textbf{artificial $\Omegat$} & \textbf{artificial $\Omegaf$} & \textbf{artificial $\Omegae$} & \textbf{unit}\\ \hline\hline
	all& $\varrho$ & 1.6018 {/\color{blue}\textit{(2)}}& 0.51{/\color{blue}\textit{(2)}} & 1{/\color{blue}\textit{(6)}} & $\frac{\textup{kg}}{\textup{m}^3}$\\ \hline
	& $\lambda$   & 12.3041{/\color{blue}\textit{(6)}} & / & 2.87{/\color{blue}\textit{(10)}} & $\frac{\textup{N}}{\textup{m}^2}$\\
	elastic& $\mu$		   & 6.4070{/\color{blue}\textit{(4)}} & / & 1.69{/\color{blue}\textit{(6)}} &  $\frac{\textup{N}}{\textup{m}^2}$\\ 
	&$\zeta$       & 6{/\color{blue}\textit{(2)}} & / & 3{/\color{blue}\textit{(4)}}    & $\frac{1}{\textup{s}}$\\ \hline
	& $c$		   & 2{/\color{blue}\textit{(1)}} & 1{/\color{blue}\textit{($\sqrt{2}$)}} & / & $\frac{\textup{m}}{\textup{s}}$\\
	& $b$  & $\frac{1}{2}${/\change{dkgreen}{\color{blue}\textit{$(\frac{11}{2})$}}} & 1{/\color{blue}\textit{(1)}} & / &$\frac{\textup{m}^2}{\textup{s}}$\\
	acoustic& $k_1$ 	   & $\frac{1}{16}${/\color{blue}$(\frac{1}{4})$} & $\frac{1}{10}${/\color{blue}$(\frac{1}{20})$} & / &$\frac{\textup{s}^2}{\textup{m}^2}$\\
	& $k_2$ 	   & 0{/\color{blue}\textit{(1)}} & 0{/\color{blue}$(\frac{1}{4})$} & / &1\\
	& $E^{\textup{a}}$ 	       & 2.8571{/\change{dkgreen}{\color{blue}\textit{(-0.25)}}} & 3.3571{/\color{blue}\textit{(-4.75)}} & / &1\\
	& $D^{\textup{a}}$ 	       & 4{/\color{blue}\textit{(1)}} & 1{/\color{blue}\textit{(2)}} & / &1\\ \hline
\end{tabular}
\centering	\caption{\textbf{Synthetic material parameters.} Black and upright numbers form Parameter set 1, blue and italic parameters in parenthesis form set 2. In each case, parameters are chosen in such a way that a given analytical solution is enforced in the present numerical experiment for a convergence study. For $\Omegat$, elastic or acoustic parameters are chosen, depending on Option 1 or Option 2 being used as a model for the tissue. (Numbers rounded to 4 digits)\label{tab:ArtMatParam}}
\end{table}}{~}

\begin{remark}
Note that for Parameter set 1, the factor $k_2$ in front of the quadratic gradient type nonlinearity is set to zero, reducing the model to Westervelt's equation in potential form (with synthetic material parameters) in the acoustic sub-domains. For Parameter set 2, the parameter $k_2$ is non-zero and hence inclusion of the quadratic gradient nonlinearity leads to Kuznetsov's equation of nonlinear acoustics~\citep{kuznetsov1971equation}, with synthetic material parameters as well.
\end{remark}

\version{\begin{figure}[h!]	
	\noindent{\hspace*{0cm}{\bf Material set 1 with Westervelt-type nonlinearity}}\\[-10cm]
		\hspace*{0cm}% This file was created by matlab2tikz.
\definecolor{mycolor1}{rgb}{0.00000,0.50000,1.00000}%
\definecolor{mycolor2}{rgb}{0.30000,0.90000,0.20000}%
\definecolor{mycolor3}{rgb}{1.00000,0.60000,0.00000}%
\begin{tikzpicture}

\begin{axis}[%
scale=0.75,
width=2.994in,
height=2.846in,
at={(0.0in,0.926in)},
scale only axis,
xmode=log,
xmin=0.310937033868723,
xmax=3.40087380793916,
xminorticks=true,
xlabel style={font=\color{white!15!black}},
xlabel={$h$},
ymode=log,
ymin=1e-05,
ymax=1.5,
yminorticks=true,
ylabel style={font=\color{white!15!black}},
ylabel={},%$\|e\|_{L^{\infty}\textup{E}}/\|(u,\psi)_{\textup{ex}}\|_{L^{\infty}\textup{E}}$},
axis background/.style={fill=white},
title style={font=\bfseries},
title={Option 1},
xmajorgrids,
xminorgrids,
ymajorgrids,
yminorgrids,
legend style={at={(1.03,0.5)}, anchor=west, legend cell align=left, align=left, draw=white!15!black}
]
\addplot [color=mycolor1, line width=2.0pt, mark=diamond, mark options={solid, mycolor1}]
  table[row sep=crcr]{%
2.72069904635133	0.512524115102187\\
1.36034952317566	0.253815045137582\\
0.906899682117109	0.16839689312084\\
0.680174761587832	0.126035026791843\\
0.544139809270265	0.10072168064436\\
0.453449841058554	0.0838841405837072\\
0.388671292335904	0.0718737397960915\\
};
%\addlegendentry{$p=1$}

\addplot [color=black!34!mycolor1, dashed, line width=2.0pt]
  table[row sep=crcr]{%
2.72069904635133	0.244862914171619\\
1.36034952317566	0.12243145708581\\
0.906899682117109	0.0816209713905398\\
0.680174761587832	0.0612157285429049\\
0.544139809270265	0.0489725828343238\\
0.453449841058554	0.0408104856952699\\
0.388671292335904	0.0349804163102314\\
};
%\addlegendentry{$C\cdot h$}

\addplot [color=mycolor2, line width=2.0pt, mark=diamond, mark options={solid, mycolor2}]
  table[row sep=crcr]{%
2.72069904635133	0.077490231476529\\
1.36034952317566	0.0190126226804041\\
0.906899682117109	0.00839710418272116\\
0.680174761587832	0.00471133312745843\\
0.544139809270265	0.00301369181063691\\
0.453449841058554	0.00209538378297584\\
0.388671292335904	0.00154426986916315\\
};
%\addlegendentry{$p=2$}

\addplot [color=black!34!mycolor2, dashed, line width=2.0pt]
  table[row sep=crcr]{%
2.72069904635133	0.0370110165040851\\
1.36034952317566	0.00925275412602127\\
0.906899682117109	0.00411233516712057\\
0.680174761587832	0.00231318853150532\\
0.544139809270265	0.0014804406601634\\
0.453449841058554	0.00102808379178014\\
0.388671292335904	0.000755326867430309\\
};
%\addlegendentry{$C\cdot h^2$}

\addplot [color=mycolor3, line width=2.0pt, mark=diamond, mark options={solid, mycolor3}]
  table[row sep=crcr]{%
2.72069904635133	0.00840216107875206\\
1.36034952317566	0.00104750528373648\\
0.906899682117109	0.00031039083619995\\
0.680174761587832	0.0001320730065555\\
0.544139809270265	6.97741203513801e-05\\
0.453449841058554	4.36966851641945e-05\\
};
%\addlegendentry{$p=3$}

\addplot [color=black!34!mycolor3, dashed, line width=2.0pt]
  table[row sep=crcr]{%
2.72069904635133	0.0040278334922863\\
1.36034952317566	0.000503479186535787\\
0.906899682117109	0.000149179018232826\\
0.680174761587832	6.29348983169735e-05\\
0.544139809270265	3.22226679382903e-05\\
0.453449841058554	1.86473772791032e-05\\
};
%\addlegendentry{$C\cdot h^3$}

\end{axis}

\begin{axis}[%
width=9.825in,
height=7.368in,
at={(0in,0in)},
scale only axis,
xmin=0,
xmax=1,
ymin=0,
ymax=1,
axis line style={draw=none},
ticks=none,
axis x line*=bottom,
axis y line*=left,
legend style={legend cell align=left, align=left, draw=white!15!black}
]
\end{axis}
\end{tikzpicture}%
\hspace*{-20.5cm}% This file was created by matlab2tikz.
\definecolor{mycolor1}{rgb}{0.00000,0.50000,1.00000}%
\definecolor{mycolor2}{rgb}{0.30000,0.90000,0.20000}%
\definecolor{mycolor3}{rgb}{1.00000,0.60000,0.00000}%
\begin{tikzpicture}

\begin{axis}[%
scale=0.75,
width=2.994in,
height=2.846in,
at={(1.094in,0.93in)},
scale only axis,
xmode=log,
xmin=0.310937033868723,
xmax=3.40087380793916,
xminorticks=true,
xlabel style={font=\color{white!15!black}},
xlabel={$h$},
ymode=log,
ymin=1e-05,
ymax=1.5,
yminorticks=true,
ylabel style={font=\color{white!15!black}},
ylabel={},%$\|e\|_{L^{\infty}\textup{E}}/\|(u,\psi)_{\textup{ex}}\|_{L^{\infty}\textup{E}}$},
axis background/.style={fill=white},
title style={font=\bfseries},
title={Option 2},
xmajorgrids,
xminorgrids,
ymajorgrids,
yminorgrids,
legend style={at={(10.03,0.5)}, anchor=west, legend cell align=left, align=left, draw=white!15!black}
]
\addplot [color=mycolor1, line width=2.0pt, mark=diamond, mark options={solid, mycolor1}]
  table[row sep=crcr]{%
2.72069904635133	0.537566498291575\\
1.36034952317566	0.250509282906733\\
0.906899682117109	0.164651613321478\\
0.680174761587832	0.122875112819625\\
0.544139809270265	0.0980784636734661\\
0.453449841058554	0.081631124604748\\
0.388671292335904	0.0699171383576819\\
};
\addlegendentry{$p=1$}

\addplot [color=black!34!mycolor1, dashed, line width=2.0pt]
  table[row sep=crcr]{%
2.72069904635133	0.244862914171619\\
1.36034952317566	0.12243145708581\\
0.906899682117109	0.0816209713905398\\
0.680174761587832	0.0612157285429049\\
0.544139809270265	0.0489725828343238\\
0.453449841058554	0.0408104856952699\\
0.388671292335904	0.0349804163102314\\
};
\addlegendentry{$C\cdot h$}

\addplot [color=mycolor2, line width=2.0pt, mark=diamond, mark options={solid, mycolor2}]
  table[row sep=crcr]{%
2.72069904635133	0.0811283430603854\\
1.36034952317566	0.0194494902075759\\
0.906899682117109	0.00852705213671816\\
0.680174761587832	0.00476542503364779\\
0.544139809270265	0.00304010331180313\\
0.453449841058554	0.0021089386100639\\
0.388671292335904	0.00155066470424568\\
};
\addlegendentry{$p=2$}

\addplot [color=black!34!mycolor2, dashed, line width=2.0pt]
  table[row sep=crcr]{%
2.72069904635133	0.0370110165040851\\
1.36034952317566	0.00925275412602127\\
0.906899682117109	0.00411233516712057\\
0.680174761587832	0.00231318853150532\\
0.544139809270265	0.0014804406601634\\
0.453449841058554	0.00102808379178014\\
0.388671292335904	0.000755326867430309\\
};
\addlegendentry{$C\cdot h^2$}

\addplot [color=mycolor3, line width=2.0pt, mark=diamond, mark options={solid, mycolor3}]
  table[row sep=crcr]{%
2.72069904635133	0.00871803677695128\\
1.36034952317566	0.00108243111017819\\
0.906899682117109	0.000319539615898647\\
0.680174761587832	0.000135239035103899\\
0.544139809270265	7.06827232084058e-05\\
};
\addlegendentry{$p=3$}

\addplot [color=black!34!mycolor3, dashed, line width=2.0pt]
  table[row sep=crcr]{%
2.72069904635133	0.0040278334922863\\
1.36034952317566	0.000503479186535787\\
0.906899682117109	0.000149179018232826\\
0.680174761587832	6.29348983169735e-05\\
0.544139809270265	3.22226679382903e-05\\
};
\addlegendentry{$C\cdot h^3$}

\end{axis}

\begin{axis}[%
width=10in,
height=7.246in,
at={(0in,0in)},
scale only axis,
xmin=0,
xmax=1,
ymin=0,
ymax=1,
axis line style={draw=none},
ticks=none,
axis x line*=bottom,
axis y line*=left,
legend style={legend cell align=left, align=left, draw=white!15!black}
]
\end{axis}
\end{tikzpicture}%
\\[-1cm]
	{\hspace*{0cm}{\bf Material set 2 with Kuznetsov-type nonlinearity}}\\[-10cm]
		\hspace*{0cm}% This file was created by matlab2tikz.
\definecolor{mycolor1}{rgb}{0.00000,0.50000,1.00000}%
\definecolor{mycolor2}{rgb}{0.30000,0.90000,0.20000}%
\definecolor{mycolor3}{rgb}{1.00000,0.60000,0.00000}%
\begin{tikzpicture}

\begin{axis}[%
scale=0.75,
width=2.994in,
height=2.846in,
at={(0.0in,0.926in)},
scale only axis,
xmode=log,
xmin=0.310937033868723,
xmax=3.40087380793916,
xminorticks=true,
xlabel style={font=\color{white!15!black}},
xlabel={$h$},
ymode=log,
ymin=1e-05,
ymax=1.5,
yminorticks=true,
ylabel style={font=\color{white!15!black}},
ylabel={},%$\|e\|_{L^{\infty}\textup{E}}/\|(u,\psi)_{\textup{ex}}\|_{L^{\infty}\textup{E}}$},
axis background/.style={fill=white},
title style={font=\bfseries},
title={Option 1},
xmajorgrids,
xminorgrids,
ymajorgrids,
yminorgrids,
legend style={at={(1.03,0.5)}, anchor=west, legend cell align=left, align=left, draw=white!15!black}
]
\addplot [color=mycolor1, line width=2.0pt, mark=diamond, mark options={solid, mycolor1}]
  table[row sep=crcr]{%
1.36034952317566	0.267510558190737\\
0.906899682117109	0.170705330771081\\
0.680174761587832	0.126440659367848\\
0.544139809270265	0.100598625246998\\
0.453449841058554	0.0835866389956336\\
0.388671292335904	0.0715199119717162\\
};
%\addlegendentry{$p=1$}

\addplot [color=black!34!mycolor1, dashed, line width=2.0pt]
  table[row sep=crcr]{%
1.36034952317566	0.12243145708581\\
0.906899682117109	0.0816209713905398\\
0.680174761587832	0.0612157285429049\\
0.544139809270265	0.0489725828343238\\
0.453449841058554	0.0408104856952699\\
0.388671292335904	0.0349804163102314\\
};
%\addlegendentry{$C\cdot h$}

\addplot [color=mycolor2, line width=2.0pt, mark=diamond, mark options={solid, mycolor2}]
  table[row sep=crcr]{%
2.72069904635133	0.0982407126251855\\
1.36034952317566	0.019186634746431\\
0.906899682117109	0.00836761851441703\\
0.680174761587832	0.00468246035727835\\
0.544139809270265	0.00299858690224098\\
0.453449841058554	0.00209195285647537\\
0.388671292335904	0.00155156858938708\\
};
%\addlegendentry{$p=2$}

\addplot [color=black!34!mycolor2, dashed, line width=2.0pt]
  table[row sep=crcr]{%
2.72069904635133	0.0370110165040851\\
1.36034952317566	0.00925275412602127\\
0.906899682117109	0.00411233516712057\\
0.680174761587832	0.00231318853150532\\
0.544139809270265	0.0014804406601634\\
0.453449841058554	0.00102808379178014\\
0.388671292335904	0.000755326867430309\\
};
%\addlegendentry{$C\cdot h^2$}

\addplot [color=mycolor3, line width=2.0pt, mark=diamond, mark options={solid, mycolor3}]
  table[row sep=crcr]{%
2.72069904635133	0.00876182552080723\\
1.36034952317566	0.00103984427368147\\
0.906899682117109	0.000308225294019168\\
0.680174761587832	0.000133144188493677\\
0.544139809270265	7.38440784462287e-05\\
};
%\addlegendentry{$p=3$}

\addplot [color=black!34!mycolor3, dashed, line width=2.0pt]
  table[row sep=crcr]{%
2.72069904635133	0.0040278334922863\\
1.36034952317566	0.000503479186535787\\
0.906899682117109	0.000149179018232826\\
0.680174761587832	6.29348983169735e-05\\
0.544139809270265	3.22226679382903e-05\\
};
%\addlegendentry{$C\cdot h^3$}

\end{axis}

\begin{axis}[%
width=9.825in,
height=7.368in,
at={(0in,0in)},
scale only axis,
xmin=0,
xmax=1,
ymin=0,
ymax=1,
axis line style={draw=none},
ticks=none,
axis x line*=bottom,
axis y line*=left,
legend style={legend cell align=left, align=left, draw=white!15!black}
]
\end{axis}
\end{tikzpicture}%
\hspace*{-20.5cm}% This file was created by matlab2tikz.
\definecolor{mycolor1}{rgb}{0.00000,0.50000,1.00000}%
\definecolor{mycolor2}{rgb}{0.30000,0.90000,0.20000}%
\definecolor{mycolor3}{rgb}{1.00000,0.60000,0.00000}%
\begin{tikzpicture}

\begin{axis}[%
scale=0.75,
width=2.994in,
height=2.846in,
at={(1.094in,0.93in)},
scale only axis,
xmode=log,
xmin=0.310937033868723,
xmax=3.40087380793916,
xminorticks=true,
xlabel style={font=\color{white!15!black}},
xlabel={$h$},
ymode=log,
ymin=1e-05,
ymax=1.5,
yminorticks=true,
ylabel style={font=\color{white!15!black}},
ylabel={},%$\|e\|_{L^{\infty}\textup{E}}/\|(u,\psi)_{\textup{ex}}\|_{L^{\infty}\textup{E}}$},
axis background/.style={fill=white},
title style={font=\bfseries},
title={Option 2},
xmajorgrids,
xminorgrids,
ymajorgrids,
yminorgrids,
legend style={at={(10.03,0.5)}, anchor=west, legend cell align=left, align=left, draw=white!15!black}
]
\addplot [color=mycolor1, line width=2.0pt, mark=diamond, mark options={solid, mycolor1}]
  table[row sep=crcr]{%
1.36034952317566	0.282470698387141\\
0.906899682117109	0.179319451768731\\
0.680174761587832	0.132679065795229\\
0.544139809270265	0.10552901013792\\
0.453449841058554	0.0876733973220079\\
0.388671292335904	0.0750131826816834\\
};
\addlegendentry{$p=1$}

\addplot [color=black!34!mycolor1, dashed, line width=2.0pt]
  table[row sep=crcr]{%
1.36034952317566	0.12243145708581\\
0.906899682117109	0.0816209713905398\\
0.680174761587832	0.0612157285429049\\
0.544139809270265	0.0489725828343238\\
0.453449841058554	0.0408104856952699\\
0.388671292335904	0.0349804163102314\\
};
\addlegendentry{$C\cdot h$}

\addplot [color=mycolor2, line width=2.0pt, mark=diamond, mark options={solid, mycolor2}]
  table[row sep=crcr]{%
2.72069904635133	0.107598751849621\\
1.36034952317566	0.0206571765559414\\
0.906899682117109	0.00901710057175072\\
0.680174761587832	0.00504785532190627\\
0.544139809270265	0.00322865257526021\\
0.453449841058554	0.00224716492724014\\
0.388671292335904	0.00165959110150751\\
};
\addlegendentry{$p=2$}

\addplot [color=black!34!mycolor2, dashed, line width=2.0pt]
  table[row sep=crcr]{%
2.72069904635133	0.0370110165040851\\
1.36034952317566	0.00925275412602127\\
0.906899682117109	0.00411233516712057\\
0.680174761587832	0.00231318853150532\\
0.544139809270265	0.0014804406601634\\
0.453449841058554	0.00102808379178014\\
0.388671292335904	0.000755326867430309\\
};
\addlegendentry{$C\cdot h^2$}

\addplot [color=mycolor3, line width=2.0pt, mark=diamond, mark options={solid, mycolor3}]
  table[row sep=crcr]{%
2.72069904635133	0.00972718315972817\\
1.36034952317566	0.0011508459818235\\
0.906899682117109	0.000340631589547223\\
0.680174761587832	0.000145633977601148\\
};
\addlegendentry{$p=3$}

\addplot [color=black!34!mycolor3, dashed, line width=2.0pt]
  table[row sep=crcr]{%
2.72069904635133	0.0040278334922863\\
1.36034952317566	0.000503479186535787\\
0.906899682117109	0.000149179018232826\\
0.680174761587832	6.29348983169735e-05\\
};
\addlegendentry{$C\cdot h^3$}

\end{axis}

\begin{axis}[%
width=7in,
height=5.25in,
at={(0in,0in)},
scale only axis,
xmin=0,
xmax=1,
ymin=0,
ymax=1,
axis line style={draw=none},
ticks=none,
axis x line*=bottom,
axis y line*=left,
legend style={legend cell align=left, align=left, draw=white!15!black}
]
\end{axis}
\end{tikzpicture}%
\\[-1.5cm]
		\begin{center}
			\includegraphics[scale=0.15]{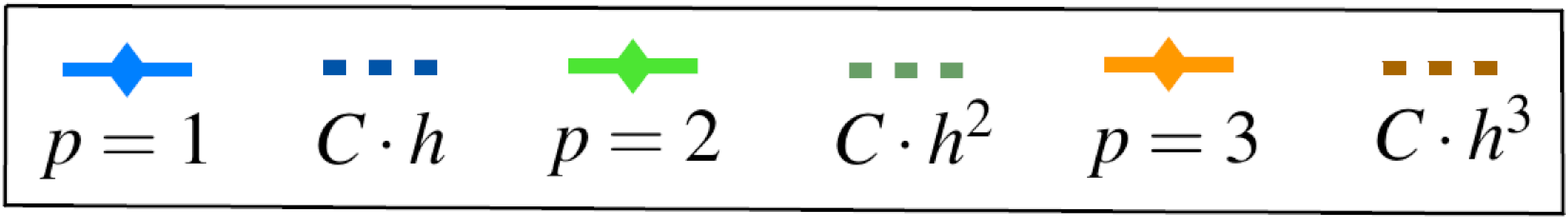}
		\end{center}
		\caption{\label{fig:ConvergenceStudy}\change{dkgreen}{Relative errors in $L^{\infty}E$-norm. \textbf{(Top)} Material set 1, \textbf{(Bottom)} Material set 2. In each we have \textbf{(Left)} convergence plots for Option 1 (elastic tissue) and \textbf{(Right)} convergence plots for Option 2 (acoustic tissue) with simulations using different polynomial degrees $p$.}}
	\end{figure}~\\}{
\begin{figure*}[h!]	
	\noindent{\hspace*{0cm}{\bf Material set 1 with Westervelt-type nonlinearity}}\\[-10cm]
	\hspace*{0cm}% This file was created by matlab2tikz.
\definecolor{mycolor1}{rgb}{0.00000,0.50000,1.00000}%
\definecolor{mycolor2}{rgb}{0.30000,0.90000,0.20000}%
\definecolor{mycolor3}{rgb}{1.00000,0.60000,0.00000}%
\begin{tikzpicture}

\begin{axis}[%
scale=0.75,
width=2.994in,
height=2.846in,
at={(0.0in,0.926in)},
scale only axis,
xmode=log,
xmin=0.310937033868723,
xmax=3.40087380793916,
xminorticks=true,
xlabel style={font=\color{white!15!black}},
xlabel={$h$},
ymode=log,
ymin=1e-05,
ymax=1.5,
yminorticks=true,
ylabel style={font=\color{white!15!black}},
ylabel={},%$\|e\|_{L^{\infty}\textup{E}}/\|(u,\psi)_{\textup{ex}}\|_{L^{\infty}\textup{E}}$},
axis background/.style={fill=white},
title style={font=\bfseries},
title={Option 1},
xmajorgrids,
xminorgrids,
ymajorgrids,
yminorgrids,
legend style={at={(1.03,0.5)}, anchor=west, legend cell align=left, align=left, draw=white!15!black}
]
\addplot [color=mycolor1, line width=2.0pt, mark=diamond, mark options={solid, mycolor1}]
  table[row sep=crcr]{%
2.72069904635133	0.512524115102187\\
1.36034952317566	0.253815045137582\\
0.906899682117109	0.16839689312084\\
0.680174761587832	0.126035026791843\\
0.544139809270265	0.10072168064436\\
0.453449841058554	0.0838841405837072\\
0.388671292335904	0.0718737397960915\\
};
%\addlegendentry{$p=1$}

\addplot [color=black!34!mycolor1, dashed, line width=2.0pt]
  table[row sep=crcr]{%
2.72069904635133	0.244862914171619\\
1.36034952317566	0.12243145708581\\
0.906899682117109	0.0816209713905398\\
0.680174761587832	0.0612157285429049\\
0.544139809270265	0.0489725828343238\\
0.453449841058554	0.0408104856952699\\
0.388671292335904	0.0349804163102314\\
};
%\addlegendentry{$C\cdot h$}

\addplot [color=mycolor2, line width=2.0pt, mark=diamond, mark options={solid, mycolor2}]
  table[row sep=crcr]{%
2.72069904635133	0.077490231476529\\
1.36034952317566	0.0190126226804041\\
0.906899682117109	0.00839710418272116\\
0.680174761587832	0.00471133312745843\\
0.544139809270265	0.00301369181063691\\
0.453449841058554	0.00209538378297584\\
0.388671292335904	0.00154426986916315\\
};
%\addlegendentry{$p=2$}

\addplot [color=black!34!mycolor2, dashed, line width=2.0pt]
  table[row sep=crcr]{%
2.72069904635133	0.0370110165040851\\
1.36034952317566	0.00925275412602127\\
0.906899682117109	0.00411233516712057\\
0.680174761587832	0.00231318853150532\\
0.544139809270265	0.0014804406601634\\
0.453449841058554	0.00102808379178014\\
0.388671292335904	0.000755326867430309\\
};
%\addlegendentry{$C\cdot h^2$}

\addplot [color=mycolor3, line width=2.0pt, mark=diamond, mark options={solid, mycolor3}]
  table[row sep=crcr]{%
2.72069904635133	0.00840216107875206\\
1.36034952317566	0.00104750528373648\\
0.906899682117109	0.00031039083619995\\
0.680174761587832	0.0001320730065555\\
0.544139809270265	6.97741203513801e-05\\
0.453449841058554	4.36966851641945e-05\\
};
%\addlegendentry{$p=3$}

\addplot [color=black!34!mycolor3, dashed, line width=2.0pt]
  table[row sep=crcr]{%
2.72069904635133	0.0040278334922863\\
1.36034952317566	0.000503479186535787\\
0.906899682117109	0.000149179018232826\\
0.680174761587832	6.29348983169735e-05\\
0.544139809270265	3.22226679382903e-05\\
0.453449841058554	1.86473772791032e-05\\
};
%\addlegendentry{$C\cdot h^3$}

\end{axis}

\begin{axis}[%
width=9.825in,
height=7.368in,
at={(0in,0in)},
scale only axis,
xmin=0,
xmax=1,
ymin=0,
ymax=1,
axis line style={draw=none},
ticks=none,
axis x line*=bottom,
axis y line*=left,
legend style={legend cell align=left, align=left, draw=white!15!black}
]
\end{axis}
\end{tikzpicture}%
\hspace*{-20.5cm}% This file was created by matlab2tikz.
\definecolor{mycolor1}{rgb}{0.00000,0.50000,1.00000}%
\definecolor{mycolor2}{rgb}{0.30000,0.90000,0.20000}%
\definecolor{mycolor3}{rgb}{1.00000,0.60000,0.00000}%
\begin{tikzpicture}

\begin{axis}[%
scale=0.75,
width=2.994in,
height=2.846in,
at={(1.094in,0.93in)},
scale only axis,
xmode=log,
xmin=0.310937033868723,
xmax=3.40087380793916,
xminorticks=true,
xlabel style={font=\color{white!15!black}},
xlabel={$h$},
ymode=log,
ymin=1e-05,
ymax=1.5,
yminorticks=true,
ylabel style={font=\color{white!15!black}},
ylabel={},%$\|e\|_{L^{\infty}\textup{E}}/\|(u,\psi)_{\textup{ex}}\|_{L^{\infty}\textup{E}}$},
axis background/.style={fill=white},
title style={font=\bfseries},
title={Option 2},
xmajorgrids,
xminorgrids,
ymajorgrids,
yminorgrids,
legend style={at={(10.03,0.5)}, anchor=west, legend cell align=left, align=left, draw=white!15!black}
]
\addplot [color=mycolor1, line width=2.0pt, mark=diamond, mark options={solid, mycolor1}]
  table[row sep=crcr]{%
2.72069904635133	0.537566498291575\\
1.36034952317566	0.250509282906733\\
0.906899682117109	0.164651613321478\\
0.680174761587832	0.122875112819625\\
0.544139809270265	0.0980784636734661\\
0.453449841058554	0.081631124604748\\
0.388671292335904	0.0699171383576819\\
};
\addlegendentry{$p=1$}

\addplot [color=black!34!mycolor1, dashed, line width=2.0pt]
  table[row sep=crcr]{%
2.72069904635133	0.244862914171619\\
1.36034952317566	0.12243145708581\\
0.906899682117109	0.0816209713905398\\
0.680174761587832	0.0612157285429049\\
0.544139809270265	0.0489725828343238\\
0.453449841058554	0.0408104856952699\\
0.388671292335904	0.0349804163102314\\
};
\addlegendentry{$C\cdot h$}

\addplot [color=mycolor2, line width=2.0pt, mark=diamond, mark options={solid, mycolor2}]
  table[row sep=crcr]{%
2.72069904635133	0.0811283430603854\\
1.36034952317566	0.0194494902075759\\
0.906899682117109	0.00852705213671816\\
0.680174761587832	0.00476542503364779\\
0.544139809270265	0.00304010331180313\\
0.453449841058554	0.0021089386100639\\
0.388671292335904	0.00155066470424568\\
};
\addlegendentry{$p=2$}

\addplot [color=black!34!mycolor2, dashed, line width=2.0pt]
  table[row sep=crcr]{%
2.72069904635133	0.0370110165040851\\
1.36034952317566	0.00925275412602127\\
0.906899682117109	0.00411233516712057\\
0.680174761587832	0.00231318853150532\\
0.544139809270265	0.0014804406601634\\
0.453449841058554	0.00102808379178014\\
0.388671292335904	0.000755326867430309\\
};
\addlegendentry{$C\cdot h^2$}

\addplot [color=mycolor3, line width=2.0pt, mark=diamond, mark options={solid, mycolor3}]
  table[row sep=crcr]{%
2.72069904635133	0.00871803677695128\\
1.36034952317566	0.00108243111017819\\
0.906899682117109	0.000319539615898647\\
0.680174761587832	0.000135239035103899\\
0.544139809270265	7.06827232084058e-05\\
};
\addlegendentry{$p=3$}

\addplot [color=black!34!mycolor3, dashed, line width=2.0pt]
  table[row sep=crcr]{%
2.72069904635133	0.0040278334922863\\
1.36034952317566	0.000503479186535787\\
0.906899682117109	0.000149179018232826\\
0.680174761587832	6.29348983169735e-05\\
0.544139809270265	3.22226679382903e-05\\
};
\addlegendentry{$C\cdot h^3$}

\end{axis}

\begin{axis}[%
width=10in,
height=7.246in,
at={(0in,0in)},
scale only axis,
xmin=0,
xmax=1,
ymin=0,
ymax=1,
axis line style={draw=none},
ticks=none,
axis x line*=bottom,
axis y line*=left,
legend style={legend cell align=left, align=left, draw=white!15!black}
]
\end{axis}
\end{tikzpicture}%
\\[-1cm]
\end{figure*}~\\

\begin{figure}[h!]	
	{\hspace*{0cm}{\bf Material set 2 with Kuznetsov-type nonlinearity}}\\[-10cm]
	\hspace*{0cm}% This file was created by matlab2tikz.
\definecolor{mycolor1}{rgb}{0.00000,0.50000,1.00000}%
\definecolor{mycolor2}{rgb}{0.30000,0.90000,0.20000}%
\definecolor{mycolor3}{rgb}{1.00000,0.60000,0.00000}%
\begin{tikzpicture}

\begin{axis}[%
scale=0.75,
width=2.994in,
height=2.846in,
at={(0.0in,0.926in)},
scale only axis,
xmode=log,
xmin=0.310937033868723,
xmax=3.40087380793916,
xminorticks=true,
xlabel style={font=\color{white!15!black}},
xlabel={$h$},
ymode=log,
ymin=1e-05,
ymax=1.5,
yminorticks=true,
ylabel style={font=\color{white!15!black}},
ylabel={},%$\|e\|_{L^{\infty}\textup{E}}/\|(u,\psi)_{\textup{ex}}\|_{L^{\infty}\textup{E}}$},
axis background/.style={fill=white},
title style={font=\bfseries},
title={Option 1},
xmajorgrids,
xminorgrids,
ymajorgrids,
yminorgrids,
legend style={at={(1.03,0.5)}, anchor=west, legend cell align=left, align=left, draw=white!15!black}
]
\addplot [color=mycolor1, line width=2.0pt, mark=diamond, mark options={solid, mycolor1}]
  table[row sep=crcr]{%
1.36034952317566	0.267510558190737\\
0.906899682117109	0.170705330771081\\
0.680174761587832	0.126440659367848\\
0.544139809270265	0.100598625246998\\
0.453449841058554	0.0835866389956336\\
0.388671292335904	0.0715199119717162\\
};
%\addlegendentry{$p=1$}

\addplot [color=black!34!mycolor1, dashed, line width=2.0pt]
  table[row sep=crcr]{%
1.36034952317566	0.12243145708581\\
0.906899682117109	0.0816209713905398\\
0.680174761587832	0.0612157285429049\\
0.544139809270265	0.0489725828343238\\
0.453449841058554	0.0408104856952699\\
0.388671292335904	0.0349804163102314\\
};
%\addlegendentry{$C\cdot h$}

\addplot [color=mycolor2, line width=2.0pt, mark=diamond, mark options={solid, mycolor2}]
  table[row sep=crcr]{%
2.72069904635133	0.0982407126251855\\
1.36034952317566	0.019186634746431\\
0.906899682117109	0.00836761851441703\\
0.680174761587832	0.00468246035727835\\
0.544139809270265	0.00299858690224098\\
0.453449841058554	0.00209195285647537\\
0.388671292335904	0.00155156858938708\\
};
%\addlegendentry{$p=2$}

\addplot [color=black!34!mycolor2, dashed, line width=2.0pt]
  table[row sep=crcr]{%
2.72069904635133	0.0370110165040851\\
1.36034952317566	0.00925275412602127\\
0.906899682117109	0.00411233516712057\\
0.680174761587832	0.00231318853150532\\
0.544139809270265	0.0014804406601634\\
0.453449841058554	0.00102808379178014\\
0.388671292335904	0.000755326867430309\\
};
%\addlegendentry{$C\cdot h^2$}

\addplot [color=mycolor3, line width=2.0pt, mark=diamond, mark options={solid, mycolor3}]
  table[row sep=crcr]{%
2.72069904635133	0.00876182552080723\\
1.36034952317566	0.00103984427368147\\
0.906899682117109	0.000308225294019168\\
0.680174761587832	0.000133144188493677\\
0.544139809270265	7.38440784462287e-05\\
};
%\addlegendentry{$p=3$}

\addplot [color=black!34!mycolor3, dashed, line width=2.0pt]
  table[row sep=crcr]{%
2.72069904635133	0.0040278334922863\\
1.36034952317566	0.000503479186535787\\
0.906899682117109	0.000149179018232826\\
0.680174761587832	6.29348983169735e-05\\
0.544139809270265	3.22226679382903e-05\\
};
%\addlegendentry{$C\cdot h^3$}

\end{axis}

\begin{axis}[%
width=9.825in,
height=7.368in,
at={(0in,0in)},
scale only axis,
xmin=0,
xmax=1,
ymin=0,
ymax=1,
axis line style={draw=none},
ticks=none,
axis x line*=bottom,
axis y line*=left,
legend style={legend cell align=left, align=left, draw=white!15!black}
]
\end{axis}
\end{tikzpicture}%
\hspace*{-20.5cm}% This file was created by matlab2tikz.
\definecolor{mycolor1}{rgb}{0.00000,0.50000,1.00000}%
\definecolor{mycolor2}{rgb}{0.30000,0.90000,0.20000}%
\definecolor{mycolor3}{rgb}{1.00000,0.60000,0.00000}%
\begin{tikzpicture}

\begin{axis}[%
scale=0.75,
width=2.994in,
height=2.846in,
at={(1.094in,0.93in)},
scale only axis,
xmode=log,
xmin=0.310937033868723,
xmax=3.40087380793916,
xminorticks=true,
xlabel style={font=\color{white!15!black}},
xlabel={$h$},
ymode=log,
ymin=1e-05,
ymax=1.5,
yminorticks=true,
ylabel style={font=\color{white!15!black}},
ylabel={},%$\|e\|_{L^{\infty}\textup{E}}/\|(u,\psi)_{\textup{ex}}\|_{L^{\infty}\textup{E}}$},
axis background/.style={fill=white},
title style={font=\bfseries},
title={Option 2},
xmajorgrids,
xminorgrids,
ymajorgrids,
yminorgrids,
legend style={at={(10.03,0.5)}, anchor=west, legend cell align=left, align=left, draw=white!15!black}
]
\addplot [color=mycolor1, line width=2.0pt, mark=diamond, mark options={solid, mycolor1}]
  table[row sep=crcr]{%
1.36034952317566	0.282470698387141\\
0.906899682117109	0.179319451768731\\
0.680174761587832	0.132679065795229\\
0.544139809270265	0.10552901013792\\
0.453449841058554	0.0876733973220079\\
0.388671292335904	0.0750131826816834\\
};
\addlegendentry{$p=1$}

\addplot [color=black!34!mycolor1, dashed, line width=2.0pt]
  table[row sep=crcr]{%
1.36034952317566	0.12243145708581\\
0.906899682117109	0.0816209713905398\\
0.680174761587832	0.0612157285429049\\
0.544139809270265	0.0489725828343238\\
0.453449841058554	0.0408104856952699\\
0.388671292335904	0.0349804163102314\\
};
\addlegendentry{$C\cdot h$}

\addplot [color=mycolor2, line width=2.0pt, mark=diamond, mark options={solid, mycolor2}]
  table[row sep=crcr]{%
2.72069904635133	0.107598751849621\\
1.36034952317566	0.0206571765559414\\
0.906899682117109	0.00901710057175072\\
0.680174761587832	0.00504785532190627\\
0.544139809270265	0.00322865257526021\\
0.453449841058554	0.00224716492724014\\
0.388671292335904	0.00165959110150751\\
};
\addlegendentry{$p=2$}

\addplot [color=black!34!mycolor2, dashed, line width=2.0pt]
  table[row sep=crcr]{%
2.72069904635133	0.0370110165040851\\
1.36034952317566	0.00925275412602127\\
0.906899682117109	0.00411233516712057\\
0.680174761587832	0.00231318853150532\\
0.544139809270265	0.0014804406601634\\
0.453449841058554	0.00102808379178014\\
0.388671292335904	0.000755326867430309\\
};
\addlegendentry{$C\cdot h^2$}

\addplot [color=mycolor3, line width=2.0pt, mark=diamond, mark options={solid, mycolor3}]
  table[row sep=crcr]{%
2.72069904635133	0.00972718315972817\\
1.36034952317566	0.0011508459818235\\
0.906899682117109	0.000340631589547223\\
0.680174761587832	0.000145633977601148\\
};
\addlegendentry{$p=3$}

\addplot [color=black!34!mycolor3, dashed, line width=2.0pt]
  table[row sep=crcr]{%
2.72069904635133	0.0040278334922863\\
1.36034952317566	0.000503479186535787\\
0.906899682117109	0.000149179018232826\\
0.680174761587832	6.29348983169735e-05\\
};
\addlegendentry{$C\cdot h^3$}

\end{axis}

\begin{axis}[%
width=7in,
height=5.25in,
at={(0in,0in)},
scale only axis,
xmin=0,
xmax=1,
ymin=0,
ymax=1,
axis line style={draw=none},
ticks=none,
axis x line*=bottom,
axis y line*=left,
legend style={legend cell align=left, align=left, draw=white!15!black}
]
\end{axis}
\end{tikzpicture}%
\\[-1.5cm]
	\begin{center}
		\includegraphics[scale=0.15]{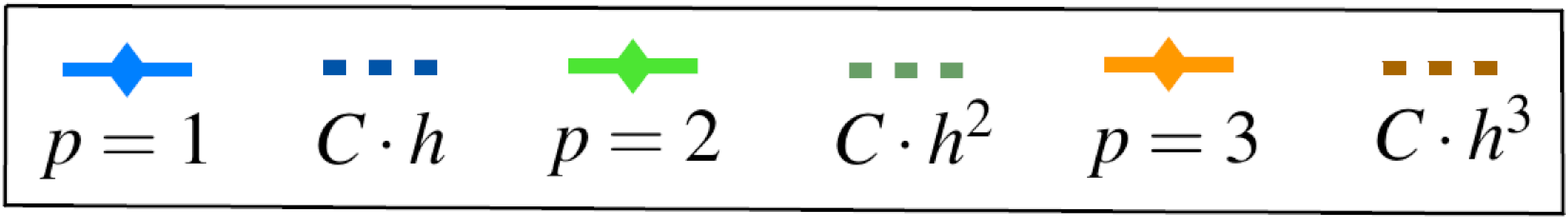}
	\end{center}
	\caption{\label{fig:ConvergenceStudy}\change{dkgreen}{Relative errors in $L^{\infty}E$-norm. \textbf{(Top)} Material set 1, \textbf{(Bottom)} Material set 2. In each we have \textbf{(Left)} convergence plots for Option 1 (elastic tissue) and \textbf{(Right)} convergence plots for Option 2 (acoustic tissue) with simulations using different polynomial degrees $p$.}}
\end{figure}~\\}

\version{\clearpage}{}

\subsection{Numerical Experiment 2: Focused ultrasound propagation into human tissue}
\begin{figure}[h!]
	\begin{center}
		\hspace*{1cm}\includegraphics[scale=0.175]{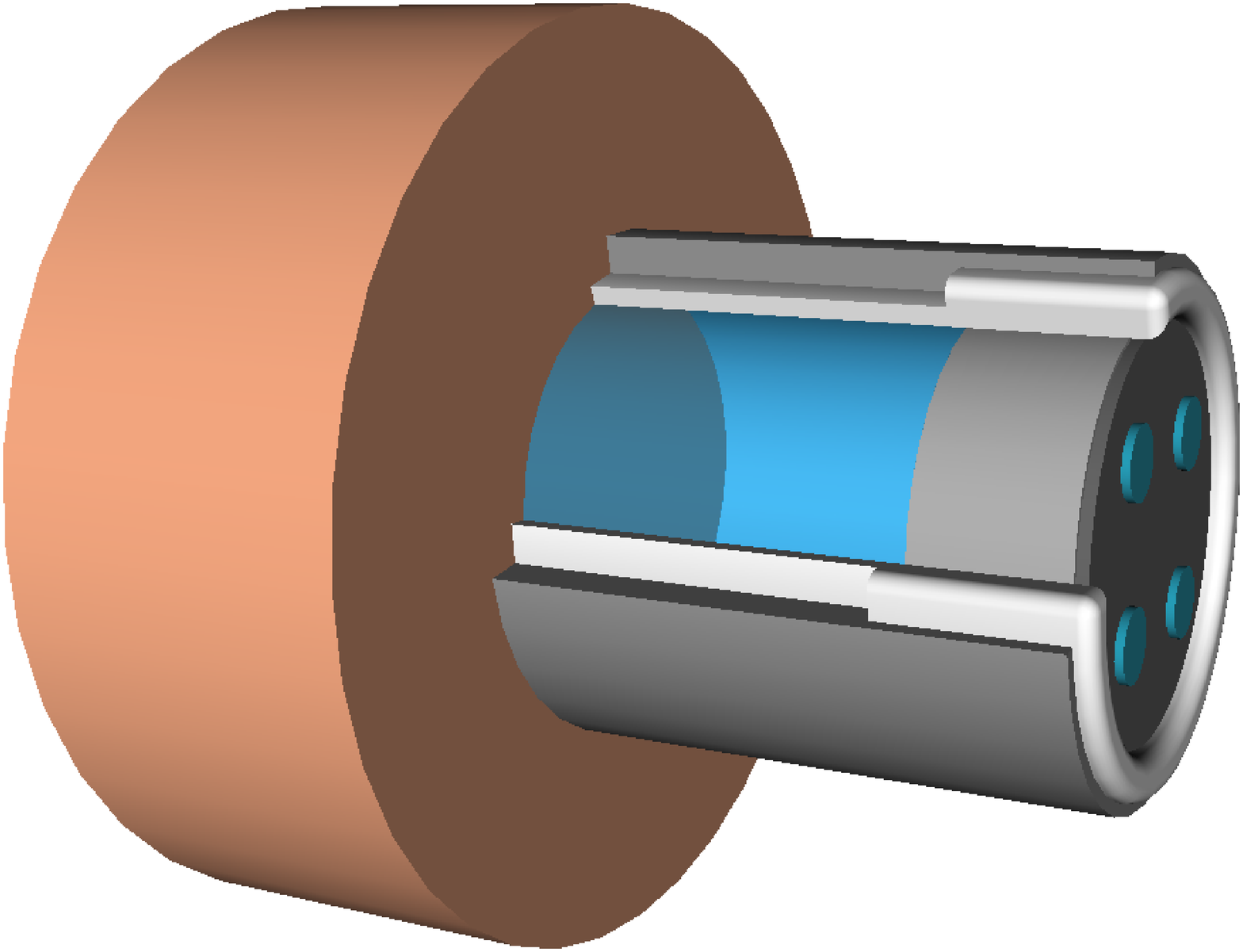}\includegraphics[scale=0.133]{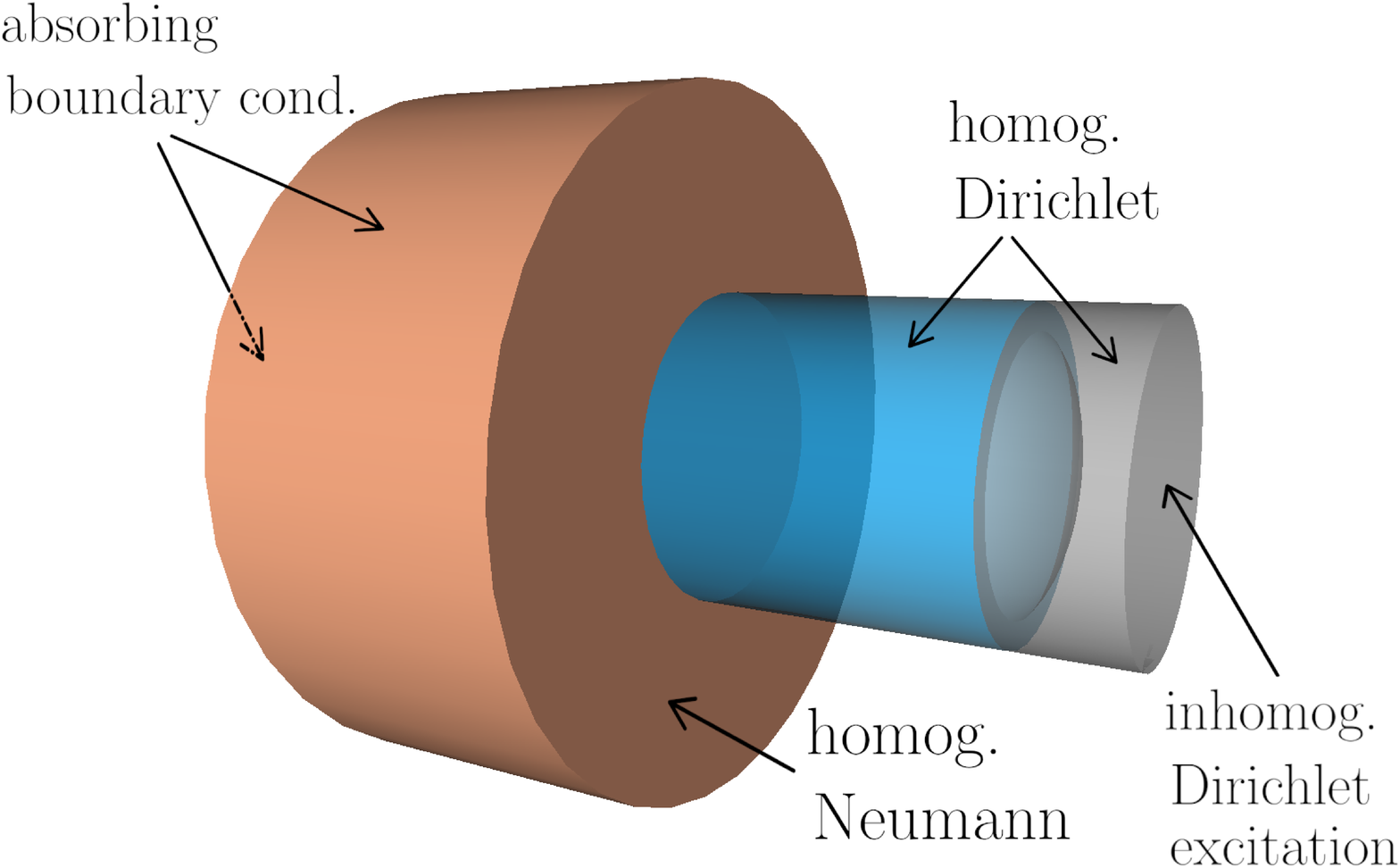}
	\end{center}
	\caption{\textbf{(left)} Depiction of the computational domain embedded into the device framework with walls ($270^{\circ}$ cut-out view), bearings, and excitator plates \textbf{(right)} actual computational domain. The remaining device parts are not numerically resolved but incorporated via boundary conditions.\label{fig:Device}}
\end{figure}

For the final numerical example, we come back to our motivation for the present work, the simulation of medical ultrasound applications involving human tissue. Due to the presence of different boundary conditions (i.e., not only homogeneous Dirichlet data), this experiment lies beyond the theory of this work. Nevertheless, it adheres to the structure of subdomains and coupling given before. Figure \ref{fig:Device} depicts the mechanical device design lying behind the subdomain partition and boundary conditions used.

\paragraph{\bf Geometry and boundary data setup}
We set up a computational domain consisting out of the following three subdomains and data.
\begin{itemize}
	\item \textbf{Curved transducer/excitation array} is modeled as an elastic body made of silicone rubber, where at the bottom side an excitation signal is applied via a (displacement) Dirichlet condition. Using our previous notation, this subdomain takes the role of the purely elastic part $\Omegae$. \version{The circular shape of the excitator has a radius of $r=0.03\,\textup{m}$ with its upper surface being curved according to the height-profile function
	$$
	h:\mathbb{B}_{0.03}(0)\subset\R^2\rightarrow \R,~ h(x,y)=G\left(1+\frac{1}{1+\exp(-KG(\sqrt{x^2+y^2}-S))\left(\frac{G}{A}-1\right)}\right),
	$$
	where we have used the values $G=0.015$, $K=10150$, $A=\frac{G}{2}$, and $S=0.025$; see also Figure \ref{fig:ExcitatorDomain}.\\
	\indent The Dirichlet data on the bottom  surface ($z=0$) is given in a space-time factorized form as \[u_d(t,x,y)=u_d^{(t)}(t)u_d^{(x)}(x,y).\] The spatial part consists out of four spatially separated ``hills'' standing for excitation via piezoelectric plates, whose vibrations induce the prescribed displacement.	More precisely,
	\begin{align*}
		u_d^{(x)}(x,y)&=k(-0.015,-0.015,0.0075,x,y)+k(-0.015,0.015,0.0075,x,y)\\
		&\hspace*{0.7cm}+k(0.015,0.015,0.0075,x,y)+k(0.015,-0.015,0.0075,x,y)
	\end{align*}
	with short-hand notation $r_m:=\sqrt{(x-x_m)^2+(y-y_m)^2}$ and $$k(x_m,y_m,r,x,y)=\begin{cases}\frac{1}{r^4}(r_m-r)^2(r_m+r)^2&\textup{in } \mathds{B}_{r}((x_m,y_m))\\
	0&\textup{else}
	\end{cases}$$
	being such a ``hill'' function; see Figure~\ref{fig:ExcitatorDomain} for a graphical depiction. The temporal part is given as a sine-pulse signal with frequency $f=100\,$kHz, angular frequency $\omega=2\pi f$, and a maximal amplitude $a=0.0025\,$m:
	$$
	u_d^{(t)}(t)=\begin{cases}
	\left(\frac{ft}{2}\right)^2\cdot a\cdot \sin(\omega t)& t\leq \frac{2}{f}\\
	\left(1-\left(\frac{f}{2}(t-\frac{2}{f})\right)^2\right)\cdot a \cdot \sin(\omega t)& \frac{2}{f}<t\leq \frac{4}{f}\\
	0&\frac{4}{f}<t
	\end{cases},
	$$
	which continuously increases over roughly the first 2.25 wavelengths to its maximal amplitude $a$ and then again decreases to zero over the next two wavelengths; see also Figure \ref{fig:ExcitatorDomain} for a plot of the time-signal.\\}{The excitation signal takes the space-time factorized form $u_d(t,x,y)=u_d^{(t)}(t)u_d^{(x)}(x,y)$, where $u_d^{(x)}$ models the spatial distribution of the excitator plates (c.f. Fig.\,\ref{fig:Device}), while $u_d^{(t)}$ is a sine-pulse signal.\\}
	\hspace*{2mm} The remaining surfaces of the excitator block have homogeneous Dirichlet conditions prescribed for the displacement, which can be interpreted as a rigid fixation of the silicon rubber part within the rest of the mechanical device; cf. Figure~\ref{fig:Device}. For a graphical depiction of different parts of the boundary and used conditions, we refer to Figure~\ref{fig:Device}.	
	\version{
	\begin{figure}[h!]
		\begin{center}
		\includegraphics[scale=0.11]{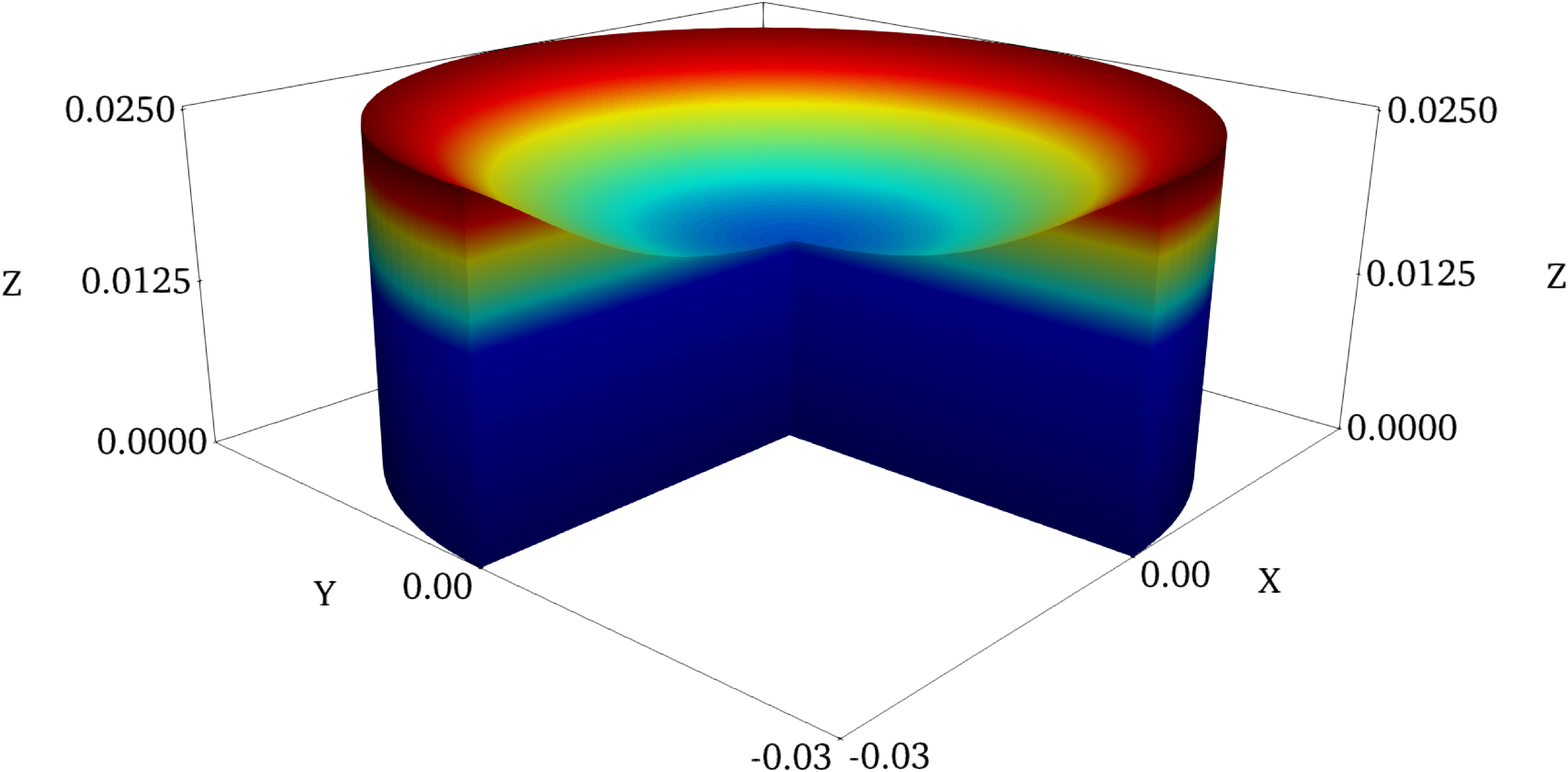}\hspace*{1cm}\includegraphics[scale=0.13]{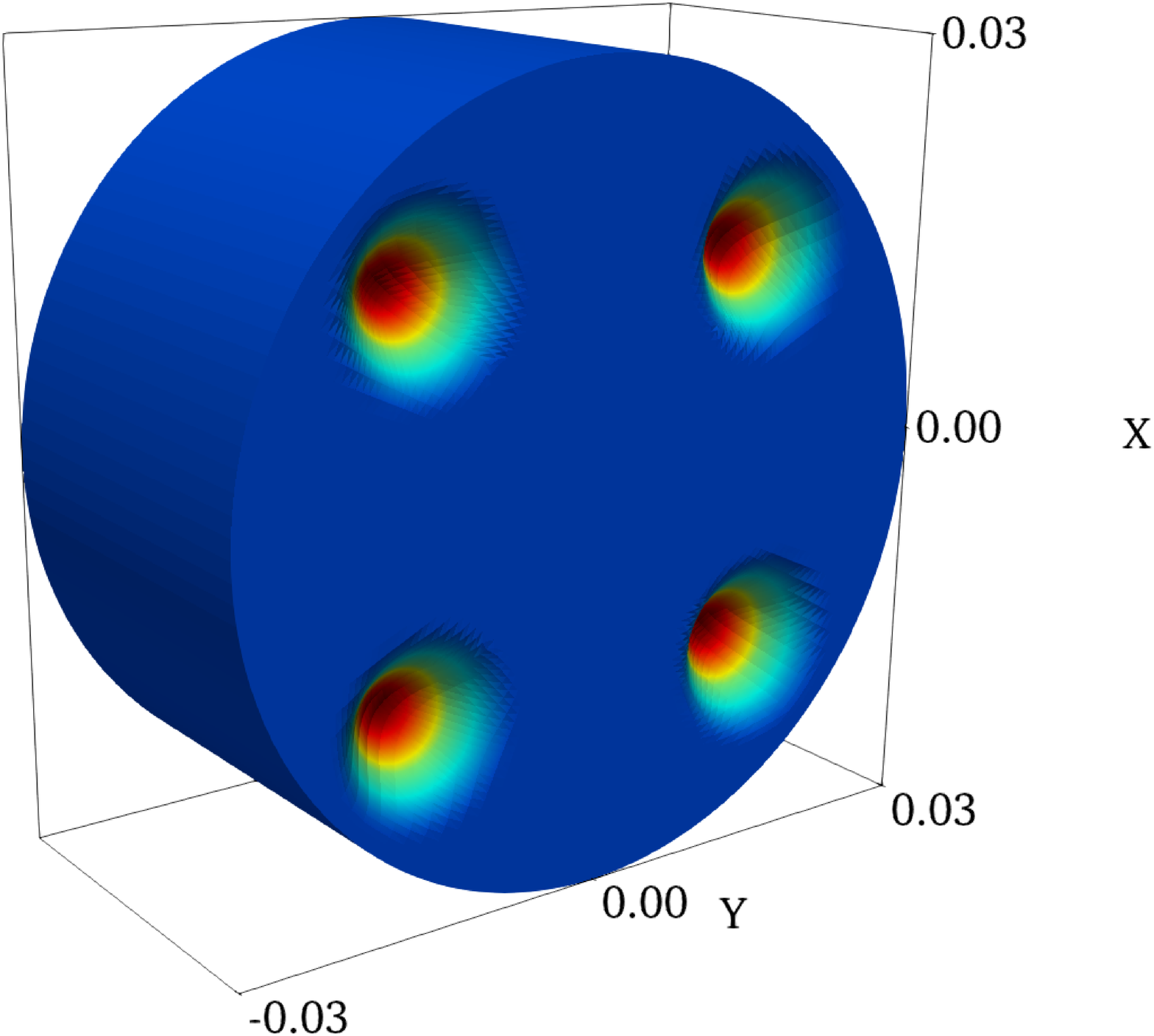}\\[-8cm]
		\hspace*{2.0cm}% This file was created by matlab2tikz.
\definecolor{mycolor1}{rgb}{0.00000,0.60000,1.00000}%
\begin{tikzpicture}

\begin{axis}[%
scale=0.30,
width=6.002in,
height=4.248in,
at={(1.106in,0.849in)},
scale only axis,
xmin=0,
xmax=8e-05,
xlabel style={font=\color{white!15!black}},
xlabel={$t$},
ymin=-1.5,
ymax=1.5,
ylabel style={font=\color{white!15!black}},
ylabel={$u_d^{(t)}(t)$},
axis background/.style={fill=white},
xmajorgrids,
ymajorgrids,
legend style={legend cell align=left, align=left, draw=white!15!black}
]
\addplot [color=mycolor1, line width=2.0pt]
  table[row sep=crcr]{%
0	0\\
1.75438596539124e-07	8.46477539195867e-06\\
3.50877192967225e-07	6.73071995095365e-05\\
5.26315789506349e-07	0.000224861128258036\\
7.0175438593445e-07	0.000525424974449296\\
8.77192982473574e-07	0.00100737325713274\\
1.05263157890167e-06	0.00170142025675812\\
1.2280701754408e-06	0.0026290815061254\\
1.40350877197992e-06	0.00380137470442021\\
1.57894736840802e-06	0.00521779660966948\\
1.75438596494715e-06	0.00686560663756841\\
2.10526315791437e-06	0.0107412771849233\\
2.80701754384882e-06	0.0193329934035564\\
2.98245614038795e-06	0.021223687714926\\
3.15789473681605e-06	0.0228309139609294\\
3.33333333335517e-06	0.0240562612162344\\
3.50877192978327e-06	0.0248004033752225\\
3.6842105263224e-06	0.0249657005699332\\
3.85964912286152e-06	0.0244589088672749\\
4.03508771928962e-06	0.0231939413991175\\
4.21052631582874e-06	0.0210946212980421\\
4.38596491225685e-06	0.018097365148307\\
4.56140350879597e-06	0.0141537351673936\\
4.73684210522407e-06	0.0092327990392933\\
4.91228070176319e-06	0.00332323823630343\\
5.0877192981913e-06	-0.00356485122032024\\
5.26315789473042e-06	-0.0113985173324608\\
5.61403508769764e-06	-0.0296507230589861\\
5.96491228066487e-06	-0.0506846810158408\\
7.01754385967757e-06	-0.117501385272946\\
7.19298245610567e-06	-0.126948288716321\\
7.36842105264479e-06	-0.135270471349936\\
7.54385964907289e-06	-0.142220525368508\\
7.71929824561202e-06	-0.147556998757941\\
7.89473684215114e-06	-0.151049210412984\\
8.07017543857924e-06	-0.152482131915433\\
8.24561403511836e-06	-0.151661250623887\\
8.42105263154647e-06	-0.148417325786155\\
8.59649122808559e-06	-0.142610947895513\\
8.77192982451369e-06	-0.134136811537011\\
8.94736842105281e-06	-0.122927613550772\\
9.12280701759194e-06	-0.108957491491476\\
9.29824561402004e-06	-0.0922449220767467\\
9.47368421055916e-06	-0.072855005555622\\
9.64912280698726e-06	-0.050901069629062\\
9.99999999995449e-06	-3.33066907387547e-16\\
1.03508771930327e-05	0.0585740903731449\\
1.08771929824281e-05	0.154893712016728\\
1.14035087719344e-05	0.250950126971487\\
1.17543859649016e-05	0.308197081960447\\
1.19298245614408e-05	0.333212371444689\\
1.21052631578689e-05	0.355133476926527\\
1.2280701754408e-05	0.373465544377021\\
1.24561403508361e-05	0.387741302532531\\
1.26315789473752e-05	0.397529548456954\\
1.28070175439143e-05	0.402443444717\\
1.29824561403424e-05	0.402148491096659\\
1.31578947368816e-05	0.396370034043914\\
1.33333333333097e-05	0.384900179459751\\
1.35087719298488e-05	0.367603979029236\\
1.36842105262769e-05	0.344424767046425\\
1.3859649122816e-05	0.315388533555088\\
1.40350877192441e-05	0.280607230537528\\
1.42105263157832e-05	0.240280920723011\\
1.43859649123224e-05	0.194698693211546\\
1.47368421052896e-05	0.089364375886493\\
1.50877192982568e-05	-0.0313503443822685\\
1.56140350877632e-05	-0.229358766943582\\
1.61403508771585e-05	-0.427727695563251\\
1.64912280701257e-05	-0.547840910558663\\
1.66666666666648e-05	-0.60140653040586\\
1.6842105263204e-05	-0.649412663777547\\
1.70175438596321e-05	-0.69098158377072\\
1.71929824561712e-05	-0.725289330655296\\
1.73684210525993e-05	-0.751579302551429\\
1.75438596491384e-05	-0.76917536705521\\
1.77192982455665e-05	-0.777494289426529\\
1.78947368421056e-05	-0.776057276610711\\
1.80701754386448e-05	-0.764500443048597\\
1.82456140350729e-05	-0.742584013919403\\
1.8421052631612e-05	-0.710200094093905\\
1.85964912280401e-05	-0.667378846544767\\
1.87719298245792e-05	-0.614292942114894\\
1.89473684210073e-05	-0.551260163189621\\
1.91228070175464e-05	-0.478744066719759\\
1.92982456140856e-05	-0.397352636927248\\
1.96491228070528e-05	-0.211075377661804\\
2.000000000002e-05	-2.22044604925031e-15\\
2.05263157895264e-05	0.324474608076422\\
2.08771929824936e-05	0.522665540730742\\
2.12280701754608e-05	0.694667719587869\\
2.14035087718889e-05	0.768115276211899\\
2.1578947368428e-05	0.831948681652858\\
2.17543859649671e-05	0.885388631980824\\
2.19298245613953e-05	0.927791799816256\\
2.21052631579344e-05	0.958658988754407\\
2.22807017543625e-05	0.977641334229398\\
2.24561403509016e-05	0.984544469830669\\
2.26315789473297e-05	0.979330606354893\\
2.28070175438688e-05	0.962118499848862\\
2.29824561402969e-05	0.933181314164582\\
2.3157894736836e-05	0.892942412694128\\
2.33333333333752e-05	0.841969142568207\\
2.35087719298033e-05	0.780964702285757\\
2.36842105263424e-05	0.710758210103202\\
2.40350877193096e-05	0.54661411611115\\
2.43859649122768e-05	0.358212006799531\\
2.50877192982113e-05	-0.0515229091355425\\
2.54385964911785e-05	-0.251982548431329\\
2.57894736842568e-05	-0.436065374645462\\
2.6140350877224e-05	-0.594846727844403\\
2.63157894736521e-05	-0.662355321243122\\
2.64912280701912e-05	-0.720885725109279\\
2.66666666666193e-05	-0.769800358919499\\
2.68421052631584e-05	-0.808594869449582\\
2.70175438596976e-05	-0.83690361660656\\
2.71929824561257e-05	-0.854503204536096\\
2.73684210526648e-05	-0.861314021656734\\
2.75438596490929e-05	-0.857399781656443\\
2.7719298245632e-05	-0.842965085879563\\
2.78947368420601e-05	-0.818351055526346\\
2.80701754385992e-05	-0.784029109281624\\
2.82456140351384e-05	-0.740592987994508\\
2.84210526315665e-05	-0.688749152476376\\
2.85964912281056e-05	-0.62930570302081\\
2.89473684210728e-05	-0.491285099138898\\
2.929824561404e-05	-0.33453151341966\\
3.03508771929417e-05	0.160107000833229\\
3.070175438602e-05	0.304582289876047\\
3.08771929824481e-05	0.368779201971147\\
3.10526315789872e-05	0.426631129382086\\
3.12280701754153e-05	0.477527049071756\\
3.14035087719544e-05	0.520966523944831\\
3.15789473683825e-05	0.55656497169808\\
3.17543859649216e-05	0.584057156657946\\
3.19298245613497e-05	0.603298869752365\\
3.21052631578889e-05	0.614266789012803\\
3.2280701754428e-05	0.617056540260482\\
3.24561403508561e-05	0.61187900449242\\
3.26315789473952e-05	0.599054944549716\\
3.28070175438233e-05	0.579008048535418\\
3.29824561403624e-05	0.55225651078285\\
3.31578947367905e-05	0.519403292611146\\
3.33333333333297e-05	0.481125224324688\\
3.36842105262969e-05	0.391299143626709\\
3.40350877192641e-05	0.28920082697274\\
3.49122807017377e-05	0.0244622549256481\\
3.52631578947049e-05	-0.0687330595147377\\
3.5438596491244e-05	-0.109963634762055\\
3.56140350876721e-05	-0.14695060500425\\
3.57894736842113e-05	-0.179304281033357\\
3.59649122807504e-05	-0.206728608122567\\
3.61403508771785e-05	-0.229024328484481\\
3.63157894737176e-05	-0.246090477046483\\
3.64912280701457e-05	-0.257924195102313\\
3.66666666666848e-05	-0.264618873378578\\
3.68421052631129e-05	-0.26636066287751\\
3.70175438596521e-05	-0.263423418108787\\
3.71929824561912e-05	-0.256162162597123\\
3.73684210526193e-05	-0.245005190455241\\
3.75438596491584e-05	-0.230444939969742\\
3.77192982455865e-05	-0.213027795210975\\
3.80701754385537e-05	-0.172010798148459\\
3.87719298245992e-05	-0.0830038589791028\\
3.89473684210273e-05	-0.0629525495000495\\
3.91228070175664e-05	-0.0449288472681201\\
3.92982456139945e-05	-0.0294237985691584\\
3.94736842105337e-05	-0.0168645846193569\\
3.96491228070728e-05	-0.00760571354457429\\
3.98245614035009e-05	-0.00192150401397806\\
4.000000000004e-05	-0\\
7.00000000000145e-05	0\\
};
%\addlegendentry{data1}

\end{axis}

\begin{axis}[%
width=7.854in,
height=5.51in,
at={(0in,0in)},
scale only axis,
xmin=0,
xmax=1,
ymin=0,
ymax=1,
axis line style={draw=none},
ticks=none,
axis x line*=bottom,
axis y line*=left,
legend style={legend cell align=left, align=left, draw=white!15!black}
]
\end{axis}
\end{tikzpicture}%
\\[-2.5cm]
		\end{center}~\\[-2.5cm]
		\caption{\textbf{(top left)} Exemplary cut-open view at the radial symmetric, curved excitator subdomain. \textbf{(top right)} View at bottom surface of excitator domain with spatial displacement Dirichlet condition $u_d^{(x)}$ applied. \textbf{(bottom)} Time signal $u_d^{(t)}$ scaling the Dirichlet displacement with respect to time (here with amplitude $a=1$).\label{fig:ExcitatorDomain}}
	\end{figure}}{}
	
	\item \textbf{Acoustic conductor pipe} is modeled as an acoustic/fluid body \version{of length $L\approx0.055\,$m}{}filled with water. This subdomain takes the role of the purely acoustic/fluid part $\Omegaf$ from before. On its bottom end, the subdomain aligns with the elastic excitator domain on the interface $\Gamma_{\textup{I}}^{\textup{e,f}}$. On its top end, the interface $\Gamma_{\textup{I}}^{\textup{f,t}}$ with the tissue domain is placed. The cylinder mantle surface has homogeneous Dirichlet conditions for the acoustic potential prescribed, modeling the devices rigid walls encasing the water (compare to Figure \ref{fig:Device}).
	\item \textbf{The tissue domain} resides on top of the acoustic conductor. It can be modeled as an elastic or acoustic body, corresponding to the two previously discussed options. In both cases, the computational tissue domain has the shape of a cylinder \version{of height $h=0.08\,\textup{m}$ and radius $r=0.06\,\textup{m}$}cut out of the larger physical tissue domain, being (part of) the human body. To avoid unphysical reflections on the mantle and top surfaces resulting from that truncation procedure, absorbing boundary conditions of Engquist--Majda type~\citep{EM} (in the acoustic case) or as in \citep{stacey1988improved,antonietti2018numerical} (in the elastic case) are employed on these surfaces. The remaining surface is the (skin) surface of the human body without the interface $\Gamma_{\textup{I}}^{\textup{f,t}}$, equipped with homogeneous Neumann conditions.
\end{itemize}

For this numerical experiment, the simulation starts from the zero initial conditions $(\psi,\dot{\psi})_{t=0}=(0,0)$, $(\u,\dot{\u})_{t=0}=(\boldsymbol{0},\boldsymbol{0})$; in other words, solid and fluid bodies being at rest. Furthermore, no additional external forces are applied; that is, $\fe=\ft=0$ and all the dynamics of the system is induced via the excitation/Dirichlet conditions described above.

\paragraph{\bf Material parameters} In contrast to the synthetic test case, parameters for realistic materials are often given not in terms of the Lamé-parameters $\lambda$ and $\mu$, but rather in terms of the elastic modulus $E$ and Poisson ratio $\nu$. Similarly in nonlinear acoustics, the coefficient of nonlinearity $\beta=1+\frac12 B/A$ is commonly used to indicate the nonlinear wave behavior. Here the ratio $B/A$ denotes the parameter of nonlinearity, which arises from the Taylor expansion of the variations of the pressure in terms of variations of the density in a given medium; cf.~\citep[\S 2]{hamilton1998nonlinear}. \\
\indent \version{The parameter values used in the upcoming experiments are given in Table~\ref{tab:MatParam}, from which the corresponding values, e.g., of $\lambda$ and $\mu$ can be computed using the well-known formulas  \[\lambda = \frac{E\nu}{(1+\nu)(1-2\nu)},\quad \mu= \frac{E}{2(1+\nu)}.\]}{} \version{For the acoustic case, we}{We} can choose between two different well-established models. The choice $k_1=\frac{2+B/A}{2c^2}$ and $k_2=0$ leads to Westervelt's equation, while $k_1 = \frac{B/A}{2c^2}, k_2 = 1$ leads to Kuznetsov's equation; see~\citep{MKaltenbacher,hamilton1998nonlinear}.\\
\version{\indent The derivation of the elastic damping coefficient $\zeta$ can be based on different reference parameters from literature. For silicone rubber, values for the loss coefficient $c_L$ are given, which relates to the quality factor $Q$ via $c_L=\frac{1}{Q}$, where $Q$ describes the so-called intensity loss $\Delta I$ of a wave per period, i.e. $Q=2\pi\frac{I_0}{\Delta I}$. From $Q$ one then finds $\zeta=\frac{\pi f}{Q}=c_L\pi f$ \citep{antonietti2018numerical}, where $f$ is the frequency used in the simulation. For tissue values for the attenuation parameter $\alpha$, depending on the used frequency $f$ are given, which yields the intensity loss after traveling through a medium for one wave length $\lambda$ as $I(\lambda) = I_0\exp(-\alpha \lambda)$. Using $\lambda = c/f$ and the correspondence of one wave length to one period, this translates to $\frac{1}{Q}=\frac{1}{2\pi}(1-\exp(-\alpha c/f))$ and hence $\zeta = \frac{f}{2}(1-\exp(-\alpha c/f))$, where in the elastic medium the wave speed $c$ was chosen to be the pressure wave speed $c=v_p=\sqrt{\frac{\lambda+2\mu}{\varrho}}$.}{}
\paragraph{\bf Discretization} As a final time for this simulation, we have chosen $T=1.5\cdot 10^{-4}\,$s which is resolved by 15000 time steps of $10\,$ns each. For the nonlinear acoustic part, we employ the Generalized $\alpha$ scheme as in \citep{antonietti2020high} to damp high-frequency oscillations that might appear in a Gibb phenomenon-like manner near a steep wavefront. In detail, we use the method parameters $\beta=4/9, \gamma=5/6, \alpha_m=0$, and $\alpha_f=1/3$. We refer to \citep{chung1993time,Genalpha_Bonaventura} for a deeper insight into this time-stepping method. For the elastic part, again the leapfrog scheme is employed. \change{blue}{In this experiment the meshes of the subdomain were chosen to be conforming at the interfaces.} Spatial discretization is done using 686937 degrees of freedom with polynomial basis functions of degree $p=3$. 

\begin{figure}[h!]
	\hspace*{-5mm}\includegraphics[trim=1cm 8cm 0cm 1.5cm, clip, scale=0.75]{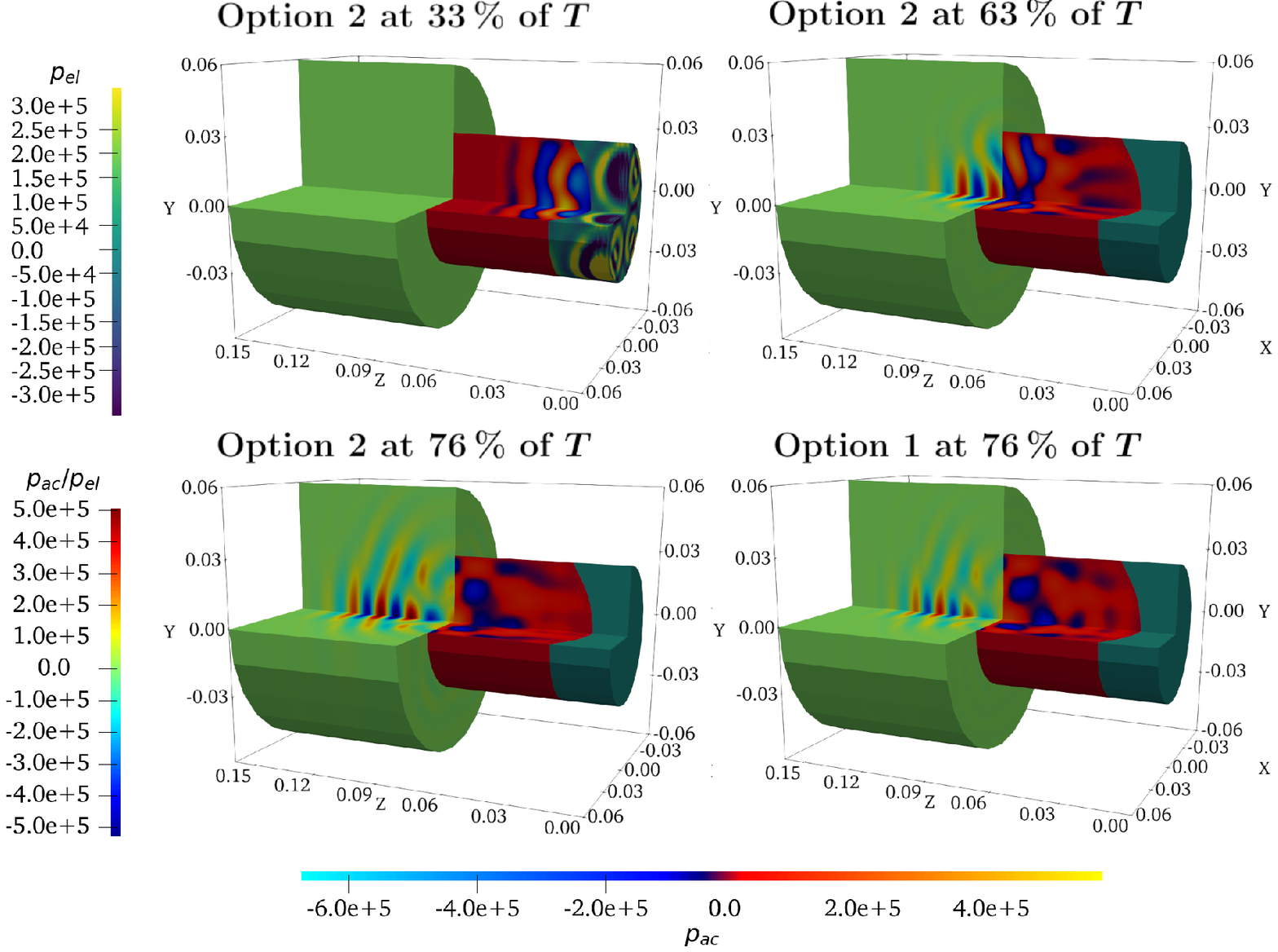}~\\
	\vspace*{-5cm}
	\caption{Pressure field over time using Westervelt's equation as acoustic model. \textbf{(top left - bottom left)} Option 2, \textbf{(bottom right)} Option 1 and Westervelt's equation only in the fluid domain.\label{Fig:FinalSimulationResults10_3D}}
\end{figure}

\paragraph{\bf Results} From the computed acoustic potential $\psi$ and displacement field $\u$ from our numerical simulation, we post-process the acoustic and elastic pressure values relevant for the application according to $p_{\textup{ac}}=\varrho\dot{\psi}$ and $p_{\textup{el}}=-\frac{1}{3}\sum_{i=1}^3\boldsymbol{\sigma}_{ii}$, respectively, depending on whether a subdomain contains an acoustic or elastic material. Herein $\boldsymbol{\sigma}(\u)=\lambda\textup{div}(\u)\mathds{1}+2\mu\varepsilon(\u)$ is the stress-tensor, as defined before. Figure~\ref{Fig:FinalSimulationResults10_3D} shows the pressure wave propagation within the computational domain over different time steps for Option 2 (acoustic) choice for the tissue using the Westervelt equation as the acoustic wave model.
\version{
\begin{table}[h!]\renewcommand{\arraystretch}{1.4}
	\normalsize\setlength\tabcolsep{10pt}
	\begin{tabular}{|c|c||c|c|c||c|}
		\cline{3-5}
		\multicolumn{2}{c||}{}&\multicolumn{3}{|c||}{\textbf{materials}} & \multicolumn{1}{c}{}\\ \hline
		\textbf{type} & \textbf{parameter} & \textbf{silicone rubber} & \textbf{water} & \textbf{tissue} & \textbf{unit}\\ \hline\hline
		all& $\varrho$ & 1100 & 998.23 & 1000& $\frac{\textup{kg}}{\textup{m}^3}$\\ \hline
		& $E$   &  $5\cdot 10^7$& / & $1\cdot 10^9$ & $\frac{\textup{N}}{\textup{m}^2}$\\
		elastic& $\nu$		   &0.49 & / &0.375& 1\\ 
		&$c_L$            &  1.0  & /  &   /   & 1 \\
		&$\alpha$             &  /    &  / &  0.2668    &  $\frac{\textup{Np}}{\textup{m}}$\\ \hline
		& $c$		   & / & 1500& 1540 & $\frac{\textup{m}}{\textup{s}}$\\
		acoustic& $b$  & / & $6\cdot 10^{-9}$ & $6.4117\cdot 10^{-4}$ &$\frac{\textup{m}^2}{\textup{s}}$\\
		& $B/A$ 	   &  / & 5 & 7.44 &1\\ \hline
	\end{tabular}
	\centering	\caption{\textbf{Material parameters.} 
		Parameters of \textbf{water}: $\varrho$ at $20\,^{\circ}$C: \textit{{http://www.matweb.com/}}, $c,b,B/A$: \textit{\citep[\S 5]{MKaltenbacher}}.\\
		Elastic parameters  for \textbf{tissue}: $E$ (order of magnitude), $\nu$ (mean value) \textit{\citep{chawla2006mechanical}},	acoustic parameters for \textbf{tissue}: $\varrho,c,b$: \citep{velasco2015finite}, $B/A, \alpha$: \textit{{https://itis.swiss/virtual-population/tissue-properties/database/} (average value of kidney tissue)}. 
		Elastic parameters for \textbf{silicon rubber}: $\varrho, E,\nu$: \textit{{https://www.azom.com/properties.aspx?ArticleID=920}}, $c_L$ as well (mean value used).\label{tab:MatParam}}
\end{table}}{}

\begin{figure}[h!]
	\vspace*{-5cm}
	\hspace*{-1.55cm}\input{Amp10Spatial230.tex}\hspace*{-22cm}\input{Amp10Temporal120.tex}\\[-1.5cm]
	\begin{center}
		\includegraphics[scale=0.12]{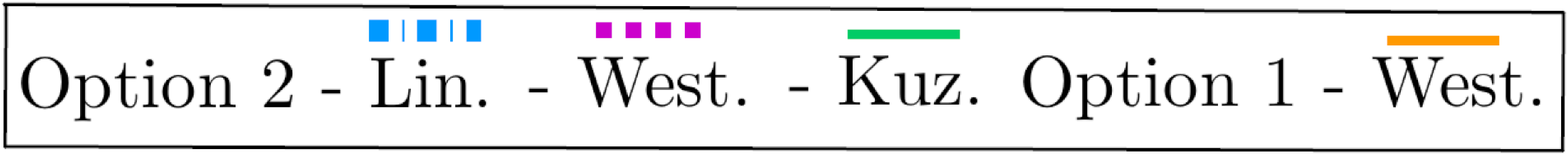}
	\end{center}
	\caption{Model comparison of pressure signal at \textbf{(left)} axis of symmetry of tissue domain at 76\,\% of simulation time. \textbf{(right)} Same comparison but over time at fixed point $p=(0,0,0.12)^{\top}\,\textup{m}$.\label{Fig:FinalSimulationResults10_1D}}
\end{figure}

Further, the pressure field within the tissue is compared with Option 1 of having an elastic tissue model. Figure \ref{Fig:FinalSimulationResults10_1D} then compares different choices for the acoustic wave equation being the linear wave equation, the Westervelt or the Kuznetsov equation. 

First it can be seen that, while qualitatively behaving similar, Options 1 and 2 differ in amplitude and propagation velocity of the wave, which is due to the material parameters for the two options \version{(see Table \ref{tab:MatParam})}{} stemming from different references with some of them being available only approximately.

\begin{figure}[h!]
	\hspace*{-5mm}\includegraphics[trim=1cm 8cm 0cm 1cm, clip, scale=0.75]{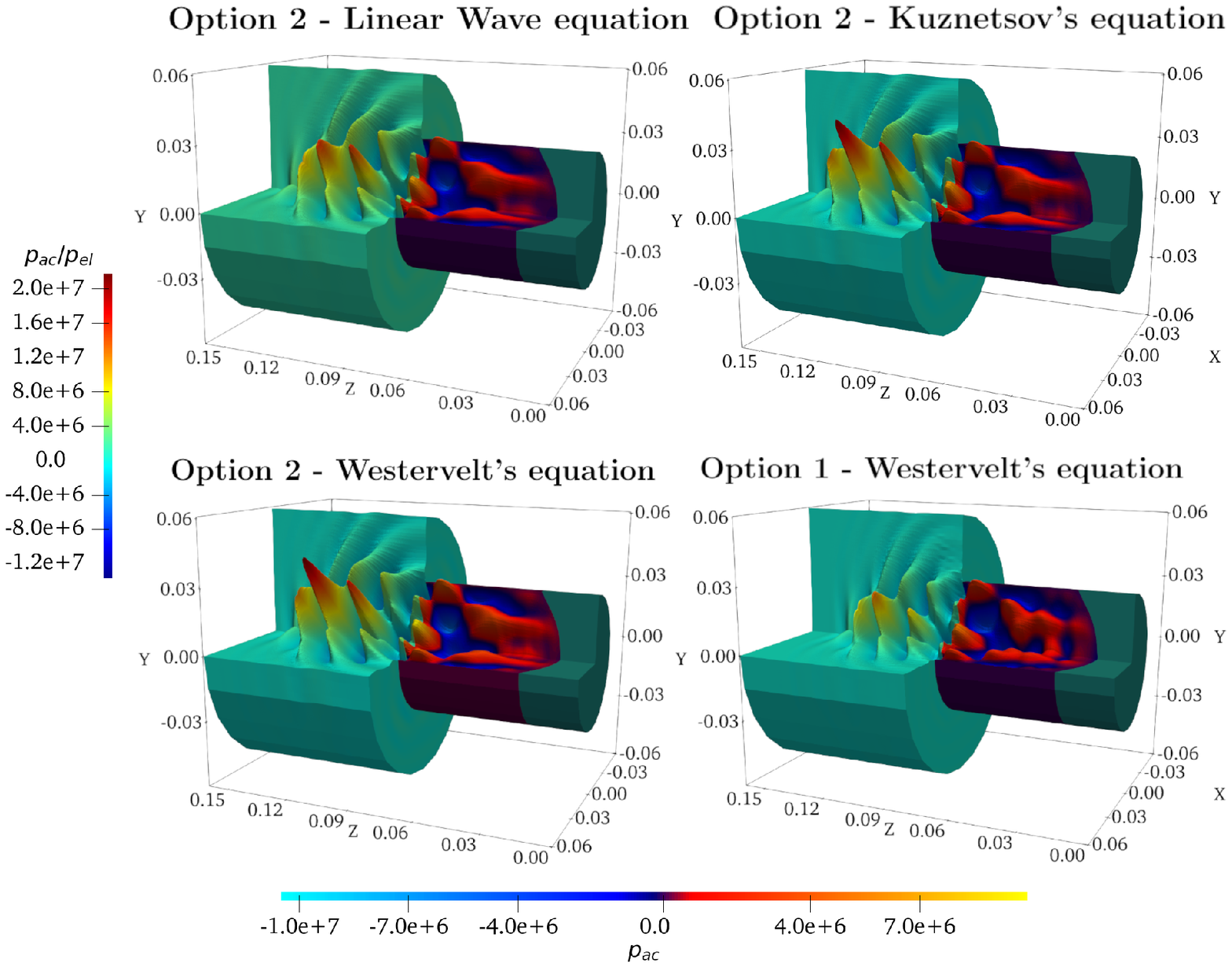}	~\\
	\vspace*{-5cm}
	\caption{Comparison of acoustic models in high pressure regime, all at 76\,\% simulation time in tissue domain. \textbf{(top left)} Linear wave equation, \textbf{(top right)} Kuznetsov's equation, \textbf{(bottom left)} Westervelt's equation. \textbf{(bottom right)} Option 1 with Westervelt's equation only in the fluid domain.\label{Fig:FinalSimulationResults300_3D}}
\end{figure}

\begin{figure}[h!]
	\vspace*{-5cm}
	\hspace*{-1.5cm}\input{Amp300Spatial230.tex}\hspace*{-30cm}\raisebox{-2mm}{\input{Amp300Temporal120.tex}}\\[-1.5cm]
	\begin{center}
		\includegraphics[scale=0.18]{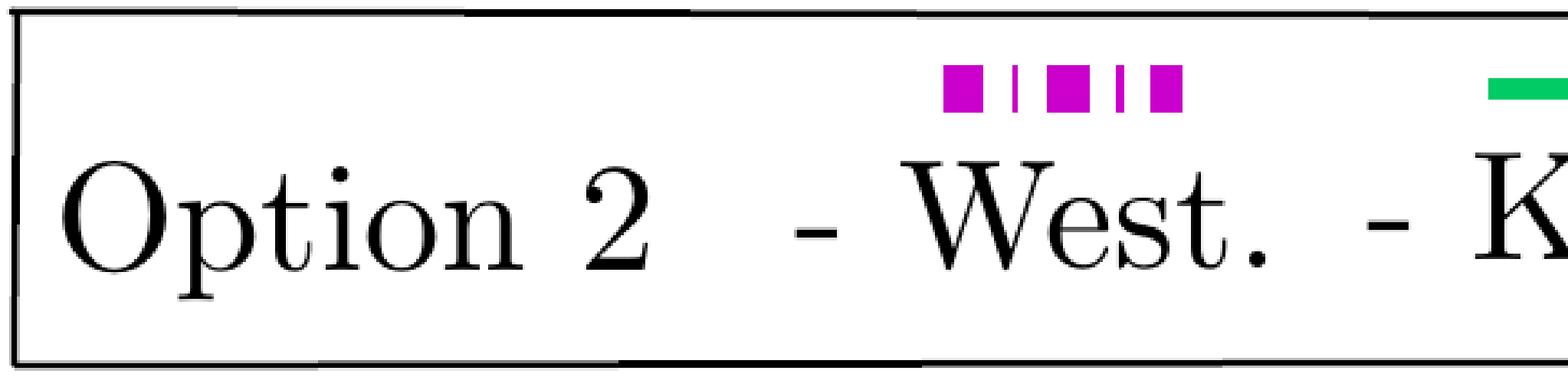}´
	\end{center}
	\caption{Model comparison of pressure signal at \textbf{(left)} axis of symmetry of tissue domain at 76\,\% of simulation time. \textbf{(right)} Same comparison but over time at fixed point $p=(0,0,0.12)^{\top}\,\textup{m}$.\label{Fig:FinalSimulationResults300_1D}}
\end{figure}

In contrast to that, one observes almost no visible difference among the different acoustic models (i.e, the linear wave equation, Westervelt's, and Kuznetsov's equations) in the given pressure regime. This observation changes once we employ higher excitation amplitudes and hence increase the influence of the nonlinear terms of the models. Figures \ref{Fig:FinalSimulationResults300_3D} and \ref{Fig:FinalSimulationResults300_1D} directly compare the acoustic models with each other for a simulation with an excitation amplitude of $a=0.075\,\textup{m}$. We note that this amplitude might be exaggerated from this application/experiment point of view. However, it shows very accurately the wave steepening effect modeled by the nonlinear terms in the higher pressure regime in contrast to the linear model. Even in Option 1 simulations with the linear elastic tissue model, a small steepening and amplitude increase of the pressure wave is visible once the nonlinear wave equations are used in the fluid region compared to a completely linear simulation.

The comparably quite small difference between the Westervelt and Kuznetsov equation's results further shows that the approximative assumption of $\nabla\psi\cdot\nabla\dot{\psi}\approx\frac{1}{c^2}\dot{\psi}\ddot{\psi}$, which is made in the derivation of the Westervelt equation as a simplification of Kuznetsov's equation and which holds with equality in the case of a linear plane wave, is still reasonable in the given regime of nonlinearity.

%\begin{figure}[h!]
%	\hspace*{-5mm}\includegraphics[trim=1cm 8cm 0cm 1cm, clip, scale=0.75]{images/Picture4.png}	~\\
%	\vspace*{-5cm}
%	\caption{Comparison of acoustic models in high pressure regime, all at 76\,\% simulation time in tissue domain. \textbf{(top left)} Linear wave equation, \textbf{(top right)} Kuznetsov's equation, \textbf{(bottom left)} Westervelt's equation. \textbf{(bottom right)} Option 1 with Westervelt's equation only in the fluid domain.\label{Fig:FinalSimulationResults300_3D}}
%\end{figure}
%
%\begin{figure}[h!]
%	\vspace*{-5cm}
%	\hspace*{-1.5cm}\input{images/Amp300Spatial230.tex}\hspace*{-30cm}\raisebox{-2mm}{\input{images/Amp300Temporal120.tex}}\\[-1.5cm]
%	\begin{center}
%		\includegraphics[scale=0.18]{images/Legend2.png}
%	\end{center}
%	\caption{Model comparison of pressure signal at \textbf{(left)} axis of symmetry of tissue domain at 76\,\% of simulation time. \textbf{(right)} Same comparison but over time at fixed point $p=(0,0,0.12)^{\top}\,\textup{m}$.\label{Fig:FinalSimulationResults300_1D}}
%\end{figure}

\section{Conclusion}
In this work, we have considered a coupled elasto-acoustic problem with thermoviscous dissipation and general nonlinearities of quadratic type in the acoustic regime. The mathematical model was in particular motivated by the medical applications of high-intensity ultrasound, which tend to involve different elastic and acoustic subdomains with possible jumps in the material parameters. \\
\indent The elasto-acoustic interface was resolved by a coupling based on force-exchange via Neumann conditions, while the acoustic-acoustic interfaces were treated using a discontinuous Galerkin approach, resulting in a fully coupled initial boundary-value problem. We discretized the problem in space using hexahedral elements of degree $p$ being conform within each material subdomain, but with the option of non-matching grids across interfaces. By a careful choice of test functions, taking into account the thermoviscous dissipation term, we proved stability and error bounds for the linearized approximation in a suitable energy norm. Under regularity and smallness assumptions on the exact solution in an appropriate norm and smallness of the mesh size parameter $h$, we derived an error estimate in the energy norm for the nonlinear problem by employing \change{orange}{the Banach fixed-point theorem}. 

Additionally, we have conducted extensive convergence studies to support our theoretical findings for different options with respect to modeling the tissue and material parameters of the discussed model. Finally, a 3D simulation using realistic material data and geometries showing the propagation of ultrasound waves into human tissue closes the loop to the original motivation of the work.

Future research on the topic will be concerned with the analysis and simulation of problems including nonlinear effects not only in the acoustic but also elastic regime, and the simulation of ultrasound heating by an additional coupling with a temperature model. 

%\section*{Declaration of competing interest}
%The authors declare that they have no known competing financial interests or personal relationships that could have appeared to influence the work reported in this paper.

%\clearpage

%\vspace*{-2cm}

\section*{Acknowledgments}
\change{dkgreen}{We would like to thank the reviewers for their careful reading of our manuscript and for the insightful remarks.} M. Muhr and B. Wohlmuth acknowledge the financial support provided by the Deutsche Forschungsgemeinschaft under the grant number WO 671/11-1.

\version{\section*{Appendix A}
We present here the postponed proof of Theorem~\ref{Thm:LinError}. \\

\begin{proof}[Proof of Theorem \ref{Thm:LinError}]
%Arguments making use of the concrete form of the nonlinearities considered here will follow later within a fixed-point argument.\\
Insertion of a solution $(\u,\psi)$ into the semi-discrete form and subtraction of the actual semi-discrete form therefrom yields the following equation for the error $e$:
	\begin{equation*}
	\begin{aligned} 
		&\producte{\varrhoe\ddot{e}_{\u}}{\v_h}+\producte{2\varrhoe\zeta\dot{e}_{\u}}{\v_h}+\producte{\varrhoe\zeta^2 e_{\u}}{\v_h}+\formeh(e_{\u},\v_h)\\
		&\qquad +\producta{c^{-2}\ddot{e}_{\psi}}{\phi_h}+\formah(\tilde{e}_{\psi},\phi_h)+\mathcal{I}(\varrhof\dot{\tilde{e}}_{\psi},\v_h)-\mathcal{I}(\phi_h,\dot{e}_{\u})\\
		&=\producta{\fa-\fah}{\phi_h}
	\end{aligned}
	\end{equation*}
	for all $(\v_h,\phi_h)\in\Veh\times\Vah$, a.e.\ in time. We next decompose the error into the difference of $(e_{\u,I},e_{\psi,I})$ and $(e_{\u,h}, e_{\psi,h})$, see \eqref{eq:ErrorDecomp}. We  proceed similarly to the stability proof and test the resulting error equation for $(e_{\u,h}, e_{\psi,h})$ with $\v_h=\dot{e}_{\u,h}$ and $\phi_h=\varrhof\dot{\tilde{e}}_{\psi,h}$.\\
	
	This yields
	\begin{equation*}\label{error_eq}
	\begin{aligned}
	&\producte{\varrhoe\ddot{e}_{\u,h}}{\dot{e}_{\u,h}}+\producte{2\varrhoe\zeta\dot{e}_{\u,h}}{\dot{e}_{\u,h}}+\producte{\varrhoe\zeta^2 e_{\u,h}}{\dot{e}_{\u,h}}+\formeh(e_{\u,h},\dot{e}_{\u,h})\\
	&\qquad +\producta{c^{-2}\ddot{e}_{\psi,h}}{\varrhof\dot{\tilde{e}}_{\psi,h}}+\formah(\tilde{e}_{\psi,h},\varrhof\dot{\tilde{e}}_{\psi,h})+\mathcal{I}(\varrhof\dot{\tilde{e}}_{\psi,h},\dot{e}_{\u,h})-\mathcal{I}(\varrhof\dot{\tilde{e}}_{\psi,h},\dot{e}_{\u,h})\\
	&=\producte{\varrhoe\ddot{e}_{\u,I}}{\dot{e}_{\u,h}}+\producte{2\varrhoe\zeta\dot{e}_{\u,I}}{\dot{e}_{\u,h}}+\producte{\varrhoe\zeta^2 e_{\u,I}}{\dot{e}_{\u,h}}+\formeh(e_{\u,I},\dot{e}_{\u,h})\\
	&\qquad +\producta{c^{-2}\ddot{e}_{\psi,I}}{\varrhof\dot{\tilde{e}}_{\psi,h}}+\formah(\tilde{e}_{\psi,I},\varrhof\dot{\tilde{e}}_{\psi,h})+\mathcal{I}(\varrhof\dot{\tilde{e}}_{\psi,I},\dot{e}_{\u,h})-\mathcal{I}(\varrhof\dot{\tilde{e}}_{\psi,h},\dot{e}_{\u,I})\\
	&\qquad-\producta{\fa-\fah}{\varrhof\dot{\tilde{e}}_{\psi,h}}.
	\end{aligned}
	\end{equation*}
We denote the left-hand side of this equation by LHS and the right-hand side by RHS. %To further transform it, we apply the chain rule ``backwards'' to obtain the following expression being equal to the left hand side. The bilinear forms $a_h^{\textup{e}}(\cdot,\cdot)$ and  $a_h^{\textup{a}}(\cdot,\cdot)$ were resolved according to their definitions here and it was used that with the concrete choice of test-functions, 
	Noting that the interface terms on the left cancel out with our particular choice of test functions, we are left with
	\begin{align*}
	\textup{LHS}=&\,\begin{multlined}[t]\frac{1}{2}\frac{\textup{d}}{\textup{d}t}\left(\norme{\sqrt{\varrhoe}\dot{e}_{\u,h}}^2+\norme{\sqrt{\varrhoe}\zeta e_{\u,h}}^2+\norme{\sqrt{\C}\boldsymbol{\varepsilon}(e_{\u,h})}^2\right)+\norme{\sqrt{2\varrhoe\zeta}\dot{e}_{\u,h}}^{2}\\
	+\frac{1}{2}\frac{\textup{d}}{\textup{d}t}\left(\norma{c^{-1}\sqrt{\varrhof}\dot{e}_{\psi,h}}^2+\norma{\sqrt{\varrhof}\nabla\tilde{e}_{\psi,h}}^2-(\varrhof\llbrace \nabla \tilde{e}_{\psi,h}\rrbrace,\llbracket \tilde{e}_{\psi,h}\rrbracket)_{\Fhaa}\right.\\
+\left.\|{\sqrt{\chi\varrhof}\,\llbracket \tilde{e}_{\psi,h}\rrbracket}\|_{\Fhaa}^2\right)+\norma{\sqrt{b\varrhof}c^{-2}\ddot{e}_{\psi,h}}^2. \end{multlined}
	\end{align*}	
	Above we have also used the following identity:
	\begin{equation*}
		\frac{\textup{d}}{\textup{d}t}(\llbrace \nabla \tilde{e}_{\psi,h}\rrbrace,\llbracket \tilde{e}_{\psi,h}\rrbracket)_{\Fhaa}=(\llbrace \nabla\tilde{e}_{\psi,h}\rrbrace,\llbracket\dot{\tilde{e}}_{\psi,h}\rrbracket)_{\Fhaa}+(\llbrace \nabla\dot{\tilde{e}}_{\psi,h}\rrbrace,\llbracket{\tilde{e}}_{\psi,h}\rrbracket)_{\Fhaa}
	\end{equation*}
	and the fact that $\varrho \f =\textup{const}$. We next integrate LHS with respect to time and identify the energy norms $\normEe{e_{\u,h}}^2$ and $\normEf{e_{\psi,h}}^2$ (scaled by appropriate constants) within the result to obtain
	\begin{align*}
		\int_0^t\textup{LHS}(s)~\textup{d}\tau \gtrsim&\,\begin{multlined}[t] \normEe{e_{\u,h}(t)}^2+\int_0^t\norme{\sqrt{2\varrhoe\zeta}\dot{e}_{\u,h}}^{2}~\textup{d}\tau+\normEa{e_{\psi,h}{(t)}}^2\\-(\llbrace \nabla \tilde{e}_{\psi,h}{(t)}\rrbrace,\llbracket \tilde{e}_{\psi,h}{(t)}\rrbracket)_{\Fhaa}\\
	-\normEe{e_{\u,h}(0)}^2-\normEa{e_{\psi,h}(0)}^2+(\llbrace \nabla \tilde{e}_{\psi,h}(0)\rrbrace,\llbracket \tilde{e}_{\psi,h}(0)\rrbracket)_{\Fhaa}.\end{multlined}
	\end{align*}
	Note that all terms evaluated at time $t=0$ vanish due to our choice of the discrete initial conditions. Applying Lemma~\ref{Lemma:Ineq} yields the following estimate:
	\begin{equation*}\label{eq:LHSestimate}
		\normE{{e}_h{(t)}}^2=\normEe{e_{\u,h}{(t)}}^2+\normEa{e_{\psi,h}{(t)}}^2 \lesssim \int_0^t\textup{LHS}~\textup{d}\tau
	\end{equation*}
for all $t \in [0,T]$. \\
\indent We proceed with estimating the terms on the right-hand side. By the Cauchy--Schwarz and Young's inequalities, we have 
	\begin{equation*}
		\producte{\varrho\ddot{e}_{\u,I}}{\dot{e}_{\u,h}}\lesssim  \norme{\ddot{e}_{\u,I}}^2+\norme{\dot{e}_{\u,h}}^2.
	\end{equation*}
%where $C(\epsilon^{-1})$ is some constant depending on $\epsilon^{-1}$ with $\epsilon$ being chosen later. 
Noting that $\dot{\tilde{e}}_{\psi}=\dot{e}_{\psi}+\frac{b}{c^2}\ddot{e}_{\psi}$, by the triangle and Young's $\epsilon$ inequalities, we infer
	\begin{equation*} \label{intermediate_est}
		\producta{c^{-2}\ddot{e}_{\psi,I}}{\varrhof\dot{\tilde{e}}_{\psi,h}}\lesssim C(\epsilon^{-1})\norma{\ddot{e}_{\psi,I}}^2+ \norma{\dot{e}_{\psi,h}}^2+\epsilon\norma{\ddot{e}_{\psi,h}}^2,
	\end{equation*}
%where we note that for simplicity of notation at each occurrence of Young's $\epsilon$ inequality we choose the same $\epsilon$, being the minimum of all individually possible $\epsilon$, while a tracking of all different $\epsilon$ would be possible and could lead to sharper bounds. 
where $\epsilon>0$ can be chosen as conveniently small later on. To simplify exposition, we choose the same $\epsilon>0$ in all applications of Young's $\epsilon$-inequality. We similarly treat the elastic damping terms with the parameter $\zeta$:
	\begin{align*}
		\producte{2\varrhoe\zeta\dot{e}_{\u,I}}{\dot{e}_{\u,h}}&\lesssim  \norme{\dot{e}_{\u,I}}^2+ \norme{\dot{e}_{\u,h}}^2,\\
		\producte{\varrhoe\zeta^2e_{\u,I}}{\dot{e}_{\u,h}}&\lesssim \norme{e_{\u,I}}^2+\norme{\dot{e}_{\u,h}}^2.
	\end{align*}
Integrating RHS with respect to time yields
\begin{align} 
	\int_0^t\textup{RHS}~\textup{d}\tau\lesssim\int_0^t& \left(\norme{{e}_{\u,I}}^2+\norme{\dot{e}_{\u,I}}^2+\norme{\ddot{e}_{\u,I}}^2+C(\epsilon^{-1})\norma{\ddot{e}_{\psi,I}}^2\right. \label{rhs_est} \\
	&\quad\,\,\left. +\norme{\dot{e}_{\u,h}}^2+\norma{\dot{e}_{\psi,h}}^2+\epsilon\norma{\ddot{e}_{\psi,h}}^2\right)~\textup{d}\tau+\mathcal{A}+\mathcal{C}+\mathcal{F}, \notag
\end{align}
where we have introduced the short-hand notation
\begin{align*}
	\mathcal{A}&:=\int_0^t \formeh(e_{\u,I},\dot{e}_{\u,h})+\formah(\tilde{e}_{\psi,I},\varrhof\dot{\tilde e}_{\psi,h})~\textup{d}\tau,\\
	\mathcal{C}&:=\int_0^t\mathcal{I}(\varrhof\dot{\tilde{e}}_{\psi,I},\dot{e}_{\u,h})-\mathcal{I}(\varrhof\dot{\tilde{e}}_{\psi,h},\dot{e}_{\u,I})~\textup{d}\tau,\\
	\mathcal{F}&:=\int_0^t \producta{\fa-\fah}{\varrhof{\dot{\tilde e}}_{\psi,h}}~\textup{d}\tau.
\end{align*}
\noindent We continue by estimating these three terms separately.\vspace*{1mm}
	
	\begin{itemize}
	\item \textbf{Estimating the bilinear forms $\boldsymbol{\mathcal{A}}$.}
	Integrating by parts with respect to time yields
	\begin{equation*}
		\mathcal{A}=\left[\formeh(e_{\u,I},{e}_{\u,h})+\formah(\tilde{e}_{\psi,I},\varrhof{\tilde e}_{\psi,h})\right] {{\Big \vert_{0}^t}}-\int_0^t \left( \formeh(\dot{e}_{\u,I},{e}_{\u,h})+\formah(\dot{\tilde{e}}_{\psi,I},\varrhof{\tilde e}_{\psi,h})\right)~\textup{d}\tau
	\end{equation*}
	with the expression under the integral being the same as the one in the brackets except for the additional time-derivative within the first arguments of the bilinear forms. Hence, we conduct the following estimates in a generic setting.\\´
\indent The Cauchy--Schwarz inequality applied to the bulk as well as interface forms, followed by an application of Young's $\epsilon$ inequality yields
	\begin{align*}
		\formeh(\v_I,\v_h)+\formah(\tilde{\phi}_I,\varrhof\tilde{\phi}_h)&\lesssim C(\epsilon^{-1})\left(\norme{\boldsymbol{\varepsilon}(\v_I)}^2+\norma{\nabla\tilde{\phi}_{I}}^2\right.\\
		&\left.\qquad\qquad\quad+\|{\chi^{-1/2}\llbrace\nabla\tilde{\phi}_I\rrbrace}\|_{\mathcal{F}_h^{\textup{a,a}}}^2+\normFa{\tilde{\phi}_I}^2\right)\\
		&\qquad+\epsilon\left(\norme{\boldsymbol{\varepsilon}(\v_h)}^2+\norma{\nabla\tilde{\phi}_{h}}^2\right.\\
		&\left.\qquad\qquad\quad+\|{\chi^{-1/2}\llbrace\nabla\tilde{\phi}_h\rrbrace}\|_{\mathcal{F}_h^{\textup{a,a}}}^2+\normFa{\tilde{\phi}_h}^2\right)
	\end{align*}
	everywhere in time. When estimating the terms within the time integral it is sufficient to take $\epsilon=1/2$. By means of Lemma~\ref{Lemma:ChiIneq}, the term $\|{\chi^{-1/2}\llbrace\nabla\tilde{\phi}_h\rrbrace}\|^2_{\mathcal{F}_h^{\textup{a,a}}}$ can be estimated by $\beta^{-1}\norma{\nabla\tilde{\phi}_{h}}^2$ up to a constant. Therefore,
	\begin{align*}
		\mathcal{A}\lesssim&\,C(\epsilon^{-1})\left(\norme{\boldsymbol{\varepsilon}(e_{\u,I})(t)}^2+\norma{\nabla\tilde{e}_{\psi,I}(t)}^2+\|{\chi^{-1/2}\llbrace\nabla\tilde{e}_{\psi,I}(t)\rrbrace}\|_{\mathcal{F}_h^{\textup{a,a}}}^2\right.\\
		&\left.+\normFa{\tilde{e}_{\psi,I}(t)}^2\right)+\epsilon\left(\norme{\boldsymbol{\varepsilon}(e_{\u,h})(t)}^2+\norma{\nabla\tilde{e}_{\psi,h}(t)}^2+\normFa{\tilde{e}_{\psi,h}(t)}^2\right)\\
		&+\int_0^t\left(\norme{\boldsymbol{\varepsilon}(\dot{e}_{\u,I})}^2+\norma{\nabla\dot{\tilde{e}}_{\psi,I}}^2+\|{\chi^{-1/2}\llbrace\nabla\dot{\tilde{e}}_{\psi,I}\rrbrace}\|_{\mathcal{F}_h^{\textup{a,a}}}^2+\normFa{\dot{\tilde{e}}_{\psi,I}}^2\right.\\
		&\qquad\qquad\qquad+\left. \norme{\boldsymbol{\varepsilon}(e_{\u,h})}^2+\norma{\nabla\tilde{e}_{\psi,h}}^2+\normFa{\tilde{e}_{\psi,h}}^2\right)~\textup{d}\tau,
	\end{align*}
	where the terms evaluated at $t=0$ all vanish due to $e_{\u,h}(0)=\tilde{e}_{\psi,h}(0)=0$ thanks to our choice of the approximate initial conditions in \eqref{approx_IC}.\\[1mm]

	\item \change{blue}{\textbf{Estimating the interface terms $\boldsymbol{\mathcal{C}}$.} One employs the Cauchy--Schwarz inequality on the elasto-acoustic interface $\Gamma_{\textup{I}}^{\textup{e,a}}$. Together with the inverse trace estimate from Lemma~\ref{lem:InverseEstimate} on the elastic side, the interpolation estimate on faces from Lemma~\ref{lem:InterfaceEstimate} on the acoustic and our mesh assumptions, which basically guarantee comparability of the different mesh sizes $h_{\textup{e}},h_{\textup{f}}$ and $h_{\textup{t}}$ with one common parameter $h$, one obtains
	\begin{align*}
		\mathcal{I}(\varrhof\dot{\tilde{e}}_{\psi,I},\dot{e}_{\u,h})&=\productFea{\varrhof \dot{\tilde{e}}_{\psi,I}\normale}{\dot{e}_{\u,h}}\leq \|\varrhof\dot{\tilde{e}}_{\psi,\textup{I}}\|_{L^2(\Gamma_{\textup{I}}^{\textup{e,a}})}\|\dot{e}_{\u,h}\|_{L^2(\Gamma_{\textup{I}}^{\textup{e,a}})}\\
		&=\sqrt{\sum_{F\in\mathcal{F}_h^{\textup{e,a}}}\|\varrhof\dot{\tilde{e}}_{\psi,\textup{I}}\|_{L^2(F)}^2}~\sqrt{\sum_{G\in\mathcal{G}_{h}^{\textup{e,a}}}\|\dot{e}_{\u,h}\|_{L^2(G)}^2}\\
		&\lesssim h^{s-1/2}\sqrt{\sum_{\kappa\in\mathcal{T}_{h_{\textup{f}},\textup{f}}}\|\dot{\tilde{\psi}}\|_{H^s(\kappa)}^2}~ h^{-1/2}\sqrt{\sum_{\kappa\in\mathcal{T}_{h_{\textup{e}},\textup{e}}}\|\dot{e}_{\u,h}\|_{L^2(\kappa)}^2}= h^{s-1}\normsa{\dot{\tilde{\psi}}}\|\dot{e}_{\u,h}\|_{\Omegae}.
	\end{align*}
	Analogously $\mathcal{I}(\varrhof\dot{\tilde{e}}_{\psi,h},\dot{e}_{\u,I})\lesssim h^{s-1}\|\dot{\u}\|_{H^s(\Omega_{\textup{e}})}\|\dot{\tilde{e}}_{\psi,h}\|_{\Omegaa}$.}

Hence by Young's inequality
	\begin{align*}
		\mathcal{C}&\lesssim \int_0^t h^{2(s-1)}\left(\|\dot{{\u}}\|_{H^s(\Omegae)}^2+\normsa{\dot{{\psi}}}^2+\int_0^t\normsa{\ddot{{\psi}}}^2\right)+\left(\norme{\dot{e}_{\u,h}}^2+\norma{\dot{\tilde{e}}_{\psi,h}}^{\change{dkgreen}{2}}\right)~\textup{d}\tau.
	\end{align*}
	
	\item \textbf{Estimating the source term $\boldsymbol{\mathcal{F}}$.} Again by Cauchy--Schwarz and Young's inequalities, we have
	\begin{align*}
	\mathcal{F}\lesssim \|\fa-\fah\|_{L^2(0,t;\Omegaa)}^2+\int_0^t\|\dot{e}_{\psi,h}\|_{\Omegaa}^2\textup{d}\tau+\int_0^t\norma{\ddot{e}_{\psi,h}}^2~\textup{d}\tau.
	\end{align*}
\end{itemize}

We can now incorporate the estimates for $\mathcal{A},\mathcal{C}$, and $\mathcal{F}$ into \eqref{rhs_est}. Identifying components of the energy norms $\normEe{e_{\u,h}}^2$ and $\normEa{e_{\psi,h}}^2$ within $\epsilon$ terms yields

\begin{align*}
\int_0^t\textup{RHS}~\textup{d}\tau\lesssim&\,\int_0^t\bigg(\norme{{e}_{\u,I}}^2+\norme{\dot{e}_{\u,I}}^2+\norme{\ddot{e}_{\u,I}}^2+C(\epsilon^{-1})\norma{\ddot{e}_{\psi,I}}^2\\
&\quad\qquad\qquad+\norme{\boldsymbol{\varepsilon}(\dot{e}_{\u,I})}^2+\norma{\nabla\dot{\tilde{e}}_{\psi,I}}^2+\|{\chi^{-1/2}\llbrace\nabla\dot{\tilde{e}}_{\psi,I}\rrbrace}\|_{\mathcal{F}_h^{\textup{a,a}}}^2\\
&\quad\qquad\qquad+\normFa{\dot{\tilde{e}}_{\psi,I}}^2+h^{2(s-1)}\left(\|\dot{{\u}}\|_{H^s(\Omegae)}^2+\normsa{\dot{{\psi}}}^2+\normsa{\ddot{{\psi}}}^2\right)\bigg. \\
&\qquad+\normEe{e_{\u,h}}^2+\normEa{e_{\psi,h}}^2\bigg)~\textup{d}\tau\\
&+C(\epsilon^{-1})\left(\norme{\boldsymbol{\varepsilon}(e_{\u,I})(t)}^2+\norma{\nabla\tilde{e}_{\psi,I}(t)}^2+\|{\chi^{-1/2}\llbrace\nabla\tilde{e}_{\psi,I}(t)\rrbrace}\|_{\mathcal{F}_h^{\textup{a,a}}}^2\right.\\
&\qquad\qquad\quad\left.+\normFa{\tilde{e}_{\psi,I}(t)}^2\right)+\epsilon\left(\normEe{e_{\u,h}(t)}^2+\normEa{e_{\psi,h}(t)}^2\right)\\
&+\|\fa-\fah\|_{L^2(0,t;\Omegaa)}^2.
\end{align*}
Let us denote
\begin{align*}
	C_1(t,\u,\psi,\dot{\psi})&:=\|{\u}(t)\|_{H^{s}(\Omegae)}^2+\normsa{{{\psi}(t)}}^2+\normsa{\dot{{\psi}}(t)}^2, \\
	C_2(\u,\dot{\u},\ddot{\u},\dot{\psi},\ddot{\psi})&:=\|\u\|_{H^s(\Omegae)}^2+\|\dot{\u}\|_{H^s(\Omegae)}^2+\|\ddot{\u}\|_{H^s(\Omegae)}^2+\normsa{\dot{\psi}}^2+\normsa{\ddot{\psi}}^2.
\end{align*}
Then by means of the interpolation error estimates of Lemmas \ref{lem:InterpolationOperator}, \ref{lem:InterpolationEstimate}, and the above derived bounds, we infer 
\begin{align*}
\normEe{e_{\u,h}(t)}^2+\normEa{e_{\psi,h}(t)}^2\lesssim&\int_0^t\textup{LHS}~\textup{d}\tau=\int_0^t\textup{RHS}~\textup{d}\tau\\
\lesssim&\,\int_0^t \left(h^{2(s-1)}C_2(\u,\dot{\u},\ddot{\u},\dot{\psi},\ddot{\psi}) \right.
+\left. \normEe{e_{\u,h}}^2+\normEa{e_{\psi,h}}^2\right)~\textup{d}\tau\\
&\quad+h^{2(s-1)}C_1(t,\u,\psi,\dot{\psi})
+\epsilon\left(\normEe{e_{\u,h}(t)}^2+\normEa{e_{\psi,h}(t)}^2\right)\\
&\quad+\|\fa-\fah\|_{L^2(0,t;\Omegaa)}^2.
\end{align*}
If $C$ is the hidden constant above, we can choose $\epsilon \in \change{orange}{(0, \frac{1}{2C})}$ and absorb the $\epsilon$ terms on the right by the left-hand side energies. Thus
\begin{align*}
	\normEe{e_{\u,h}(t)}^2+\normEa{e_{\psi,h}(t)}^2&\lesssim h^{2(s-1)}\left(\sup_{t\in(0,T)}C_1(t,\u,\psi,\dot{\psi})+\int_0^tC_2(\u,\dot{\u},\ddot{\u},\dot{\psi},\ddot{\psi})~\textup{d}\tau\right)\\
	&\qquad+\|\fa-\fah\|_{L^2(0,t;\Omegaa)}^2+\int_0^t\normEe{e_{\u,h}}^2+\normEa{e_{\psi,h}}^2~\textup{d}\tau,
\end{align*}
which by Gronwall's inequality implies 
\begin{align*}
	\normEe{e_{\u,h}(t)}^2+\normEa{e_{\psi,h}(t)}^2&\lesssim \\
	&\hspace*{-4cm}e^{Ct}\left(h^{2(s-1)}\left(\sup_{t\in(0,T)}C_1(t,\u,\psi,\dot{\psi})+\int_0^tC_2(\u,\dot{\u},\ddot{\u},\dot{\psi},\ddot{\psi})~\textup{d}\tau\right)+\|\fa-\fah\|_{L^2(0,t;\Omegaa)}^2\right).
\end{align*}
Taking the supremum over all $t\in(0,T)$ results in the final bound \eqref{lin_error_est}.
\end{proof}}{}

\bibliography{references}{}
\bibliographystyle{apalike}

\end{document}